\documentclass[a4paper,11pt]{article}
\usepackage[margin=25mm]{geometry}
\usepackage{amsmath,amssymb,amsthm,mathtools}
\usepackage{mathrsfs}
\usepackage{booktabs,longtable,enumitem}
\usepackage{cite}
\usepackage[colorlinks=true,linkcolor=blue,citecolor=blue,urlcolor=blue]{hyperref}
\numberwithin{equation}{section}
\newtheorem{theorem}{Theorem}[section]
\newtheorem{lemma}{Lemma}[section]
\newtheorem{proposition}[lemma]{Proposition}
\newtheorem{corollary}[lemma]{Corollary}
\theoremstyle{definition}
\newtheorem{remark}[lemma]{Remark}

\newcommand{\R}{\mathbb R}
\newcommand{\dd}{\,\mathrm d}
\newcommand{\diver}{\operatorname{div}}
\newcommand{\Dt}{\mathrm D_t}
\newcommand{\Om}{\Omega}
\newcommand{\norm}[1]{\left\|#1\right\|}
\newcommand{\be}{\begin{equation}}
\newcommand{\ee}{\end{equation}}
\title{Global Strong Solutions and Large-Time Behavior for Radial
Compressible Navier--Stokes Flow with Density-Dependent Viscosities
on a Disk and in the Plane}
\author{Kexin Li\\[2pt]
\small Department of Mathematics, University of Michigan, Ann Arbor, Michigan, United States\\
\small Corresponding author\enspace \href{mailto:kexinli@umich.edu}{\texttt{kexinli@umich.edu}}}
\date{}
\begin{document}
\maketitle

\begin{abstract}
The two-dimensional barotropic compressible Navier--Stokes equations admit
a unique global strong solution for arbitrarily large compatible radial initial data
without swirl and with density bounded away from zero, both on a disk with a no-slip wall
and in the whole plane. The shear and second viscosities have a common density exponent
between one half and one, inclusive. For each exponent in this range, the bulk
viscosity coefficient may take any nonnegative value. This is the full
range allowed by nonnegative viscous dissipation in a two-dimensional
Newtonian fluid with positive shear viscosity.
Every pressure exponent greater than or equal to one is admissible.
In particular, these results apply to Marche's viscous shallow water
equations on a flat bottom in the absence of rotation, bottom friction,
and capillarity. These equations retain the positive second viscosity
derived from free-surface flow and thus fall outside the
Bresch--Desjardins relation.
The density remains uniformly bounded
above and away from zero for all time. On the disk, the density and
velocity converge exponentially in the Sobolev norm of order three to the
spatial mean of the density and to zero, respectively. In the plane, the
solution converges in the same norm to the prescribed far-field state
with positive density and zero velocity. The same conclusions hold for
classes of non-power viscosities, including laws with different growth
rates at low and high density.
\end{abstract}



\section{Introduction}\label{sec:introduction}
On a disk $\Omega=B_R(0)$ or in the plane $\Omega=\R^2$, the
two-dimensional compressible Navier--Stokes equations with
density-dependent viscosities take the form
\begin{equation}\label{eq:system}
\begin{cases}
\rho_t+\diver(\rho U)=0,\\
\rho(U_t+U\cdot\nabla U)+\nabla P(\rho)
=\diver(2\mu(\rho) DU)+ \nabla(\lambda(\rho)\diver U),
\end{cases}
\end{equation}
Here $t>0$ is time, $x\in\Omega$ is position, and $\rho$ and $U$ are
the density and velocity, respectively. We set
$P(\rho)\triangleq\kappa\rho^\gamma$ with $\kappa>0$ and $\gamma\ge1$,
and $DU\triangleq(\nabla U+\nabla U^{\mathsf T})/2$. The viscous stress is
$2\mu(\rho)DU+\lambda(\rho)(\diver U)I$, where $I$ is the identity matrix,
and $\mu$ and $\lambda$ are the shear and second viscosity coefficients,
respectively. Motivated by viscous shallow water models, in which the
effective viscosities are proportional to the water depth~\cite{Marche07,BN11},
we consider the common power laws
\be\label{eq:viscosities}
\mu(\rho)\triangleq a\rho^\alpha,\qquad
\lambda(\rho)\triangleq b\rho^\alpha,\qquad
a>0,\quad a+b\ge0,\quad \frac12\le\alpha\le1.
\ee
In two dimensions the Newtonian viscous stress decomposes as
\[
 2a\rho^\alpha\left(DU-\tfrac12(\diver U)I\right)
 +(a+b)\rho^\alpha(\diver U)I .
\]
The shear and bulk viscosities are therefore $a\rho^\alpha$ and
$(a+b)\rho^\alpha$, respectively. For positive shear viscosity ($a>0$),
the full range allowed by nonnegative viscous dissipation is
$a+b\ge0$~\cite[Eqs.~(4.44) and~(5.26)]{Jensen15}. The endpoint $b=-a$ corresponds
to zero bulk viscosity, as occurs in the kinetic theory of a dilute gas
of smooth elastic hard disks~\cite[Table~II]{MS21}.

We prescribe the initial data
\begin{equation}\label{eq:initialdata}
(\rho,U)|_{t=0}=(\rho_0,U_0),\qquad x\in\Omega,
\end{equation}
and, on the disk, impose the no-slip boundary condition
\begin{equation}\label{eq:boundary}
U|_{\partial\Omega}=0,\qquad t>0.
\end{equation}

An important example of \eqref{eq:system} is the viscous shallow
water model, in which $\rho$ denotes the water depth and $U$ the
depth-averaged horizontal velocity. The shallow water, or Saint-Venant,
equations describe flows whose depth is small compared with the
horizontal length scale. They arise from the incompressible
Navier--Stokes equations in a thin fluid layer, as described by Gerbeau and
Perthame~\cite{GP}, Marche~\cite{Marche07}, and Bresch and Noble~\cite{BN11}.
The hydrostatic pressure is $\kappa\rho^2$, where $2\kappa>0$ is the
gravitational coefficient.
For a flat bottom and in the absence of rotation, friction, and
capillarity, Marche's model~\cite{Marche07} has viscous operator
\[
2a\diver(\rho DU)+2a\nabla(\rho\diver U),
\]
and takes the form \eqref{eq:system} with
\begin{equation}\label{eq:marcheparameters}
\alpha=1,\qquad \gamma=2,\qquad b=2a.
\end{equation}
The same viscous operator appears in the shallow water system rigorously
justified by Bresch and Noble~\cite[Eq.~(6)]{BN11}, with $a$ equal to
the reciprocal Reynolds number. In the thin-layer derivation,
incompressibility relates the vertical strain to the horizontal
divergence. The vertical normal stress at the free surface thus
contributes the term $2a\nabla(\rho\diver U)$ to the depth-averaged
momentum equation.
Here we impose a fixed no-slip lateral wall on the disk.

A central difficulty in the analysis of these models is to
control the density when the viscosities depend on it.
A major advance was the discovery of the Bresch--Desjardins
(BD) entropy, which provides additional control of density
gradients beyond the mechanical energy estimate \cite{BD03,BD06,BDL}.
For the viscous stress considered here, this entropy structure
requires
\begin{equation}\label{eq:BDrelation}
\lambda(\rho)=2\bigl(\rho\mu'(\rho)-\mu(\rho)\bigr).
\end{equation}
For \eqref{eq:viscosities}, this condition is $b=2a(\alpha-1)$.
The BD entropy is used in the construction of global weak solutions
by Li and Xin (\href{https://arxiv.org/abs/1504.06826v2}{arXiv:1504.06826v2})
and in \cite{VY16,BVY22}. L\"u, Zhang, and Zhong~\cite{LZZ19}
treated the related quantum Navier--Stokes equations on the
three-dimensional torus. For strong solutions with large radial data,
the BD relation enters the estimates described below.

For the linear BD law $\mu(\rho)=a\rho$, $\lambda(\rho)=0$,
Haspot~\cite{Haspot17} obtained global strong solutions with a large
irrotational component under partial smallness assumptions.
For large radial data in the plane, Haspot's Theorem~2.4
(\href{https://arxiv.org/abs/1901.03143v1}{arXiv:1901.03143v1}) covers $\gamma=1$ and $\gamma\ge2$ under Besov and finite-energy
assumptions, with both initial velocities $u_0,w_0$ in
$L^1\cap L^\infty$. Here $w_0$ denotes the initial effective velocity.
In bounded domains, Zhang~\cite{XZweak25} constructed weak solutions
with bounded density and additional regularity away from the center.
On a fixed two-dimensional ball with no slip, Huang, Meng, and
Zhang (\href{https://arxiv.org/abs/2512.15029v2}{arXiv:2512.15029v2}, Theorem~2.3) proved global classical solvability
for arbitrary radial $H^3$ data with strictly positive density and
$\gamma\ge3/2$, together with time-independent density bounds and
convergence to equilibrium. Gu and Huang
(\href{https://arxiv.org/abs/2512.24799v1}{arXiv:2512.24799v1}, Corollary~1.1) extended these conclusions to every $\gamma>1$ in the
constant-entropy case.
Chen, Zhang, and Zhu (\href{https://arxiv.org/abs/2512.18545v2}{arXiv:2512.18545v2}, Theorem~2.4)
treated the Cauchy
problem for the linear BD law with strictly positive initial and
far-field densities. Cao, Li, and Zhu~\cite{CLZ24} studied the exterior
problem, and Chen, Zhang, and Zhu (\href{https://arxiv.org/abs/2601.06620v3}{arXiv:2601.06620v3})
treated the vacuum free boundary.
For nonlinear viscosity laws on the BD line, Zhang~\cite{XZ} proved
global strong solvability on a ball for large positive radial data,
in particular for $4/5\le\alpha<1$ and $\gamma>1$ in two dimensions.
On the disk, Huang, Meng, and Zhang proved global classical solvability
for $1/2<\alpha<1$ and $\gamma>1$, with uniform density bounds and
convergence for $0.54<\alpha<1$
(\href{https://arxiv.org/abs/2512.15029v2}{Theorems~2.1 and~2.2} of their preprint).
On the plane, Chen, Zhang, and Zhu
(\href{https://arxiv.org/abs/2605.17121v1}{arXiv:2605.17121v1}, Theorems~2.1 and~2.2) treated $1/2<\alpha<1$ and $\gamma>1$ with either vacuum or positive
density at infinity. In the positive far-field case they proved global
regularity for arbitrary radial $H^3$ perturbations of a constant state
with $\inf\rho_0>0$.
Remark~\ref{rem:bzero} compares our results with the large-time
estimates of Lei and Liang on both domains.  Lei (\href{https://arxiv.org/abs/2604.18306v3}{arXiv:2604.18306v3})
and Huang and Li (\href{https://arxiv.org/abs/2609.11304v1}{arXiv:2609.11304v1}) obtained further radial
results under the BD relation. Fan, Huang, and  Lei
(\href{https://arxiv.org/abs/2604.08946v1}{arXiv:2604.08946v1}) studied the coupled Navier--Stokes--Poisson system.

The linear law $\mu(\rho)=a\rho$ satisfies \eqref{eq:BDrelation} at
$\lambda(\rho)=0$. The positive second viscosity in
\eqref{eq:marcheparameters} increases mechanical dissipation and
introduces a term outside the BD cancellation for density derivatives.
Fabrie and Marche~\cite{FM09} proved stability of global weak solutions
for a periodic shallow water model with friction and capillarity.
Bresch, Desjardins, and
M\'etivier~\cite[Section~3.1]{BDM06} discussed the difficulty of obtaining a
density-gradient estimate with the full physical stress.
In \href{https://arxiv.org/abs/1504.06826v2}{Remark~1.3} of their preprint, Li and Xin
identified global weak existence for general initial data with this
stress as an open problem.
This motivates the study of global solutions retaining
the second viscosity in the physical model.

For a shallow water system with the stress in \eqref{eq:marcheparameters},
Bresch and Noble~\cite[Theorem~1]{BN11} recalled global strong solvability
for small perturbations of a constant state in the plane.
Li, Pan, and Zhu~\cite{LPZ17} established local well-posedness and a continuation
criterion for linear viscosities, including Marche's coefficients, in
the two-dimensional Cauchy problem with far-field vacuum and additional
assumptions on the logarithmic density gradient. Zhai and Ye~\cite{ZY20}
proved global existence with large energy for
$\mu(\rho)=\nu\rho$, $\lambda(\rho)=\nu\rho$, $\nu>0$, under smallness
conditions involving the height perturbation and the compressible part
of the velocity. In the plane a radial velocity without swirl is
irrotational, so the latter condition constrains the whole velocity.

For other shallow water models, Sundbye~\cite{Sundbye96}
treated small Dirichlet data. With the full-gradient diffusion
$\nu\diver(\rho\nabla U)$, Wang and Xu~\cite{WX} proved local
Sobolev and global small-data results, and Lian, Huang, and
Liu~\cite{LHL12} studied radial global solutions.
Qi and Wang~\cite{QW} obtained global solutions for classes of
large Cauchy data under partial smallness assumptions.
In one dimension, Constantin, Drivas, Nguyen, and
Pasqualotto~\cite{CDNP} used the BD entropy and active potentials
to prove global positive smooth solutions for large periodic data.
Related results on radial weak solutions and free-boundary problems
appear in \cite{GJX08,LZ16,XWZ20}.

For the full heat-conducting system in several dimensions, Chen, Huang,
and Zhu~\cite{CHZ24} constructed global finite-relative-energy weak
solutions of the Cauchy problem for large, possibly discontinuous,
spherically symmetric initial data. Their solutions include the origin
and allow the formation of a vacuum there.

On the same two-dimensional disk with no slip, Huang, Su, Yan, and
Yu~\cite{HSYY} studied constant shear viscosity and a second viscosity
proportional to $\rho^m$. For $\gamma>1$ they obtained global radial
strong solutions for $m>1$, and a time-independent upper bound for the
density when $\max\{1,(\gamma+2)/4\}<m\le\gamma$. For the same type of
viscosities, Huang and Meng (\href{https://arxiv.org/abs/2506.16261v1}{arXiv:2506.16261v1}, Theorem~1.1) allowed swirl
and initial vacuum and obtained a time-independent density upper bound
for $m\ge1$ and $\gamma>1$. For constant shear viscosity, Jiu and
Zhang~\cite{JZ26} also proved global strong solvability in the plane
for large radial data in weighted spaces, including positive far-field
density. With $\mu(\rho)=a\rho^\alpha$, the radial
momentum equation contains the additional term
$-2a(\rho^\alpha)_r u/r$ in \eqref{eq:radialproblem}. Controlling this
term near the center is necessary even for a velocity without swirl.
For constant viscosities, Huang, Li, and
Xin~\cite{HLX} allowed vacuum and large oscillations under a
small-energy condition.

Large-time behavior has also been studied for related compressible
flow models. In one-dimensional unbounded domains, Li and Xu~\cite{LiXu26} obtained
time-independent bounds for the specific volume and temperature and
proved asymptotic stability for large data with degenerate heat conductivity.
For constant viscosities, Fang, Zi, and
Zhang~\cite{FZZ12} proved exponential decay of finite-energy weak
solutions on bounded domains, assuming a time-independent density
upper bound. For constant shear viscosity and a power-law second
viscosity, Fan, Li, and Wang (\href{https://arxiv.org/abs/2310.15520v1}{arXiv:2310.15520v1})
studied convergence to steady
states for large data under slip boundary conditions.
Du, Guo, and Xu~\cite{DGX26} considered power-law viscosities in a
three-dimensional bounded domain with slip boundary conditions and
an external potential force. When the initial density has a sufficiently
large lower bound, they proved global strong solvability and convergence
to the steady state in $L^p$. Fan, Huang, and Ni~\cite{FHN26}
studied power-law viscosities in
the three-dimensional Cauchy problem. They obtained global classical
solutions and optimal decay for small initial perturbations in
$L^1\cap L^2$ and viscosity exponents close to one.

For the common power laws \eqref{eq:viscosities}, we prove global
strong solvability throughout the range $b\ge-a$,
including the stress in Marche's model \eqref{eq:marcheparameters}.
The result holds on the disk and in the plane for
$1/2\le\alpha\le1$ and every fixed $\gamma\ge1$. The initial density
is bounded away from zero, and the radial data without swirl may be
arbitrarily large. The density remains bounded above and away from
zero uniformly in time. Both density and velocity converge in $H^3$
to equilibrium. On the disk this convergence is exponential. On the
plane the equilibrium density is the prescribed positive far-field value.

Throughout the paper, we write $\norm{\cdot}_{L^p}$
for the $L^p(\Omega)$ norm.
On a time interval $[0,T]$, we use the notation
$C H^k\triangleq C([0,T];H^k(\Omega))$ and
$L^pH^k\triangleq L^p(0,T;H^k(\Omega))$,
with the interval specified by context.

On the disk we seek strong solutions in the class
\begin{equation}\label{eq:class}
\begin{gathered}
\rho\in C H^3,\qquad
U\in C(H^3\cap H^1_0)\cap L^2H^4,\\
U_t\in C H^1_0\cap L^2H^2,\qquad
U_{tt}\in L^2L^2.
\end{gathered}
\end{equation}
The equations hold in distributions, and the initial data
are attained as strong traces.

We restrict to
\[
\rho(t,x)=\rho(t,r),\qquad U(t,x)=u(t,r)\frac{x}{r},\qquad r\triangleq |x|.
\]
For radial velocities without swirl, $\nabla\diver U=\Delta U$.
The radial diffusion coefficient and the two weight exponents are given by
\be\label{eq:theta}
\nu(\rho)=(2a+b)\rho^{\alpha-1},\qquad \beta=\frac{b}{2a+b},\qquad
\theta=1-\frac{2a\alpha}{2a+b}=1-\alpha(1-\beta).
\ee
The radial and tangential components of the viscous stress are
$\rho^\alpha[(2a+b)u_r+bu/r]$ and $\rho^\alpha[bu_r+(2a+b)u/r]$.
The relation $\theta-\beta=(1-\alpha)(1-\beta)$ will be used below.
The momentum equation becomes
\begin{equation}\label{eq:heat}
U_t=\nu(\rho)\Delta U+F_\alpha(\rho,U),
\end{equation}
where
\begin{equation}\label{eq:F}
F_\alpha(\rho,U)\triangleq
-U\cdot\nabla U+
[\alpha\rho^{\alpha-1}\mathcal S(U)-P'(\rho)I]\nabla\log\rho,
\end{equation}
with
\begin{equation}\label{eq:stress}
\mathcal S(U)\triangleq 2aDU+b(\diver U)I.
\end{equation}
The momentum equation determines the initial velocity time derivative,
$U_t(0)=V_0\triangleq\nu(\rho_0)\Delta U_0+F_\alpha(\rho_0,U_0)$.
Compatibility with the no-slip boundary condition requires
\begin{equation}\label{eq:compat}
V_0\in H^1_0(\Omega).
\end{equation}

\begin{remark}[Boundary compatibility]
The positive $H^3$ data give $V_0\in H^1$. At the wall $u_0(R)=0$,
so the zero trace of $V_0$ in \eqref{eq:compat} reduces to
\begin{equation}\label{eq:radialcompat}
\nu(\rho_0(R))\left(
u_0''+\frac{u_0'}R+\alpha(\log\rho_0)'u_0'
\right)(R)
=P'(\rho_0(R))(\log\rho_0)'(R).
\end{equation}
Since $U_t(0)=V_0$, this is the first-order wall compatibility
required by \eqref{eq:class}. Together with $U_0|_{\partial\Omega}=0$,
it suffices for the local construction in Appendix~\ref{app:localconstruction}.
On $\R^2$, the hypotheses of Theorem~\ref{ws:theorem} already give
$V_0\in H^1(\R^2)$.
\end{remark}
We first state the result on the disk.

\begin{theorem}\label{thm:global}
Let $\Omega=B_R(0)$, where $R>0$. Let $\kappa>0$, $\gamma\ge1$, and
assume that $a,b,\alpha$ satisfy \eqref{eq:viscosities}.
Suppose that the initial data are radial without swirl and satisfy
\[
\rho_0,U_0\in H^3(\Omega),\qquad
\inf_{x\in\Omega}\rho_0(x)\ge m_0>0,\qquad U_0\in H^1_0(\Omega),
\]
together with the compatibility condition \eqref{eq:compat}.
Then problem \eqref{eq:system}--\eqref{eq:boundary} has a unique
global positive strong solution $(\rho,U)$, which is radial without swirl.
For every $0<T<\infty$, it belongs
to \eqref{eq:class} on $[0,T]$, with $\rho_t\in C([0,T];H^2)$.
There are constants $m,M,C,\lambda_*>0$, independent of time, such that
\begin{equation}\label{eq:mainbound}
0<m\le\rho(t,x)\le M<\infty,
\qquad t\ge0,\quad x\in\overline\Omega,
\end{equation}
and
\begin{equation}\label{eq:mainH3}
\begin{aligned}
&\sup_{t\ge0}\bigl(\norm{\rho(t)}_{H^3}^2+\norm{U(t)}_{H^3}^2
                              +\norm{U_t(t)}_{H^1}^2\bigr)\\
&\quad+\int_0^\infty\bigl(\norm U_{H^4}^2+\norm{U_t}_{H^2}^2
                           +\norm{U_{tt}}_{L^2}^2\bigr)\dd t\le C.
\end{aligned}
\end{equation}
Moreover, setting
\[
\bar\rho=\frac1{|\Omega|}\int_\Omega\rho_0\dd x
        =\frac2{R^2}\int_0^R\rho_0(r)r\dd r,
\]
the solution converges exponentially to $(\bar\rho,0)$:
\begin{equation}\label{eq:maindecay}
\begin{aligned}
&\norm{\rho(t)-\bar\rho}_{H^3}^2+\norm{U(t)}_{H^3}^2
 +\norm{U_t(t)}_{H^1}^2+\norm{\rho_t(t)}_{H^2}^2
 \le Ce^{-\lambda_* t},\\
&\int_t^\infty
 \bigl(\norm U_{H^4}^2+\norm{U_t}_{H^2}^2
                   +\norm{U_{tt}}_{L^2}^2\bigr)\dd s
 \le Ce^{-\lambda_* t},\qquad t\ge0.
\end{aligned}
\end{equation}
The constants depend only on $R,a,b,\kappa,\alpha,\gamma,m_0^{-1}$
and the initial $H^3$ norms. No smallness condition is imposed on
the initial norms.
\end{theorem}

On the plane the equilibrium density is prescribed at infinity. We use
the same letter $\bar\rho$ for this positive constant. The Cauchy
problem consists of \eqref{eq:system} and \eqref{eq:initialdata} with
$\Omega=\R^2$ and
\begin{equation}\label{ws:farfield}
  (\rho,U)(t,x)\longrightarrow(\bar\rho,0)\qquad (|x|\to\infty).
\end{equation}

\begin{theorem}[The Cauchy problem]\label{ws:theorem}
Let $\bar\rho,\kappa>0$, $\gamma\ge1$, and let $a,b,\alpha$ satisfy
\eqref{eq:viscosities}. Suppose that the data are radial without swirl,
$\rho_0\ge m_0>0$, and $\rho_0-\bar\rho,U_0\in H^3(\R^2)$.
Then the Cauchy problem has a unique global strong solution. On every
finite interval it satisfies \eqref{eq:class} with $\rho-\bar\rho$
in place of $\rho$ and with all spaces taken on $\R^2$. Moreover,
$\rho_t\in C([0,T];H^2(\R^2))$ for each $T<\infty$.
There exist constants $m,M,C>0$, independent of time, such that
\begin{equation}\label{ws:main-estimate}
\begin{gathered}
 m\le\rho(t,x)\le M,\qquad
 \sup_{t\ge0}\bigl(\|\rho(t)-\bar\rho\|_{H^3}
             +\|U(t)\|_{H^3}+\|U_t(t)\|_{H^1}\bigr)\le C,\\
 \sup_{t\ge0}\int_t^{t+1}
  \bigl(\|U\|_{H^4}^2+\|U_t\|_{H^2}^2+\|U_{tt}\|_{L^2}^2\bigr)\dd s
 \le C .
\end{gathered}
\end{equation}
Moreover,
\begin{equation}\label{ws:main-limit}
 \lim_{t\to\infty}
 \bigl(\|\rho(t)-\bar\rho\|_{H^3}+\|U(t)\|_{H^3}
       +\|U_t(t)\|_{H^1}+\|\rho_t(t)\|_{H^2}\bigr)=0 .
\end{equation}
\end{theorem}

\begin{corollary}[Marche's viscous shallow water model]\label{cor:shallowwater}
For the physical coefficients \eqref{eq:marcheparameters}, the water
depth $\rho$ and depth-averaged velocity $U$ satisfy
\begin{equation}\label{eq:marchemodel}
\begin{cases}
\rho_t+\diver(\rho U)=0,\\
\rho(U_t+U\cdot\nabla U)+2\kappa\rho\nabla\rho
=2a\diver(\rho DU)+2a\nabla(\rho\diver U).
\end{cases}
\end{equation}
On the no-slip disk assume the hypotheses of Theorem~\ref{thm:global},
and in the plane those of Theorem~\ref{ws:theorem}. Then
\eqref{eq:marchemodel} has a unique global strong solution for radial
data without swirl of arbitrary size. The water depth remains bounded
above and away from zero uniformly in time. On the disk the depth
converges to its spatial mean and the velocity to zero exponentially
in $H^3$ by \eqref{eq:maindecay}. In the plane the depth converges to
its prescribed far-field value and the velocity to zero in $H^3$ by
\eqref{ws:main-limit}.
\end{corollary}

\begin{remark}[The physical viscosity range and the BD relation]\label{rem:bzero}
For the power laws \eqref{eq:viscosities} with $a>0$, the bulk viscosity
is $(a+b)\rho^\alpha$. Thus nonnegative viscous dissipation in two
dimensions allows every $b\ge-a$, with zero bulk viscosity at $b=-a$.
For fixed $a$ and $\alpha$, the BD relation \eqref{eq:BDrelation}
selects the single value $b=2a(\alpha-1)$, which equals $-a$ only when
$\alpha=1/2$.

Under the BD relation, Lei and Liang
(\href{https://arxiv.org/abs/2609.37860v1}{arXiv:2609.37860v1}, Theorem~1.3)
proved global classical solvability and uniform density bounds on the
disk and in the plane for large radial data with positive density,
$1/2<\alpha\le1$, and $\gamma>1$. For some $q\in(2,4]$ and every
$2<q_0<q$, their estimates (1.36)--(1.37) give exponential decay of
$\rho-\bar\rho$ in $W^{1,q_0}$ and of $U$ in $H^2$ on the disk.
In the plane, (1.35) gives convergence to zero of $\rho-\bar\rho$ in
$L^s$ for every $2<s<\infty$, of $\nabla\rho$ in $L^2\cap L^{q_0}$,
of $\nabla U$ in $H^1$, and of $\rho_t,U_t$ in $L^2$.

To our knowledge, Theorems~\ref{thm:global} and~\ref{ws:theorem} are
the first to establish global strong well-posedness for arbitrarily large
radial data without swirl throughout the full physical coefficient range
of \eqref{eq:viscosities}, both on a no-slip disk and in the whole plane.
For each $1/2\le\alpha\le1$ and $\gamma\ge1$, they allow every $b\ge-a$.
The initial density is bounded away from zero, and on the disk the data
satisfy the wall compatibility condition \eqref{eq:compat}.
In both domains $(\rho-\bar\rho,U)$ tends to zero in the full $H^3$ norm,
with exponential convergence on the disk.
\end{remark}

The following two results use a common relation between the viscosities.
For a fixed $-1\le\theta<1$, set
\begin{equation}\label{eq:viscosity-structure}
 A(s):=2\mu(s)+\lambda(s)=\frac{2s\mu'(s)}{1-\theta},\qquad s>0 .
\end{equation}
Every power law in \eqref{eq:viscosities} satisfies \eqref{eq:viscosity-structure}
with $\theta$ given by \eqref{eq:theta}. For $\mu>0$, nonnegative bulk
viscosity is equivalent to $s\mu'/\mu\ge(1-\theta)/2$. Equality at
every density forces $\mu(s)=Cs^{(1-\theta)/2}$.

\begin{proposition}[A structural condition on the viscosities]\label{prop:general-laws}
Let $a,b,\alpha$ satisfy \eqref{eq:viscosities}, and take $\theta$ from
\eqref{eq:theta}. Assume \eqref{eq:viscosity-structure} and
\begin{equation}\label{eq:general-law-condition}
 \begin{gathered}
  A(s)=(2a+b)s^\alpha L(s),\qquad \mu(0+)=0,\qquad
  \mu(s)+\lambda(s)\ge0,\\
  0<c\le L(s)\le C,\qquad
  0\le-(1+s)\frac{L'(s)}{L(s)}\le\varepsilon_0,\qquad s>0,
 \end{gathered}
\end{equation}
where $L$ is smooth on $(0,\infty)$. For $\varepsilon_0>0$
sufficiently small depending only on $\alpha$ and $b/a$,
the conclusions of Theorems~\ref{thm:global} and~\ref{ws:theorem} hold
for these viscosities and every fixed $\gamma\ge1$, under the same
initial and far-field assumptions. On the disk the wall compatibility
condition reads
\[
 \left[\frac{A(\rho_0)}{\rho_0}\left(u_0''+\frac{u_0'}R\right)
       +\bigl(A'(\rho_0)u_0'-P'(\rho_0)\bigr)(\log\rho_0)'\right]_{r=R}=0 .
\]
The estimates may depend on the chosen viscosity functions.
\end{proposition}

\begin{remark}[Power-law and non-power viscosities]
Given $L$, the viscosities are
\[
 \mu(s)=a\alpha\int_0^s y^{\alpha-1}L(y)\dd y,\qquad
 \lambda(s)=(2a+b)s^\alpha L(s)-2\mu(s).
\]
For $L\equiv1$ this gives every power law in \eqref{eq:viscosities},
including $b=-a$. At this endpoint, integration by parts gives
\[
 \mu(s)+\lambda(s)=a\int_0^s y^\alpha L'(y)\dd y\le0.
\]
Thus $\mu+\lambda\ge0$ forces $L$ to be constant at this endpoint.
If $a+b>0$, take $L(s)=1+\delta/(1+s)$ with
$0<\delta<\min\{\varepsilon_0,(a+b)/a\}$. Indeed,
$A(s)-\mu(s)\ge(a+b-a\delta)s^\alpha>0$.

When $b>0$ or $\alpha=1/2$, the proof also allows the weaker slope
condition $|sL'(s)/L(s)|\le\varepsilon_0$, retaining the explicit
assumption $\mu+\lambda\ge0$. For example, $L(s)=1+\delta s/(1+s)$
satisfies these conditions at $\alpha=1/2$, $b=-a$ for small $\delta>0$.
\end{remark}

The next condition allows the shear viscosity to grow at different rates at
low and high density. Both results are proved in Section~\ref{sec:general-laws}.

\begin{corollary}[A viscosity class without a smallness condition]\label{cor:mixed-laws}
Let $\mu$ be smooth on $(0,\infty)$ and define $\lambda$ by
\eqref{eq:viscosity-structure} for a fixed $\theta$. Assume
\begin{equation}\label{eq:mixed-law-condition}
 \frac25\le\theta\le\frac45,\qquad
 \mu(s)\ge c\sqrt{s},\qquad
 \frac25\le\frac{s\mu'(s)}{\mu(s)}\le1\qquad(s>0)
\end{equation}
for some $c>0$. Then the conclusions of Proposition~\ref{prop:general-laws}
hold on both domains for every fixed $\gamma\ge1$, with the same wall
compatibility condition on the disk.

For example, let $\mu(s)=a\sqrt{s}+d s$ with arbitrary $a,d>0$.
At $\theta=1/2$ this gives $\lambda(s)=2d s$. At $\theta=2/5$ it gives
$\lambda(s)=(-a\sqrt{s}+4d s)/3$, which changes sign while
$\mu(s)+\lambda(s)>0$ for all $s>0$.
\end{corollary}

\begin{remark}[Viscosity and pressure ranges]
For the stress convention in \eqref{eq:system}, the two-dimensional
dissipation decomposes as
\[
2a\rho^\alpha|DU|^2+b\rho^\alpha(\diver U)^2
=2a\rho^\alpha\left|DU-\tfrac12(\diver U)I\right|^2
 +(a+b)\rho^\alpha(\diver U)^2.
\]
For $a>0$, this form is nonnegative for every strain at positive density
if and only if $b\ge-a$. At $b=-a$ only the trace-free term remains.
For every such $b$ and $1/2\le\alpha\le1$, Theorems~\ref{thm:global}
and~\ref{ws:theorem} allow every fixed $\gamma\ge1$. This includes the isothermal law $\gamma=1$, the usual
ideal-gas exponents $1<\gamma\le2$ in two-dimensional models
\cite[p.~2]{LMS22}, the values $7/5$ for air and $5/3$ for a
monatomic gas~\cite{FSW12}, and the shallow water law
$\gamma=2$~\cite[Eq.~(6)]{BN11}.
\end{remark}

\begin{remark}[Large radial data and the symmetry center]\label{rem:largedatanonbd}
Both domains contain the symmetry center, where the area weight $r$
in the pressure integral gives little control of density concentration.
For $b>0$, the weight $r^{-\theta}$ in \eqref{eq:wz} and conservation
of the mass inside a material circle yield the bound on $r^2\rho^\gamma$
in Lemma~\ref{pos:mass-pressure}. On the disk this bound reaches the wall.
On the plane it is used only inside the material circle
\eqref{ws:interface}. Outside this circle, \eqref{ws:exterior}
controls the deviation from the positive far-field density.
At $\alpha=1/2$, the identity $\rho(w-u)^2=(2a+b)^2|(\log\rho)_r|^2$
provides the logarithmic density estimate at the center in
Section~\ref{pr:half-all-pressure}. These estimates allow arbitrarily
large initial energy and density variation in both domains.
\end{remark}

\begin{remark}[Equilibrium and large-time behavior]
On the disk, conservation of mass fixes the equilibrium density at the
initial mean. With the uniform density upper bound, the pressure test
in Lemma~\ref{ext:early-decay} gives exponential decay of the relative
mechanical energy. Damping of $\nabla\log\rho$ and the parabolic estimates
in Section~\ref{ext:higher-section} give exponential convergence in $H^3$
as in \eqref{eq:maindecay}. On the plane, the equilibrium density is the
prescribed far-field value $\bar\rho$. The low-frequency estimate
\eqref{ws:zero-decay} gives convergence in $L^2$, and the higher-order
estimates give convergence in $H^3$ by \eqref{ws:main-limit}.
This conclusion requires no $L^1$ assumption on the initial perturbation.
The whole-space theorem gives no quantitative decay rate under these hypotheses.
\end{remark}

\begin{remark}[Open problems and further directions]\label{rem:generalviscosities}
For the viscosities \eqref{eq:viscosities}, extending the large-data
result to radial flows with swirl or to general nonradial flows is a
natural question. Initial vacuum makes both viscosities vanish and
the logarithmic derivative in \eqref{eq:w} singular. In shallow
water this corresponds to an initially dry region. For the non-BD
Cauchy problem, allowing vacuum at infinity also removes the positive
far-field state used in the exterior estimate \eqref{ws:exterior}.

Another question is whether the bounds persist for different powers,
for instance $\mu(\rho)=a\rho^\alpha$ and $\lambda(\rho)=b\rho^\delta$
with $a,b>0$ and $\delta>0$, $\delta\ne\alpha$. For these laws
$2\rho\mu'/(2\mu+\lambda)$ varies with the density, so the fixed
radial weight used in Proposition~\ref{prop:general-laws} no longer
cancels the geometric term.
\end{remark}

We now describe the main ideas of the proof. Under radial symmetry,
problem \eqref{eq:system}--\eqref{eq:boundary} reduces to
\begin{equation}\label{eq:radialproblem}
\begin{cases}
\displaystyle \rho_t+(\rho u)_r+\frac{\rho u}{r}=0,
&0<r<R,\ t>0,\\[2mm]
\displaystyle \rho(u_t+uu_r)+(\kappa\rho^\gamma)_r
=(2a+b)\bigl[\rho^\alpha(u_r+u/r)\bigr]_r
  -\frac{2a}{r}(\rho^\alpha)_r u,
&0<r<R,\ t>0,\\[2mm]
u(t,0)=u(t,R)=0,&t\ge0,\\[1mm]
(\rho,u)(0,r)=(\rho_0,u_0)(r),&0\le r\le R.
\end{cases}
\end{equation}
Radial $H^3$ regularity gives $\rho_r(t,0)=0$ and $u(t,r)=O(r)$
as $r\downarrow0$. The condition $u(t,R)=0$ expresses the no-slip boundary
condition. Throughout the estimates, $[0,\tau]\subset[0,T]$ denotes a
positive strong-solution interval. A constant $C_T$ may depend on the
prescribed finite time $T$, the initial data, and the fixed parameters.
The constants in the a priori estimates are independent of $\tau$.

Define the effective velocity
\begin{equation}\label{eq:w}
w\triangleq u+(2a+b)\rho^{\alpha-1}(\log\rho)_r .
\end{equation}
Write $\Dt\triangleq\partial_t+u\partial_r$.
The continuity equation and \eqref{eq:radialproblem}$_2$ give
\begin{equation}\label{eq:weq}
\Dt w=\left(\theta\frac{u}{r}-K\right)(w-u),
\qquad K\triangleq \frac{\kappa\gamma}{2a+b}\rho^{\gamma-\alpha}.
\end{equation}
Thus $r^{-\theta}$ cancels the term linear in $w$ that contains $u/r$:
\be\label{eq:wz}
\Dt(r^{-\theta}w)+Kr^{-\theta}w
=Kr^{-\theta}u-\theta r^{-\theta-1}u^2 .
\ee
In the linear case $\theta=\beta$. For the power laws
\eqref{eq:viscosities}, the sign of $\theta$ determines the sign of the
last term in \eqref{eq:wz}. When $\theta>0$, completing the
square gives
\[
Kr^{-\theta}u-\theta r^{-\theta-1}u^2
\le \frac{K}{4\theta}\bigl(Kr^{1-\theta}\bigr).
\]
A bound on $Kr^{1-\theta}$ therefore yields an upper comparison for
$r^{-\theta}w$.

Testing radial momentum with the weight $r^{1-\theta}$ and using
the wall flux bounds $r^{1-\theta}\rho^\alpha$ uniformly in time.
For $b>0$ and $\alpha>1/2$, the positive part of $r^{-\theta}w$ is
tested with a power of the enclosed mass, which is constant along
particle paths. This gives a time-independent bound on $r^2\rho^\gamma$.
In the inward velocity moment of order $5+b/a$, this bound absorbs
the large-density pressure contribution into the positive transport
term. A weighted quadratic energy then controls the spatial modulus
of the density, and the velocity comparisons give finite-time
positivity down to the center.

For $b\le0$, the sign of $\theta$ distinguishes three estimates.
When $\theta<0$, the favorable sign in \eqref{eq:wz} is used in
a lower comparison. When $\theta=0$, the effective velocity satisfies
a relaxation equation and yields the BD estimate. For $\theta>0$
and $\alpha>1/2$, a weighted quadratic estimate contains an inward
cubic term. On $\{\rho\ge1\}$ this term is controlled by the mechanical
dissipation and the weighted quadratic energy. On $\{\rho<1\}$ it is
absorbed by the weighted dissipation and the transport term of an
inward velocity moment. The half-power endpoint uses the corresponding
weighted control of $(\log\rho)_r$.

Once a positive density band on each finite time interval and an
unweighted density-gradient estimate are available, maximal regularity
for \eqref{eq:heat} gives continuation. On the disk, the pressure test
gives exponential decay of the mechanical energy. At $b=-a$ we apply
this test after establishing uniform density bounds and recovering
control of the full velocity gradient. The damped density-gradient
equation then gives \eqref{eq:maindecay} together with parabolic regularity.

Guo, Xu, and Zhang~\cite{GXZ25} matched interior and exterior solutions
of a three-dimensional BD system across a stress-free material sphere.
In the present Cauchy problem the positive far-field density makes
the total mass infinite. Relative energy controls the mass inside
each circle, so a sufficiently distant material circle stays in a fixed
annulus. Pressure damping keeps the density there close to its far-field
value. Averaging over the mass labels selects a material circle with
$\int_0^\tau |u_r(t,R(t))|^2\dd t\le C$ as in \eqref{ws:interface},
including at $b=-a$. Once the core estimates give a uniform spatial modulus,
mechanical dissipation forces
every local large-time limit to have constant pressure. The far field
identifies the constant. A low-frequency estimate then gives convergence
in $L^2(\R^2)$ and the higher-order estimates give \eqref{ws:main-limit}.

Section~\ref{sec:preliminary} gives the energy and flux identities.
The next three sections treat $b>0$, $b\le0$, and $\alpha=1/2$ on the disk.
The density bounds for $\alpha=1$, $b=-a$ are completed for both domains
in Subsection~\ref{ws:finite-high}.
Sections~\ref{sec:higher-continuation} and~\ref{ext:higher-section}
prove global continuation and exponential convergence.
Section~\ref{ws:section} adapts these estimates to the plane and treats
the low frequencies. Section~\ref{sec:general-laws} proves the extensions
to non-power viscosities. The local theory is given in
Appendix~\ref{app:localconstruction}.

\section{Preliminaries}\label{sec:preliminary}
We work on the disk until Section~\ref{ws:section}, except where an estimate is proved for both domains.

\subsection{Local strong solutions}\label{sec:local}
\begin{proposition}\label{prop:local}
Let $R,a,\kappa>0$, $b\ge-a$, $1/2\le\alpha\le1$, and $\gamma\ge1$.
Suppose that the initial data are radial and without swirl, that
$(\rho_0,U_0)\in H^3(\Omega)$, and that $\rho_0\ge m_0>0$.
Assume also that $U_0\in H^1_0(\Omega)$ and that \eqref{eq:compat} holds.
Then \eqref{eq:system}--\eqref{eq:boundary} admits a local radial solution
of class \eqref{eq:class}, with $\rho\ge m_0/2$ and $\rho_t\in CH^2$. Its lifespan depends only
on the parameters, $m_0^{-1}$, and the initial $H^3$ norms.
The solution is unique in the positive strong class.
\end{proposition}
In Appendix~\ref{app:localconstruction} we construct the solution by
iteration. Each transport and parabolic step preserves the prescribed
initial data, and the velocity has initial time derivative $V_0$ at every step.

\subsection{Radial weights at the origin}
Let $[0,T_*)$ be the maximal existence interval of the solution from
Proposition~\ref{prop:local}. Fix $\tau<T_*$. Unless otherwise specified,
$C$ depends only on $R$, the fixed coefficients, and the initial data,
and is independent of time. The constant $C_T$ may also depend on a
fixed $T\ge\tau$. All estimates before continuation are understood on
positive strong-solution intervals $[0,\tau]$.
On $[0,\tau]$, positivity, Sobolev embedding, and radial symmetry give,
for $\max\{0,\theta\}<\chi<1$,
\begin{equation}\label{eq:centerregularity}
\begin{gathered}
 \log\rho\in CH^3,\qquad U,\log\rho\in C C^{1,\chi},\\
 U(t,0)=0,\qquad \nabla\log\rho(t,0)=0,\\
 u=O_\tau(r),\qquad (\log\rho)_r=O_\tau(r^\chi),\qquad
 r^{-\theta}w=O_\tau(r^{\chi-\theta}) .
\end{gathered}
\end{equation}
For the nonlinear velocity tests with power $\ell>2$ and radial
weight $r^{1-\sigma}$, $\sigma<\ell$, the stress and pressure
boundary terms at $\epsilon$ are $O_\tau(\epsilon^{\ell-\sigma})$.
The transport boundary term is $O_\tau(\epsilon^{\ell+2-\sigma})$.
For the $p$-moment of $r^{-\theta}w$, the boundary flux is
$O_\tau(r^{2+p(\chi-\theta)})$. Its geometric source has integrand
$O_\tau(r^{2-\theta+(p-1)(\chi-\theta)})$.
The quadratic weighted test is justified by the same bounds with
$\ell=2$ and $\sigma=2\theta<2$.
Thus integration over $[\epsilon,R]$ followed by $\epsilon\downarrow0$
gives the stated identities on the disk.

The pointwise velocity comparisons below may likewise be made first
on an annulus. For each $\epsilon>0$ the radial profiles on
$[\epsilon,R]$ belong to $C([0,\tau];H^3([\epsilon,R]))$, and
$u_t$ belongs to $C([0,\tau];H^1([\epsilon,R]))$. The one-dimensional
Sobolev embeddings give continuous $u_{rr}$ and $u_t$ there.
At the inner endpoint we use \eqref{eq:centerregularity} before
letting $\epsilon$ tend to zero.

At $t=0$, \eqref{eq:centerregularity} bounds $r^{-\theta}w_0$ in
$L^\infty$ and hence in every finite $L^p(\rho_0 r\dd r)$ on the disk.
It also gives
\[
 \int_0^R\rho_0|u_0|^\ell r^{1-\sigma}\dd r
 \le C_{\ell,\sigma}\int_0^R r^{\ell+1-\sigma}\dd r<\infty
 \qquad(\ell>\sigma\ge0).
\]
In the plane the same bounds hold on each fixed ball. The exterior
energy in \eqref{ws:exterior} is finite since $U_0,\nabla\rho_0\in L^2$
and the initial density is bounded above and away from zero.

\subsection{Energy and effective velocity}
Differentiating $\Dt\log\rho=-(u_r+u/r)$ in $r$ and using
$[\Dt,\partial_r]=-u_r\partial_r$, we obtain from \eqref{eq:w}
\[
\Dt(w-u)
=-\nu(\rho)\partial_r(u_r+u/r)
-(w-u)[\alpha u_r+(\alpha-1)u/r].
\]
Adding the momentum equation \eqref{eq:radialproblem}$_2$ divided by $\rho$
yields \eqref{eq:weq}. Multiplying by $r^{-\theta}$ gives \eqref{eq:wz}.
On each outer annulus $[\epsilon,R]$, \eqref{eq:w} and $\rho_t\in C H^2$
give $w\in C H^2$ and $w_t\in C H^1$ in the radial variable.
Thus $w_r$ and $w_t$ are continuous up to the wall, and $u w_r\to0$
as $r\uparrow R$ since $u(t,R)=0$. Taking the interior trace in
\eqref{eq:weq} gives
\[
 \frac{\mathrm d}{\mathrm dt}w(t,R)=-K(t,R)w(t,R),
\]
with one-sided time derivatives at $t=0,\tau$.

By the continuity equation,
\begin{equation}\label{lju1}
 \int_0^R\rho(t,r)r\dd r=\int_0^R\rho_0(r)r\dd r .
\end{equation}
Use the pressure potential
\[
 \Pi(s)=
 \begin{cases}
 \kappa s^\gamma/(\gamma-1),&\gamma>1,\\
 \kappa s\log(s/\bar\rho),&\gamma=1 .
 \end{cases}
\]
Then $s\Pi'(s)-\Pi(s)=P(s)$ and the relative internal energy density is
$h_\gamma(s)=\Pi(s)-\Pi(\bar\rho)-\Pi'(\bar\rho)(s-\bar\rho)$.
In integrals we write $\Pi$ for $\Pi(\rho)$.
Define the relative mechanical energy
$\mathcal E=\int_0^R(\rho u^2/2+h_\gamma(\rho))r\dd r$ and the
mechanical dissipation
\[
 \mathscr D=\int_0^R\rho^\alpha
  \{2a(u_r^2+u^2/r^2)+b(u_r+u/r)^2\}r\dd r .
\]
Testing \eqref{eq:radialproblem}$_2$ with $ur$ and using continuity gives
$\mathcal E'=-\mathscr D$. Mass conservation and this identity give
\begin{equation}\label{eq:energy}
\sup_{t\le\tau}\int_0^R(\rho u^2+\rho^\gamma)r\dd r
+\int_0^\tau\!\int_0^R\rho^\alpha
\bigl[2a(u_r^2+u^2/r^2)+b(u_r+u/r)^2\bigr]r\dd r\dd t\le C.
\end{equation}
The dissipation integrand is
$\rho^\alpha[a(u_r-u/r)^2+(a+b)(u_r+u/r)^2]$. For every $b\ge-a$,
\[
 \mathscr D\ge
 2\min\{a,a+b\}\int_0^R\rho^\alpha(u_r^2+u^2/r^2)r\dd r.
\]
At $b=-a$ only the trace-free square remains. Once the density has
a positive lower bound, this square still controls $\|\nabla U\|_{L^2}^2$,
by the following identity for $V\in H^1_0(\Omega)$ on the disk
or $V\in H^1(\R^2)$ in the plane.
\begin{equation}\label{enp:deviatoric-korn}
 \int_\Omega 2|DV-\tfrac12(\diver V)I|^2\dd x
 =\int_\Omega|\nabla V|^2\dd x .
\end{equation}
Indeed $\int\partial_iV_j\,\partial_jV_i=\int(\diver V)^2$
by integration by parts and approximation. The positive lower bound
of the viscosity is taken out of the nonnegative trace-free square
before applying the identity. An upper density bound then also
gives the corresponding weighted full-gradient estimate.
For a radial velocity, $u(t,0)=0$. Integrating $(u^2)_r=2uu_r$ gives,
on both domains,
\begin{equation}\label{eq:radial-sup}
 |u(t,r)|^2\le\int_0^r (u_y^2+u^2/y^2)y\dd y .
\end{equation}
Thus a positive density lower bound implies
$\|U(t)\|_\infty^2\le C\mathscr D(t)$, including at $b=-a$.
\subsection{Uniform weighted density and pressure estimates}
To control pressure in the moment estimates, put
\[
\mathscr F=(2a+b)\rho^\alpha(u_r+u/r)-P,\qquad
H=\frac{2a+b}{\alpha} r^{1-\theta}\rho^\alpha
-\int_r^R\rho u\xi^{1-\theta}\dd\xi
+2a\int_r^R\rho^\alpha\xi^{-\theta}\dd\xi.
\]
Since $(2a+b)(1-\theta)/\alpha=2a$, the terms in
$r^{-\theta}\rho^\alpha$ cancel on differentiation, giving
$H_r=\rho w r^{1-\theta}$. Continuity and \eqref{eq:weq} then give
\[
 \partial_r(\Dt H)=\Dt H_r+u_rH_r
 =-r^{1-\theta}P_r-\theta\rho u^2r^{-\theta}.
\]
At the wall $u=0$, so $\Dt H(t,R)=-R^{1-\theta}[(2a+b)\rho^\alpha u_r]_{r=R}$,
where $u_r(t,R)$ is the trace from inside the disk.
Integrating from $r$ to $R$ and integrating the pressure term by parts gives
\begin{equation}\label{eq:H}
\Dt H+r^{1-\theta}P+(1-\theta)\int_r^RP\xi^{-\theta}\dd\xi
=\theta\int_r^R\rho u^2\xi^{-\theta}\dd\xi-R^{1-\theta}\mathscr F(R).
\end{equation}
Testing the momentum equation with the vector field $x$ gives the wall moment identity
\begin{equation}\label{eq:wall}
\begin{aligned}
 W(t)&=\int_0^R[\rho ur-2(a+b)B_\alpha(\rho)]r\dd r,\\
 W'(t)&=R^2\mathscr F(R)+\int_0^R(\rho u^2+2P)r\dd r,
\end{aligned}
\end{equation}
where $B_\alpha(s)=(s^\alpha-s)/(\alpha-1)$ if $\alpha\ne1$ and
$B_1(s)=s\log s$. Indeed, the viscous stress contributes the volume
term $-2(a+b)\int_0^R\rho^\alpha(u_r+u/r)r\dd r$. Since
$sB_\alpha'(s)-B_\alpha(s)=s^\alpha$, this term equals
$2(a+b)\frac{\mathrm d}{\mathrm dt}\int_0^R B_\alpha(\rho)r\dd r$.
The time derivative in \eqref{eq:wall} is justified by
$\rho_t=-\diver(\rho U)\in C H^2$, $U_t\in C H^1_0$ and positivity
on the fixed strong-solution interval.
Set $Y=H+R^{-1-\theta}W$. The terms involving $\mathscr F(R)$ in
\eqref{eq:H} and $R^{-1-\theta}W'$ cancel. Write
\[
 \Pi_\theta(t,r)=r^{1-\theta}P+(1-\theta)\int_r^R P\xi^{-\theta}\dd\xi .
\]
Adding $R^{-1-\theta}W'$ to \eqref{eq:H} gives
\[
 \Dt Y=G(t,r),\qquad
 G(t,r)=\theta\int_r^R\rho u^2\xi^{-\theta}\dd\xi
 +R^{-1-\theta}\int_0^R(\rho u^2+2P)\xi\dd\xi-\Pi_\theta(t,r).
\]
In \eqref{eq:wall}, Cauchy--Schwarz bounds the momentum term by
$R(\int_0^R\rho r\dd r)^{1/2}(\int_0^R\rho u^2r\dd r)^{1/2}$.
For $\alpha<1$ we have $|B_\alpha(s)|\le C_\alpha(1+s)$, and for
$\alpha=1<\gamma$ we have $|s\log s|\le C_\gamma(1+s^\gamma)$.
When $\alpha=\gamma=1$, use
\[
 \rho\log\rho=\kappa^{-1}h_1(\rho)+(1+\log\bar\rho)\rho-\bar\rho .
\]
The relative mechanical energy identity $\mathcal E'=-\mathscr D$
preceding \eqref{eq:energy} gives $\int_0^R h_1(\rho)r\dd r\le\mathcal E(0)$.
Together with mass conservation \eqref{lju1} and \eqref{eq:energy},
these estimates give $|Y-H|=R^{-1-\theta}|W|\le C$.

\begin{lemma}[Uniform weighted bounds]\label{ext:uniform-weighted}
Let $a>0$, $b\ge-a$, $1/2\le\alpha\le1$, and $\gamma\ge1$. On any
positive strong-solution interval $[0,\tau]$,
\begin{equation}\label{ext:weighted-bound}
\begin{gathered}
 \sup_{0\le t\le\tau}\left(\sup_{0<r\le R}r^{1-\theta}\rho^\alpha
          +\int_0^R\rho^\alpha r^{-\theta}\dd r\right)\le C,\\
 \int_s^t\!\int_0^R\rho^\gamma r^{-\theta}\dd r\dd\eta
 \le C(1+t-s),\qquad 0\le s\le t\le\tau .
\end{gathered}
\end{equation}
In particular, $\rho(t,r)\le Cr^{-(1-\beta)}$. The constant is
independent of $\tau$.
\end{lemma}
\begin{proof}
Write $H=\frac{2a+b}{\alpha}V-I+2aJ$, where
\[
 V=r^{1-\theta}\rho^\alpha,\qquad
 I=\int_r^R\rho u\xi^{1-\theta}\dd\xi,\qquad
 J=\int_r^R\rho^\alpha\xi^{-\theta}\dd\xi .
\]
Let $S(t)=1+C_0+(\max_rY(t,r))_+$, where $C_0$ bounds $|H-Y|$
by \eqref{eq:wall} and the mechanical energy, and write
$\theta_+=\max\{\theta,0\}$. Since $\theta_+\le(a+b)/(2a+b)$, the dissipation decomposition gives
$\theta_+\int_0^R\rho^\alpha(u_r^2+u^2/r^2)r\dd r\le C\mathscr D$
for every $b\ge-a$. Also
$\rho^{1-\alpha}r^{1-\theta}\le
C\|V\|_\infty^{(1-\alpha)/\alpha}$, since the remaining power
of $r$ is $2a(2\alpha-1)/(2a+b)\ge0$.
The Cauchy--Schwarz inequality and \eqref{eq:energy} therefore give
\begin{equation}\label{eq:absorption}
 \begin{aligned}
 |I|^2&\le C\|V\|_\infty^{(1-\alpha)/\alpha}J,\\
 \theta_+\int_r^R\rho u^2\xi^{-\theta}\dd\xi
 &\le C\|V\|_\infty^{(1-\alpha)/\alpha}\mathscr D .
 \end{aligned}
\end{equation}
For $\alpha>1/2$, use $H\le S$ and
$|I|\le aJ+C\|V\|_\infty^{(1-\alpha)/\alpha}$.
Taking the supremum and absorbing the power $(1-\alpha)/\alpha<1$ gives
\begin{equation}\label{ext:Y-coercivity}
 \|V\|_\infty+J(t,0)\le CS(t).
\end{equation}

At $\alpha=1/2$ this power equals one. To prove \eqref{ext:Y-coercivity}
at the endpoint, set $x=\log(R/r)$ and $f(x)=r\sqrt\rho\,u$. In this
coordinate,
\begin{equation}\label{eq:endpointvolterra}
 H=2(2a+b)V+\int_0^x V(y)[2a-f(y)]\dd y,\qquad
 \int_0^\infty f^2\dd x=\int_0^R\rho u^2r\dd r .
\end{equation}
Since $[f-2a]_+\le f^2/(8a)$, Gronwall's inequality in $x$ gives
\begin{equation}\label{eq:endpointcoercivity}
 \|V\|_\infty\le\frac{S(t)}{2(2a+b)}
 \exp\!\left(\frac{1}{16a(2a+b)}\int_0^R\rho u^2r\dd r\right)\le CS(t).
\end{equation}
The first bound in \eqref{eq:absorption} now gives
$2aJ\le S+|I|\le S+aJ+C\|V\|_\infty$, proving \eqref{ext:Y-coercivity}.

It remains to bound $S$. Write $m(t)=\max_{0\le r\le R}Y(t,r)$.
At a point where $Y=m(t)$,
\eqref{eq:absorption} and \eqref{ext:Y-coercivity} give
$|I|\le S/2+CJ$, hence $S\le C(1+V+J)$. Since $\gamma/\alpha\ge1$,
H\"older's inequality for the finite measure $\xi^{-\theta}\dd\xi$ gives
\[
 \Pi_\theta(t,r)\ge c(V+J)^{\gamma/\alpha}\ge cS-C .
\]
Let $X(t,q)$ be the trajectory starting at $q\in[0,R]$. Its flow fixes
both endpoints. Since $\theta<1$, the tails defining $Y$ and $G$ are
continuous up to the center, and the identity
\[
 Y(t,X(t,q))-Y(s,X(s,q))=\int_s^t G(\eta,X(\eta,q))\dd\eta
\]
extends to $q=0,R$. Thus $m$ is Lipschitz. Taking a maximizing label
at $t+\delta$ and passing to a convergent subsequence as $\delta\downarrow0$
gives $D^+m(t)\le\max_{\{Y=m(t)\}}G(t,\cdot)$.
At each such maximum, \eqref{eq:absorption} bounds the positive part
of the kinetic tail. The pressure bound then gives
$G\le C-cS+C\mathscr D S$. Taking the positive part and
enlarging $C$ to cover $m\le0$ gives
\[
 D^+S+cS\le C+C\mathscr D S .
\]
In the variation-of-constants formula the factor from $s$ to $t$ is at most
$\exp(C\mathcal E(0)-c(t-s))$, since
$\int_0^\tau\mathscr D\dd t\le\mathcal E(0)$. Hence $\sup_{[0,\tau]}S\le C$.

Finally $|Y(t,0)|\le C$ by \eqref{eq:absorption} and
\eqref{ext:Y-coercivity}. In the trajectory identity at $q=0$,
the bulk term is bounded by \eqref{eq:energy}, and the positive
part of the kinetic tail is bounded by $C\mathscr D$.
Its negative part only improves the upper pressure estimate.
Since $\Pi_\theta(t,0)=(1-\theta)\int_0^R P\xi^{-\theta}\dd\xi$,
integration over $[s,t]$ proves the pressure bound.
The pointwise density estimate
follows from $(1-\theta)/\alpha=1-\beta$.
\end{proof}

\subsection{Decay of the relative mechanical energy}

For $b\ge-a$ the pressure test below uses
\[
 \int_0^R\rho^\alpha(u_r^2+u^2/r^2)r\dd r\le C\mathscr D .
\]
For positive bulk viscosity this follows from the dissipation decomposition.
At $b=-a$ we apply the test after the independent density bounds in
\eqref{bdlt:uniform-density}, \eqref{enp:subcritical-density} and the argument
at \eqref{enp:center-power}. These bounds use the mechanical energy
identity and give the displayed inequality by \eqref{enp:deviatoric-korn}.

To estimate the density fluctuations, use a bounded function $G$ and
a radial test field with mean-zero divergence. Set
\[
 G(s)=\frac{s}{s+\bar\rho},\qquad
 v(t,r)=\frac1r\int_0^r[G(\rho(t,\xi))-\langle G(\rho(t))\rangle]\xi\dd\xi,
 \qquad \langle f\rangle=\frac2{R^2}\int_0^R f r\dd r .
\]
Then $v(0)=v(R)=0$ and $v_r+v/r=G(\rho)-\langle G(\rho)\rangle$.
The pressure term obtained by testing radial momentum with $vr$ is
\[
 \mathcal C(t)=\int_0^R P(\rho)[G(\rho)-\langle G(\rho)\rangle]r\dd r
 \ge\int_0^R[P(\rho)-P(\bar\rho)][G(\rho)-G(\bar\rho)]r\dd r .
\]
Indeed, $\langle\rho\rangle=\bar\rho$, and Jensen's inequality gives
$\langle P(\rho)\rangle\ge P(\bar\rho)$ and
$\langle G(\rho)\rangle\le G(\bar\rho)$.

The field $v$ satisfies the estimates
\begin{equation}\label{ext:bounded-test-stress}
 \int_0^R\rho^\alpha(v_r^2+v^2/r^2)r\dd r\le C\mathcal C,
 \qquad \|v\|_{L^\infty}^2\le C\int_0^R h_\gamma(\rho)r\dd r .
\end{equation}
Since $v(0)=v(R)=0$, expanding the square gives
\[
 \int_0^R(v_r^2+v^2/r^2)r\dd r
 =\int_0^R(v_r+v/r)^2r\dd r
 \le\int_0^R|G(\rho)-G(\bar\rho)|^2r\dd r\le C\mathcal C .
\]
Here the mean minimizes the squared distance to a constant, and
$|G(s)-G(\bar\rho)|^2\le C[P(s)-P(\bar\rho)][G(s)-G(\bar\rho)]$.
This proves the weighted estimate on $\{\rho\le2\bar\rho\}$.
On its complement, $|v_r|+|v/r|\le C$ and
$\rho^\alpha\le C[P(\rho)-P(\bar\rho)][G(\rho)-G(\bar\rho)]$.
For the supremum bound apply Cauchy--Schwarz in the definition of $v$,
using $|G(s)-G(\bar\rho)|^2\le C h_\gamma(s)$.

The coarse estimate in Lemma~\ref{ext:uniform-weighted} gives
$\rho^{1-\alpha}r^2\le C$, hence
\begin{equation}\label{pr:iso-kinetic-dissipation}
 \int_0^R\rho u^2r\dd r\le C\int_0^R\rho^\alpha u^2/r\dd r
 \le C\mathscr D .
\end{equation}
Since $G(s)-sG'(s)=G(s)^2$, the continuity equation gives
\[
 v_t+G(\rho)u=\frac1r\int_0^rG(\rho)^2(u_\xi+u/\xi)\xi\dd\xi
             -\frac r2\langle G(\rho)^2(u_r+u/r)\rangle .
\]
The right-hand side is bounded by $C\mathscr D^{1/2}$, by
Cauchy--Schwarz and $G(s)^4\le Cs^\alpha$. In the momentum test,
\eqref{pr:iso-kinetic-dissipation} thus bounds the convection and
$v_t$ terms by $C\mathscr D$. The stress term is at most
$C(\mathscr D\mathcal C)^{1/2}$ by \eqref{ext:bounded-test-stress}.
The endpoint terms vanish since $v(0)=v(R)=0$. We obtain
\begin{equation}\label{pr:iso-cross-energy}
 \frac{\dd}{\dd t}\int_0^R\rho uvr\dd r\ge\tfrac12\mathcal C-C\mathscr D .
\end{equation}
The integral on the left is bounded in absolute value by $C\mathcal E$,
by mass conservation and \eqref{ext:bounded-test-stress}. For small
fixed $\varepsilon>0$, the mechanical identity $\mathcal E'=-\mathscr D$
now gives
\begin{equation}\label{pr:iso-modified-energy}
 \Phi=\mathcal E-\varepsilon\int_0^R\rho uvr\dd r,\qquad
 \tfrac12\mathcal E\le\Phi\le\tfrac32\mathcal E,\qquad
 \Phi'+c(\mathscr D+\mathcal C)\le0 .
\end{equation}

\begin{lemma}[Decay of the relative mechanical energy]\label{ext:early-decay}
Let $a>0$, $b\ge-a$, $1/2\le\alpha\le1$, and $\gamma\ge1$.
At $b=-a$, assume the independent time-uniform density bounds cited
above. On a positive strong-solution interval $[0,\tau]$, if $\gamma>1$
or if the density has a time-independent upper bound, then
\begin{equation}\label{ext:energy-decay}
 \int_t^\tau\mathscr D(s)\dd s
 =\mathcal E(t)-\mathcal E(\tau)
 \le\mathcal E(t)\le Ce^{-\lambda t},\qquad 0\le t\le\tau .
\end{equation}
In general, when $\gamma=1$,
\begin{equation}\label{pr:iso-energy}
 \mathcal E(t)+\int_t^\tau\mathscr D(s)\dd s
 \le C e^{-c\sqrt{1+t}},\qquad 0\le t\le\tau .
\end{equation}
For every $\gamma\ge1$ and $0<s\le1$,
\begin{equation}\label{pr:iso-fractional-coefficients}
 \int_0^\tau \mathscr D(t)^s\dd t\le C_s .
\end{equation}
The constants are independent of $\tau$. For a global solution to which
these estimates apply, $\int_t^\infty\mathscr D(s)\dd s=\mathcal E(t)$.
\end{lemma}
\begin{proof}
For $\gamma>1$, the product $[P(s)-P(\bar\rho)][G(s)-G(\bar\rho)]$
controls $h_\gamma(s)$ on $[0,\infty)$. Both vanish quadratically
at $\bar\rho$, and their ratio tends to $(\gamma-1)/2$ at infinity.
For $\gamma=1$ the same comparison holds on any bounded density
interval. Thus $\int_0^R h_\gamma(\rho)r\dd r\le C\mathcal C$.
Equation~\eqref{pr:iso-kinetic-dissipation} controls the kinetic
energy, so \eqref{pr:iso-modified-energy} gives $\Phi'+\lambda\Phi\le0$.
Since $\mathcal E'=-\mathscr D$, this proves \eqref{ext:energy-decay}.
At $b=-a$ the independent endpoint density bounds cited above give this
exponential estimate also for $\gamma=1$.

It remains to treat $\gamma=1$ and $b>-a$ before a uniform upper
density bound is available. Mass conservation gives
\begin{equation}\label{pr:iso-covariance}
 \mathcal C=\frac{\kappa}{2}\int_0^R
        \frac{(\rho-\bar\rho)^2}{\rho+\bar\rho}\,r\dd r .
\end{equation}
Since $\beta>-1$, choose $\delta>0$ with $(1-\beta)(1+2\delta)<2$.
The coarse bound in Lemma~\ref{ext:uniform-weighted} gives
\begin{equation}\label{pr:iso-coarse}
 \sup_{t\le\tau}\int_0^R\rho^{1+2\delta}r\dd r\le C .
\end{equation}
For $\ell\ge1$ and $s\le\bar\rho e^\ell$,
$h_1(s)\le C\ell(s-\bar\rho)^2/(s+\bar\rho)$.
For $s>\bar\rho e^\ell$, $h_1(s)\le Cs^{1+\delta}
\le Ce^{-\delta\ell}s^{1+2\delta}$. Splitting at this level and using
\eqref{pr:iso-covariance}, \eqref{pr:iso-coarse} and the kinetic bound
\eqref{pr:iso-kinetic-dissipation}, we obtain
\[
 \Phi\le C\mathcal E\le C\ell(\mathscr D+\mathcal C)+Ce^{-\delta\ell}.
\]
Set $\ell=\sqrt{1+t}$. In this time variable \eqref{pr:iso-modified-energy}
gives $\dd\Phi/\dd\ell+c\Phi\le Ce^{-\delta\ell}$. Integration gives
$\mathcal E(t)\le Ce^{-c'\sqrt{1+t}}$, and the mechanical identity gives
\eqref{pr:iso-energy}.

For all $\gamma\ge1$, H\"older's inequality on each interval
$[n,\min(n+1,\tau)]$ bounds the integral of $\mathscr D^s$ there by
$\mathcal E(n)^s$. Summing over integers $0\le n<\tau$ proves
\eqref{pr:iso-fractional-coefficients}.
\end{proof}
\section{A priori estimates for positive second viscosity}
\label{sec:positive-density}
Assume $b>0$, so that $0<\beta\le\theta<1$. The pressure bound below
holds for $1/2\le\alpha\le1$. In the remainder of this section we take
$\alpha>1/2$. We return to $\alpha=1/2$ in Section~\ref{pr:half-all-pressure}.

\subsection{A pointwise bound for the radial pressure}
\begin{lemma}\label{pos:mass-pressure}
Let $b>0$, $1/2\le\alpha\le1$ and $\gamma\ge1$. On every positive
strong-solution interval $[0,\tau]$ one has
\begin{equation}\label{eq:radial-pressure-bound}
 \sup_{0\le t\le\tau}\sup_{0<r\le R} r^2\rho(t,r)^\gamma\le C
\end{equation}
with $C$ independent of $\tau$.
\end{lemma}
\begin{proof}
For $\gamma\le2+b/a$, the assertion follows from
$\rho\le Cr^{-2/(2+b/a)}$ in Lemma~\ref{ext:uniform-weighted}.
Otherwise $\gamma(1-\theta)>2\alpha\ge1$. Write $z=2\alpha-1$ and let $m(t,r)=\int_0^r\rho(t,y)y\dd y$.
By the continuity equation, $\Dt m=0$. By H\"older and \eqref{eq:energy},
\[
 m\le Cr^{2(1-1/\gamma)},\qquad
 s=\frac{\gamma\theta}{\gamma-1}\in(0,1),\qquad m^sr^{-2\theta}\le C .
\]
To bound the positive part of the effective velocity, set
$Z=\frac12\int_0^R m^s\rho w_+^2 r^{1-2\theta}\dd r$ and
\[
 A=\int_0^R m^s\rho^{\gamma+1-\alpha}w_+^2 r^{1-2\theta}\dd r,\qquad
 I=\int_0^R m^s\rho^{\gamma+1-\alpha}u^2 r^{1-2\theta}\dd r .
\]
Testing \eqref{eq:wz} with $m^sr^{-\theta}w_+$ and integrating
with respect to the transported mass $\rho r\dd r$ gives
\[
 Z'+\frac{\kappa\gamma}{2(2a+b)}A\le\frac{\kappa\gamma}{2(2a+b)}I .
\]
The geometric term is nonpositive and $2w_+u\le w_+^2+u^2$.
First test on $m\ge\varepsilon$. The inner boundary is material, the
wall flux vanishes, and \eqref{eq:centerregularity} justifies
$\varepsilon\downarrow0$.
Since $(2a+b)\rho^{\alpha-1}[(\log\rho)_r]_+=(w-u)_+\le w_++|u|$,
\[
 S_+=\int_0^R m^s\rho^{\gamma+\alpha-1}[(\log\rho)_r]_+^2
                  r^{1-2\theta}\dd r\le \frac{2}{(2a+b)^2}(A+I).
\]

We bound $I$ using $S_+$ only where $r^2\rho^\gamma$ is large.
Fix $r$ and write $L=\rho(t,r)$. If $\rho\ge L/2$ on $[r/2,r]$,
then $\int_{r/2}^r\rho^\gamma y\dd y\ge 3r^2L^\gamma/(8\,2^\gamma)$.
Choose $K_0>R^2$ large enough that the last lower bound contradicts
\eqref{eq:energy} whenever $r^2L^\gamma>K_0$. In that case $L>1$.
Let $r_0$ be the largest radius in $[r/2,r)$ for which $\rho(t,r_0)=L/2$.
On $[r_0,r]$ we have $\rho\ge L/2$, so $m_y=\rho y\ge Lr/4$ and
\[
 0<m(r)-m(r_0)
 \le (L/2)^{1-\gamma}\int_{r_0}^r\rho^\gamma y\dd y
 \le C L^{1-\gamma}.
\]
Changing variables from $y$ to the enclosed mass gives
\[
 \int_{r_0}^r m(y)^{-s}\dd y
 \le \frac{4}{Lr}\frac{m(r)^{1-s}-m(r_0)^{1-s}}{1-s}
 \le \frac{4}{(1-s)Lr}\,[m(r)-m(r_0)]m(r)^{-s}.
\]
Here $x^{1-s}\ge x$ for $0\le x\le1$. Inserting the mass increment
bound gives $C r^{-1}L^{-\gamma}m(r)^{-s}$. Alternatively,
$[m(r)-m(r_0)]m(r)^{-s}\le[m(r)-m(r_0)]^{1-s}$ gives
$C r^{-1}L^{-\gamma+s(\gamma-1)}$.
Both constants are independent of $K_0$ and $\tau$.
Integrating the positive derivative of $\rho^{(\gamma+\alpha-1)/2}$ from
$r_0$ to $r$ and using Cauchy--Schwarz therefore gives
\[
 L^{\gamma+\alpha-1}\le C S_+ r^{2\theta-1}\int_{r_0}^r m(y)^{-s}\dd y,
 \qquad m(r)^s r^{2-2\theta}L^{2\gamma+\alpha-1}\le C S_+ .
\]
In this region $m^sr^{-2\theta}L^{\gamma-\alpha}
\le CS_+/(r^2 L^{\gamma+z})\le CS_+/K_0$.
On $r^2L^\gamma\le K_0$ we have
$m^sr^{2-2\theta}L^{\gamma-z}\le C(K_0+R^2)$, since
$m^sr^{-2\theta}\le C$, $\gamma-z\ge0$ and $z\ge0$.
Using the kinetic energy bound for the first region and the angular
dissipation for the second yields
\[
 I\le C K_0^{-1} S_+ + C_{K_0}\int_0^R\rho^\alpha u^2\,\frac{\dd r}{r}.
\]
The coefficient of $K_0^{-1}S_+$ is independent of $K_0$. Choosing
$K_0$ sufficiently large and using $S_+\le 2(A+I)/(2a+b)^2$, we obtain
\[
 I\le\tfrac13(A+I)+C_{K_0}\int_0^R\rho^\alpha u^2\,\frac{\dd r}{r},
 \qquad
 Z'+\frac{\kappa\gamma}{4(2a+b)}A
 \le C_{K_0}\int_0^R\rho^\alpha u^2\,\frac{\dd r}{r}.
\]
The initial data and $m^sr^{-2\theta}\le C$ bound $Z(0)$.
Integrating and using \eqref{eq:energy} yields
$\sup_{t\le\tau} Z(t)+\int_0^\tau(A+I+S_+)\dd t\le C$.

Finally \eqref{eq:w}, $m^sr^{-2\theta}\le C$ and the kinetic energy
bound imply
\[
 \int_0^R m^s\rho^z[(\log\rho)_r]_+^2r^{1-2\theta}\dd r
 \le C(1+Z)\le C .
\]
On $[r_0,r]$ the net increase in the primitive of $h^{z/2-1}$ is
\[
 \int_{L/2}^L h^{z/2-1}\dd h
 \ge (L/2)^{z/2}\log 2 .
\]
For $z=0$ the integral equals $\log 2$.
Cauchy--Schwarz with the preceding positive-gradient bound and the
second estimate for $\int_{r_0}^r m(y)^{-s}\dd y$ now gives
\[
 L^z\le C r^{2\theta-1}\int_{r_0}^r m(y)^{-s}\dd y
      \le C r^{2\theta-2}L^{-\gamma+s(\gamma-1)} .
\]
Since $s+\gamma(1-s)=\gamma(1-\theta)$, this bounds
$(r^2L^\gamma)^{1-\theta}L^z$. Here $L>1$, so the claimed bound
holds on this region as well.
\end{proof}

\subsection{A weighted quadratic estimate and an inward velocity moment}
Assume $b>0$, $1/2<\alpha\le1$ and $\gamma\ge1$. We shall show that
the weighted density-gradient integral $F$ defined below grows at most
linearly in time. Put $z=2\alpha-1$ and $\omega=r^{-2\theta}$.
Multiply the physical kinetic energy by $c_u=\theta/(\alpha\beta)$
to cancel the mixed viscous term, and set
\[
 \begin{aligned}
 E&=\int_0^R\omega\left\{\frac{\rho}{2}(c_u u^2+w^2)
                   +(c_u+1)\Pi\right\}r\dd r+C_*,\\
 F&=\int_0^R\omega\rho^z|(\log\rho)_r|^2r\dd r,\qquad
 S=\int_0^R\omega\rho^{\gamma+\alpha-1}|(\log\rho)_r|^2r\dd r,\\
 D_\omega&=\int_0^R\omega\rho^\alpha(u_r^2+u^2/r^2)r\dd r .
 \end{aligned}
\]
Take $C_*=0$ if $\gamma>1$. At $\gamma=1$ choose $C_*$ so the
sum of the weighted potential and $C_*$ is nonnegative, using
$s\log(s/\bar\rho)\ge-\bar\rho/e$ and $\theta<1$. By \eqref{eq:w},
$F+\int_0^R\omega\rho u^2r\dd r\le CE$.

We first record the viscous form for arbitrary constant multipliers
$c_u\omega u$ and $c_w\omega w$. Write $x=u_r$ and $y=u/r$.
In the effective-velocity test, integrate the term
$-c_w\theta\int\omega\rho(w-u)u^2\dd r$ by parts using
$\rho(w-u)=(2a+b)(\rho^\alpha)_r/\alpha$. After factoring out
$\omega\rho^\alpha r$, the viscous integrand is
\begin{equation}\label{eq:weighted-viscous-form}
 (2a+b)\left\{
 c_u[x^2+2(\beta-\theta)xy+(1-2\beta\theta)y^2]
 -\frac{2c_w\theta}{\alpha}(xy-\theta y^2)\right\}.
\end{equation}
For $c_w=1$ and $c_u=\theta/(\alpha\beta)$, the coefficient of
$(x-y)(x+y)$ vanishes. The energy identity becomes
\begin{equation}\label{pht:quadratic-identity}
\begin{aligned}
 E'&+a\int_0^R\omega\rho^\alpha\left[
 (1+\theta)\left(c_u+\frac{\theta}{1-\theta}\right)(x-y)^2
 +\frac{\theta(2a+b)}{b}(x+y)^2\right]r\dd r
 +\kappa\gamma(2a+b)S\\
 &=-\theta(c_u+1)\int_0^R\omega\rho u^3\dd r
   -2\theta(c_u+1)\int_0^R\omega(P+\Pi)u\dd r .
\end{aligned}
\end{equation}
Both square coefficients are positive since $b>0$. For the pressure we used
\[
 \frac{\dd}{\dd t}\int_0^R\omega\Pi r\dd r
 =-\int_0^R\omega P(u_r+u/r)r\dd r-2\theta\int_0^R\omega\Pi u\dd r .
\]
The wall fluxes vanish by no slip and the central fluxes by
\eqref{eq:centerregularity}.

We use the following estimate twice, with $k=\gamma+\alpha-1$ and
$k=2\alpha-1$. Let $N=1+\max_r\rho$. For $0<\theta<1$ and $k>0$,
\begin{equation}\label{eq:weighted-density-maximum}
 N^{k+\gamma\theta}\le C_k\left(1+
 \int_0^R|\partial_r\rho^{k/2}|^2 r^{1-2\theta}\dd r\right).
\end{equation}
For $L=\max\rho$ larger than twice the uniform wall bound and $2$,
choose a maximum point $p$ and the first $q>p$ with $\rho(t,q)=L/2$.
Then $q^2-p^2\le CL^{-\gamma}$ by \eqref{eq:energy}. Cauchy--Schwarz
and the concavity of $x^\theta$ give
\[
 (1-2^{-k/2})^2 L^k
 \le \int_p^q|(\rho^{k/2})_r|^2r^{1-2\theta}\dd r\,
       \frac{q^{2\theta}-p^{2\theta}}{2\theta}
 \le C L^{-\gamma\theta}\int_0^R|(\rho^{k/2})_r|^2 r^{1-2\theta}\dd r .
\]
This proves \eqref{eq:weighted-density-maximum}, including $p=0$.
Bounded $L$ is covered by the constant on the right.
With $k=\gamma+\alpha-1$ this reads
\begin{equation}\label{pht:pressure-peak}
 N^{\gamma+\alpha-1+\gamma\theta}\le C(1+S).
\end{equation}

By Young's inequality the absolute value of the pressure term in
\eqref{pht:quadratic-identity} is at most
\[
 \varepsilon D_\omega+
 C_\varepsilon\int_0^R\frac{(P+\Pi)^2}{\rho^\alpha}r^{1-2\theta}\dd r .
\]
For $\gamma>1$ the quotient is a constant times $\rho^{2\gamma-\alpha}$.
Take $\eta=0$ when $\gamma>1$ and $\eta=z/2$ when $\gamma=1$. In the latter
case $(P+\Pi)^2/\rho^\alpha=\kappa^2\rho^{2-\alpha}(1+\log(\rho/\bar\rho))^2$.
For $0<\rho\le1$ this is bounded by $C\rho^{2-\alpha-\eta}\le C$,
and for $\rho\ge1$ by $C\rho^{2-\alpha+\eta}$. Here we used
$|\log s|^2\le C_\eta\max(s^\eta,s^{-\eta})$ and $2-\alpha-\eta>0$.
Put $\sigma=2\gamma-\alpha+\eta\ge\gamma$ and $d=RN^{-\gamma/2}$.
Splitting at $d$ and using \eqref{eq:energy} gives
\[
 \int_0^R \rho^\sigma r^{1-2\theta}\dd r
 \le C N^\sigma d^{2-2\theta}
      +N^{\sigma-\gamma}d^{-2\theta}\int_d^R\rho^\gamma r\dd r
 \le C N^{\sigma-\gamma+\gamma\theta}.
\]
The exponent is below that in \eqref{pht:pressure-peak} by $z-\eta>0$.
After absorption, with $v=u_-=\max\{-u,0\}$,
\begin{equation}\label{pht:reduced-quadratic}
 E'+c_0(D_\omega+S)\le C\int_0^R\omega\rho v^3\dd r+C .
\end{equation}

We first estimate this cubic term on $\rho\ge1$. In the coordinate
$y=\log r$ consider an interval directed toward the center with
length $\ell=c_\alpha\omega\rho^z/(F+\omega\rho^z)\le1$, where
$\rho$ and $\omega$ on the right are evaluated at its initial point.
Cauchy--Schwarz applied to $\rho^{z/2}$ shows that the density stays
between half and twice its initial value there if $c_\alpha$ is
sufficiently small. The one-dimensional trace inequality then gives
\begin{equation}\label{eq:fix-localtrace}
 |u|\le C \mathscr D^{1/2}\rho^{-\alpha/2}
       (1+F^{1/2}\omega^{-1/2}\rho^{-z/2}).
\end{equation}
Put $M=\int_0^R\omega\rho u^2 r\dd r\le CE$. On $\rho\ge1$,
both $\rho^{1-\alpha/2}$ and $\rho^{1-(\alpha+z)/2}$ are at most
$\rho^{(1+\alpha)/2}$. Multiplying \eqref{eq:fix-localtrace} by
$\omega\rho u^2$ and using the angular terms of $D_\omega$ and $\mathscr D$
gives
\begin{equation}\label{dic:high}
 \int_{\rho\ge1}\omega\rho v^3\dd r
 \le C(\mathscr D\,D_\omega M)^{1/2}+C\mathscr D(FM)^{1/2}
 \le\varepsilon D_\omega+C_\varepsilon \mathscr D(F+M).
\end{equation}

Choose $J=5+b/a$ and $\delta=J-2-2\theta>0$. Define
\[
 M_J=\frac1J\int_0^R\rho v^Jr^{1+\delta}\dd r,\qquad
 U_J=\int_0^R\rho v^{J+1}r^\delta\dd r,\qquad
 D_J=\int_0^R\rho^\alpha v^{J-2}(u_r^2+u^2/r^2)r^{1+\delta}\dd r .
\]
On $\rho<1$ the identity
\[
 \omega\rho v^3
 =\rho^{(1-\alpha)(J-2)/(J-1)}
  (\omega\rho^\alpha v^2/r)^{(J-2)/(J-1)}
  (\rho v^{J+1}r^\delta)^{1/(J-1)}
\]
bounds the integral by $\varepsilon D_\omega+C_\varepsilon U_J$.

Test radial momentum by $-v^{J-1}r^{1+\delta}$ against $\dd r$.
On $u<0$ the viscous quadratic form, divided by
$(2a+b)\rho^\alpha v^{J-2}r^{1+\delta}$, is
\[
 (J-1)u_r^2+(J\beta+\delta)u_r\,u/r+(1+\beta\delta)(u/r)^2.
\]
With $J=5+b/a$, its determinant multiplied by $(2+b/a)^2$ is
$4(1+b/a)(4+b/a)-(1+2\alpha)^2>0$, since $b>0$ and $\alpha\le1$.
The transported weight contributes $\delta U_J/J$ on the left.
Integration by parts in the pressure gives a nonpositive angular
term $-(1+\delta)\int_0^RPv^{J-1}r^\delta\dd r$. The boundary terms
vanish by $u(R)=0$ and $u=O(r)$ at the center. Hence
\[
 M_J'+c_J D_J+\frac{\delta}{J}U_J
 \le (J-1)\int_0^R P v^{J-2}|u_r|r^{1+\delta}\dd r .
\]

On the set $q=rv\le1$, Cauchy--Schwarz bounds the pressure contribution by
\[
 C(\mathscr D\,M_J)^{1/2}
 \left(\sup_{q\le1}r^{2-2\theta}
              \rho^{2\gamma-\alpha-1}q^{J-4}\right)^{1/2}.
\]
The wall density is bounded. Integrating
$\partial_r\rho^{(2\gamma+\alpha-1)/2}$ from $r$ to $R$ and using
Cauchy--Schwarz gives
\begin{equation}\label{eq:weighted-pressure-trace}
 \rho(t,r)^{(2\gamma+\alpha-1)/2}
 \le C+C S^{1/2}\left(\int_r^R\rho^\gamma y^{2\theta-1}\dd y\right)^{1/2}
 \le C+C r^{\theta-1}S^{1/2}.
\end{equation}
Because $2\gamma-\alpha-1\ge0$, we have
$\rho^{2\gamma-\alpha-1}\le 1+\rho^{2\gamma+\alpha-1}$.
By \eqref{eq:weighted-pressure-trace},
$r^{2-2\theta}\rho^{2\gamma-\alpha-1}\le C(1+S)$.
Since $q^{J-4}\le1$ on $\{q\le1\}$, the pressure contribution there is at most
\[
 C(\mathscr D M_J)^{1/2}(1+S)^{1/2}
 \le \varepsilon(1+S)+C_\varepsilon\mathscr D M_J .
\]
On $\{q>1\}$, Young's
inequality bounds the pressure contribution by $c_JD_J/2$ plus the
integral of $C_J\rho^{2\gamma-\alpha}v^{J-2}r^{1+\delta}$.
The ratio of the latter integrand to that of $U_J$, up to $C_J$, is
\[
 \frac{r^4\rho^{2\gamma-\alpha-1}}{q^3}
 =\frac{(r^2\rho^\gamma)^2}{\rho^{1+\alpha}q^3}
 =\frac{(r^2\rho^\gamma)^{\,2-(1+\alpha)/\gamma}
                  r^{\,2(1+\alpha)/\gamma}}{q^3}.
\]
Since $2-(1+\alpha)/\gamma\ge0$, \eqref{eq:radial-pressure-bound}
bounds the ratio by $C r^{2(1+\alpha)/\gamma}/q^3$. Thus $C_J$
times the ratio is at most $\delta/(2J)$ for $r$ sufficiently small
or $q$ sufficiently large. Elsewhere in $q>1$, the radius is bounded
away from zero, so both $v=q/r$ and $\rho$ are bounded. The remaining
integral is therefore bounded uniformly in time.
Thus
\[
 M_J'+\tfrac12c_JD_J+\frac{\delta}{2J}U_J
 \le\varepsilon S+C_\varepsilon(1+\mathscr D M_J).
\]

Use \eqref{dic:high} and the low-density cubic estimate in
\eqref{pht:reduced-quadratic}. Absorb the small multiple of $D_\omega$,
add a fixed large multiple $\Lambda$ of the moment inequality to
absorb $U_J$, and then take its $\varepsilon$ small enough to absorb
$\Lambda\varepsilon S$. For $\mathcal Y=1+E+\Lambda M_J$ this gives
\begin{equation}\label{pht:joint-energy}
 \mathcal Y'+c_*(D_\omega+S+D_J+U_J)\le C \mathscr D\mathcal Y+C .
\end{equation}
The initial data bound $\mathcal Y(0)$, and
$\int_0^\tau\mathscr D\dd t\le\mathcal E(0)$. Gronwall's inequality
therefore gives $F(t)\le C\mathcal Y(t)\le C(1+t)$ for $0\le t\le\tau$.

\subsection{Density bounds and continuation}
By \eqref{pht:joint-energy} and mass conservation,
\[
 \begin{aligned}
 |\rho(t,r)^{z/2}-\rho(t,s)^{z/2}|
 &\le C(1+t)^{1/2}
       \left|\int_s^r y^{2\theta-1}\dd y\right|^{1/2},\\
 \|\rho(t)\|_{L^\infty}
 &\le C(1+t)^{1/z},\qquad
 [\rho(t)]_{C^{0,\min\{\theta,1/2\}}(\overline\Omega)}
 \le C(1+t)^A
 \end{aligned}
\]
for some fixed $A>0$. The upper bound follows by taking $s$ with
$\rho(t,s)=\bar\rho$. The modulus follows by applying the mean
value inequality to the power $2/z$ on the resulting density range.
Fix $T<\infty$ and consider a positive strong solution on $[0,\tau]$ with
$\tau\le T$. By the upper bound just proved,
$K=\kappa\gamma\rho^{\gamma-\alpha}/(2a+b)\le C_T$ since $\gamma\ge\alpha$.
Completing the square in \eqref{eq:wz} and integrating along trajectories
gives $w\le C_{+,T}r^\theta$, with $C_{+,T}$ independent of $\tau$.
The physical velocity satisfies
\[
\Dt u=\nu(u_{rr}+u_r/r-u/r^2)
+(w-u)[\alpha u_r+\alpha\beta u/r-K].
\]
At a contact with a comparison function $f$ the scalar residual is
\begin{equation}\label{eq:velocity-residual}
 \mathscr P[f]=f_t+ff_r-\nu(f_{rr}+f_r/r-f/r^2)
       -(w-f)[\alpha(f_r+\beta f/r)-K].
\end{equation}
Direct substitution gives
\begin{align*}
 \mathscr P[Br^\theta]
&=\nu B(1-\theta^2)r^{\theta-2}
+\theta B^2r^{2\theta-1}\\
&\quad+(B-r^{-\theta}w)r^{2\theta-1}
[\alpha(\theta+\beta)B-Kr^{1-\theta}],\\
 \mathscr P[-Ar^{-\beta}]
&=-\nu A(1-\beta^2)r^{-\beta-2}
-\beta A^2r^{-2\beta-1}+K(w+Ar^{-\beta}).
\end{align*}
Since $K\le C_T$ and $w\le C_{+,T}r^\theta$, a sufficiently large
$B=B_T$ makes both factors in the last term of $\mathscr P[Br^\theta]$
nonnegative. For the lower barrier, the term
$-\beta A^2r^{-2\beta-1}$ dominates $K(w+Ar^{-\beta})$ uniformly on
$(0,R]$ for large $A=A_T$, since $\beta>0$. Increase $A$ and $B$
to enclose the initial data. Both comparisons hold at the wall, and
$u=O_\tau(r)$ ensures that they hold at every sufficiently small inner
radius $\epsilon$. Comparison on $[\epsilon,R]$ and passage to the limit give
\begin{equation}\label{eq:barriers}
w\le C_{+,T}r^\theta,\qquad -A_T r^{-\beta}\le u\le B_T r^\theta.
\end{equation}

Let $M_T$ be the preceding upper bound for $\rho$. Mass conservation
gives a fixed $r_a>0$ such that
$\int_0^{r_a}\rho(t,y)y\dd y\le R^2\bar\rho/4$ for all $0\le t\le\tau$.
For each such time there is $s\in[r_a,R]$ with $\rho(t,s)\ge\bar\rho/2$.
Choose $\phi_\alpha'(\rho)=(2a+b)\rho^{\alpha-2}$.
Since $\partial_r\phi_\alpha(\rho)=w-u$, \eqref{eq:barriers} gives,
for $0\le r\le s$,
\[
 \phi_\alpha(\rho(t,r))\ge
 \phi_\alpha(\bar\rho/2)-\int_0^R(C_{+,T} y^\theta+A_T y^{-\beta})\dd y.
\]
Here $\beta<1$ and $\phi_\alpha(0+)=-\infty$, so $\rho$ has a
positive lower bound throughout $[0,s]$, in particular on $[0,r_a]$.
Every trajectory ending at a radius at least $r_a$ therefore carries
an enclosed mass at least $c_T r_a^2$. At all earlier times the
mass inside $\xi$ is at most $M_T\xi^2/2$, so this trajectory stays
a fixed positive distance from zero. Its velocity is bounded by
\eqref{eq:barriers}, and integration of \eqref{eq:wz} along it
bounds $|w|$ there. For each $r\in[s,R]$, integrate
$\partial_r\phi_\alpha(\rho)=w-u$ from $s$ to $r$. This gives
a positive lower bound for $\rho(t,r)$ on the annulus, with a
constant depending on $T$ and independent of $\tau$. The gradient
estimate \eqref{eq:qbound} and the continuation argument in
Section~\ref{sec:higher-continuation} thus give global existence.

On the density range obtained above,
\[
 |\rho-\bar\rho|^2\le C(1+\|\rho(t)\|_{L^\infty})^{(2-\gamma)_+}
                          h_\gamma(\rho).
\]
For $\gamma=1$ this follows from $h_1''(s)=\kappa/s$.
Lemma~\ref{ext:early-decay} gives exponential
mechanical decay for $\gamma>1$ and decay at least
$Ce^{-c\sqrt{1+t}}$ for $\gamma=1$. Consequently
$\|\rho-\bar\rho\|_{L^2}\le Ce^{-c\sqrt{1+t}}$ for all $\gamma\ge1$.
To pass to the supremum norm, put $g=\rho(t)-\bar\rho$ and
$\eta=\min\{\theta,1/2\}$. Averaging over $\Omega\cap B_d(x)$ gives,
uniformly in $x\in\overline\Omega$ and $0<d\le R$,
\[
 |g(x)|\le [g]_{C^\eta}d^\eta+C d^{-1}\|g\|_{L^2(\Omega)}.
\]
Since $[g]_{C^\eta}$ grows at most polynomially, taking
$d=R\exp(-c\sqrt{1+t}/2)$ gives the same stretched exponential
bound in $L^\infty$, with a smaller constant in the exponent.
In particular $\bar\rho/2\le\rho\le2\bar\rho$ after a fixed time.
The finite-time bounds cover the preceding interval. With this uniform
upper bound, Lemma~\ref{ext:early-decay} gives exponential mechanical
energy decay also for $\gamma=1$. A second application of the same interpolation gives
$\|\rho-\bar\rho\|_{L^\infty}\le Ce^{-ct}$. The higher-order
conclusions follow in Section~\ref{ext:higher-section}.

\section{Nonpositive second viscosity}

\subsection{Positive geometric coefficient}\label{dic:section}
Assume $1/2<\alpha<1$, $-a\le b\le0$, $\theta>0$ and $\gamma\ge1$.
Then $0<\theta\le1-\alpha<1/2$ and $\beta>-1$ by \eqref{eq:theta}.
With $\omega=r^{-2\theta}$, define
\begin{equation}\label{eq:fix-energydef}
\begin{aligned}
 E(t)&=\int_0^R\omega\left\{
 \frac{\rho}{2}[(1+\theta)u^2+(1-\theta)w^2]+2\Pi\right\}r\dd r,\\
 F(t)&=\int_0^R\omega\rho^{2\alpha-1}|(\log\rho)_r|^2r\dd r,\\
 D_\omega(t)&=\int_0^R\omega\rho^\alpha(u_r^2+u^2/r^2)r\dd r .
\end{aligned}
\end{equation}
For $\gamma=1$ add $1+2\kappa\bar\rho e^{-1}\int_0^R\omega r\dd r$
to $E$ so its shifted potential part is positive.
The relation $w-u=(2a+b)\rho^{\alpha-1}(\log\rho)_r$ shows that
$F+\int\omega\rho u^2r\dd r\le CE$.
All initial weighted integrals are finite by radial $H^3$ regularity
and $\theta\le1/2$.

For the viscous term set
\[
 C_s=a(1+2\theta),\qquad
 C_d=(1+\theta)(a+b)-a\theta=(a+b)+b\theta .
\]
Multiplication of the physical and effective velocity equations by
$(1+\theta)\omega u$ and $(1-\theta)\omega w$, respectively, gives
\begin{equation}\label{eq:fix-energyidentity}
\begin{aligned}
 E'&+\int_0^R\omega\rho^\alpha
 \bigl\{(1+\theta)C_s(u_r-u/r)^2+(1-\theta)C_d(u_r+u/r)^2\\
 &\hspace{14mm}+\theta(C_d-C_s)(u_r^2-u^2/r^2)\bigr\}r\dd r
 +(1-\theta)\int_0^R\omega\rho K(w-u)^2r\dd r\\
 &=-2\theta\int_0^R\omega\rho u^3\dd r
 -4\theta\int_0^R\omega(P+\Pi)u\dd r .
\end{aligned}
\end{equation}
This follows from \eqref{eq:weighted-viscous-form} with
$(c_u,c_w)=(1+\theta,1-\theta)$. Their sum is two, giving the cubic
and pressure terms in \eqref{eq:fix-energyidentity}. The boundary
fluxes vanish as in \eqref{pht:quadratic-identity}.

The quadratic form in \eqref{eq:fix-energyidentity} is strictly positive.
Indeed, $0<C_d\le a\le C_s$ and
\[
 \frac{C_d}{a}
 =\frac{(2\alpha-1)(1+\theta)+2\theta^2}{1-\theta}.
\]
Writing $v=C_d/C_s\le1$, we have
$v\ge2\theta^2/[(1-\theta)(1+2\theta)]>\theta^2$.
The determinant divided by $C_s^2$ is
\[
 (1-\theta^2)v-\frac{\theta^2}{4}(v-1)^2
 =v-\frac{\theta^2}{4}(v+1)^2\ge v-\theta^2>0.
\]
Consequently its integral controls a positive fixed multiple of
$D_\omega$.

We estimate $\int_0^R \omega\rho u_-^3\dd r$ with
$u_-=\max\{-u,0\}$. Set $z=2\alpha-1>0$.
By \eqref{eq:theta}, $2\theta<1+\beta$.
Set
\[
 M=\int_0^R\omega\rho u^2r\dd r,\qquad
 S=\int_0^R\omega\rho^{\gamma+\alpha-1}|(\log\rho)_r|^2r\dd r .
\]
Then $F+M\le CE$.

On $\{\rho\ge1\}$, the proof of \eqref{dic:high} uses only $z>0$ and
$0<\theta<1$, so it applies here as well.

\paragraph{The cubic term at low density.}
Choose $j$ and $\delta$ by
\begin{equation}\label{dic:moment-choice}
 j=2+2\theta+\sigma,\qquad \delta=\sigma/2,\qquad
 0<\sigma<\min\left\{1-2\theta,\frac{1+\beta}{4},
                    \frac{z+2\alpha\theta}{1+\gamma}\right\}.
\end{equation}
Such a $\sigma$ exists since $\theta<1/2$, $\beta>-1$ and $z>0$.
Let
\[
 \begin{aligned}
 M_j&=\frac1j\int_0^R\rho u_-^j r^{1+\delta}\dd r,&
 U_j&=\int_0^R\rho u_-^{j+1}r^\delta\dd r,\\
 D_j&=\int_0^R\rho^\alpha u_-^{j-2}(u_r^2+u^2/r^2)
                                      r^{1+\delta}\dd r .
 \end{aligned}
\]
On $\{\rho<1\}$, the relation $j-2-2\theta=2\delta$ gives
\[
 \omega\rho u_-^3
 \le r^{\delta/(j-1)}
   (\omega\rho^\alpha u_-^2/r)^{(j-2)/(j-1)}
   (\rho u_-^{j+1}r^\delta)^{1/(j-1)} .
\]
On $r<r_0$ choose $r_0$ small and apply Young's inequality.
On $r\ge r_0$ the weights are bounded above and away from zero,
so $u_-^3\le\varepsilon' u_-^{j+1}+C_{\varepsilon'}u_-^2$ gives
\begin{equation}\label{dic:low}
 \int_{\rho<1}\omega\rho u_-^3\dd r
 \le \varepsilon(D_\omega+U_j)+C_\varepsilon\int_0^R\rho u^2 r\dd r .
\end{equation}

\paragraph{The pressure in the inward moment.}
Testing with $-u_-^{j-1}r^{1+\delta}$ gives the normalized viscous
form $(j-1)X^2+(j\beta+\delta)XY+(1+\beta\delta)Y^2$ in
$(X,Y)=(u_r,u/r)$. From \eqref{dic:moment-choice},
$j-2<5(1+\beta)/4$ and $\delta-\beta(j-2)<3(1+\beta)/2$. Its determinant obeys
\[
 (j-1)(1-\beta^2)-\tfrac14\{\delta-\beta(j-2)\}^2
 \ge (1+\beta)-\tfrac9{16}(1+\beta)^2
 \ge\tfrac7{16}(1+\beta)>0 .
\]
Integrating the pressure term by parts and applying Young's inequality
to the term containing $u_r$, we obtain
\[
 M_j'+c_jD_j+\frac{\delta}{j}U_j
 \le C_j\int_0^R\rho^{2\gamma-\alpha}u_-^{j-2}r^{1+\delta}\dd r .
\]
The term $\delta U_j/j$ comes from $\Dt(r^\delta)=\delta ur^{\delta-1}$
on $\{u<0\}$. The no-slip condition and $u=O(r)$ at the center remove
the boundary terms.

Fix $N_0>2\bar\rho+1$. On $\{\rho\le N_0\}$ the last integral
is at most $C\mathcal E^{(j-2)/2}$ by H\"older, since $2<j<3$.
On $\{\rho>N_0\}$ factor out the angular dissipation density
$\rho^\alpha u_-^j r^{\delta-1}$ to the power $(j-2)/j$. Young's
inequality with conjugate exponents $j/(j-2)$ and $j/2$ gives a
bound by $\varepsilon D_j$ plus a constant times
\[
 I=\int_{\rho>N_0}\rho^L r^{s-1}\dd r,\qquad
 L=\alpha+j(\gamma-\alpha),\quad s=j+\delta .
\]
Write $Q=\gamma+\alpha-1+\gamma\theta$. The choices above give
\begin{equation}\label{dic:pressure-margins}
 \begin{aligned}
 (2\gamma+\alpha-1)-(L-\gamma)
    &=2\gamma+\alpha-1-(j-1)(\gamma-\alpha)>3\alpha-1>0,\\
 Q-L+\gamma s/2
    &=z+2\alpha\theta+(\alpha-\gamma/4)\sigma>0 .
 \end{aligned}
\end{equation}
The first uses $j<3$. For the second use \eqref{dic:moment-choice}.

Let $N=1+\max_{[0,R]}\rho$. We estimate the density powers in $I$
using
\begin{equation}\label{dic:two-traces}
 N^Q+\sup_{0<r\le R} r^{2(1-\theta)}\rho^{2\gamma+\alpha-1}
 \le C(1+S).
\end{equation}
The first inequality follows from \eqref{eq:weighted-density-maximum}
with $k=\gamma+\alpha-1$, and the second from
\eqref{eq:weighted-pressure-trace}. Their proofs use only $0<\theta<1$,
the wall density bound, and the bound for $\int_0^R\rho^\gamma r\dd r$.
Since $s>2$ and $L>\gamma$, take
$k=\min\{(L-\gamma)/(2\gamma+\alpha-1),(s-2)/(2(1-\theta))\}$.
Factoring off $(r^{2(1-\theta)}\rho^{2\gamma+\alpha-1})^k$ leaves
nonnegative powers of $\rho$ and $r$. By \eqref{dic:two-traces},
\[
 \rho^{L-\gamma}r^{s-2}\le C(1+S)^\eta,\qquad
 \eta=k+\frac{L-\gamma-k(2\gamma+\alpha-1)}{Q}
 =\max\left\{\frac{L-\gamma}{2\gamma+\alpha-1},\frac{L-\gamma s/2}{Q}\right\}<1 .
\]
Here $(2\gamma+\alpha-1)-Q=\gamma(1-\theta)$ and the strict bound follows
from \eqref{dic:pressure-margins}. On $\{\rho>N_0\}$ the relative
mechanical energy $\mathcal E$ controls $\int_{\rho>N_0}\rho^\gamma r\dd r$.
Consequently
\[
 I\le C\mathcal E(1+S)^\eta
 \le\varepsilon S+C_\varepsilon
        (\mathcal E+\mathcal E^{1/(1-\eta)}).
\]
By \eqref{ext:energy-decay} or \eqref{pr:iso-energy}, this remainder
and $\mathcal E^{(j-2)/2}$ are integrable in time. Thus
\begin{equation}\label{dic:raw}
 M_j'+c'_j D_j+\frac{\delta}{j}U_j
 \le \varepsilon S+f_\varepsilon(t),\qquad
 \int_0^\tau f_\varepsilon(t)\dd t\le C_\varepsilon .
\end{equation}

\paragraph{The weighted energy estimate.}
The cubic term on the right of \eqref{eq:fix-energyidentity} is bounded
above by $2\theta\int_0^R\omega\rho u_-^3\dd r$. To estimate the pressure
term, set $\epsilon_0=0$ when $\gamma>1$, and
choose $0<2\epsilon_0<z$ when $\gamma=1$. On $\rho>N_0$ we have
$|P+\Pi|\le C\rho^{\gamma+\epsilon_0}$. On $\{\rho\le N_0\}$,
$\rho^{-\alpha}|P+\Pi|^2$ is bounded, also when $\gamma=1$.
Cauchy--Schwarz with the angular part of $\mathscr D$ bounds this
contribution by $C\mathscr D^{1/2}$ since
$\int_0^Rr^{1-4\theta}\dd r<\infty$.
On $\{\rho>N_0\}$ Cauchy--Schwarz against the angular part of
$D_\omega$ leaves, after Young's inequality,
$C_\varepsilon\int_{\rho>N_0}\rho^{2\gamma+2\epsilon_0-\alpha}r^{1-2\theta}\dd r$.
For any radius $l>0$, the integral of $\rho^\gamma r^{1-2\theta}$
over $\{\rho>N_0\}$ is at most
$C N^\gamma l^{2-2\theta}+C l^{-2\theta}\mathcal E$.
This follows by splitting at $r=l$ and using the relative energy
on $r>l$. Taking the infimum over $l>0$ gives
\[
 \int_{\rho>N_0}\rho^{2\gamma+2\epsilon_0-\alpha}r^{1-2\theta}\dd r
 \le C N^{\gamma+2\epsilon_0-\alpha+\gamma\theta}\mathcal E^{1-\theta}
 \le C(1+S)^{1-(z-2\epsilon_0)/Q}\mathcal E^{1-\theta}.
\]
Since $0<(z-2\epsilon_0)/Q<1$, Young's inequality absorbs an
arbitrarily small part of $S$. The remainder is bounded by a
constant times
$\mathcal E^{1-\theta}+\mathcal E^{(1-\theta)Q/(z-2\epsilon_0)}$.
By \eqref{pr:iso-fractional-coefficients},
$\int_0^\tau\mathscr D^{1/2}\dd t\le C$. Add \eqref{dic:raw} to
\eqref{eq:fix-energyidentity}.
In \eqref{dic:high}, \eqref{dic:low} and \eqref{dic:raw}, choose
$\varepsilon$ small enough that the left-hand side controls
$D_\omega$, $U_j$ and $S$. We obtain
\begin{equation}\label{dic:budget}
 (E+M_j)'+c_*(D_\omega+S+D_j+U_j)\le C \mathscr D(E+M_j)+f(t),
 \qquad \int_0^\tau (\mathscr D+f)\dd t\le C .
\end{equation}
Gronwall gives
\begin{equation}\label{dic:uniform}
 \sup_{0\le t\le\tau}(E+F+M_j)(t)
 +\int_0^\tau(D_\omega+S+D_j+U_j)\dd t\le C .
\end{equation}
Applying \eqref{eq:weighted-density-maximum} with $k=2\alpha-1$ gives
\begin{equation}\label{eq:fix-peak}
 \|\rho(t)\|_{L^\infty}^{\,2\alpha-1+\gamma\theta}
 \le C(1+F(t)).
\end{equation}

\subsection{A lower bound for the density when \texorpdfstring{$\theta>0$}{theta > 0}}
Continue to assume $1/2<\alpha<1$, $b\le0$ and $\theta>0$. Then
$0<\theta\le1-\alpha<1/2$. Equations \eqref{eq:fix-peak} and
\eqref{dic:uniform} give
\begin{equation}\label{eq:fix-lowerinput}
\sup_{t\le\tau}\norm{\rho(t)}_{L^\infty}
+\sup_{t\le\tau}\int_0^R
\rho^{2\alpha-1}|(\log\rho)_r|^2r^{1-2\theta}\dd r
\le C.
\end{equation}
For $\tau\le T$, the constants below may depend on $T$, but are
independent of $\tau$ and $\inf_{[0,\tau]\times\Omega}\rho$.

Weighted Cauchy--Schwarz applied to $\rho^{\alpha-1/2}$ gives
\begin{equation}\label{eq:fix-densitymodulus}
|\rho(t,r)^{\alpha-1/2}-\rho(t,s)^{\alpha-1/2}|
\le C|r-s|^\theta,
\qquad \norm{\rho(t)}_{C^\theta([0,R])}\le C.
\end{equation}
The second estimate follows from the first and the density upper bound.

We first estimate the density near the center. Completing the square
in \eqref{eq:wz} bounds its right-hand side by
$K^2r^{1-\theta}/(4\theta)\le CK$, since $K$ is bounded. By
\eqref{eq:centerregularity}, $r^{-\theta}w_0$ is bounded as well.
An integrating factor along particle paths therefore gives
\begin{equation}\label{eq:fix-wupper}
w(t,r)\le C_+r^\theta.
\end{equation}

The enclosed mass satisfies
\begin{equation}\label{eq:fix-massequation}
M(t,r)=\int_0^r\rho(t,\xi)\xi\dd\xi,
\qquad M_r=\rho r,\qquad \Dt M=0,
\qquad M_t=\nu(\rho)(M_{rr}-M_r/r)-wM_r.
\end{equation}
Let $m_0=\min\rho_0>0$ and define $k$ and $\Phi$ by
\[
\phi_\alpha(k(r))=\phi_\alpha(m_0)
-\frac{C_+}{1+\theta}(R^{1+\theta}-r^{1+\theta}),
\qquad \Phi(r)=\int_0^rk(\xi)\xi\dd\xi.
\]
Here $0<k(0)\le k\le m_0$, since $\phi_\alpha(0+)=-\infty$, and
\[
\nu(k)(\Phi''-\Phi'/r)=C_+r^\theta\Phi'.
\]
Set $\Phi_\varepsilon(t,r)=\Phi(r)-\varepsilon(1+t)$.
Since $k\le m_0$, we have $M(0,r)>\Phi_\varepsilon(0,r)$.
The inequality stays strict at $r=0,R$ because the enclosed mass is
constant at both endpoints. At a first interior contact,
$M_r=\Phi'=rk$, so $\rho=k$, and $M_{rr}\ge\Phi''$.
The mass equation and \eqref{eq:fix-wupper} would then give
\[
 \partial_t(M-\Phi_\varepsilon)
 \ge (C_+r^\theta-w)\Phi'+\varepsilon>0,
\]
whereas first contact requires $\partial_t(M-\Phi_\varepsilon)\le0$.
Letting $\varepsilon\downarrow0$ and then taking the limit of
$2M(t,r)/r^2$ at the center gives
\begin{equation}\label{eq:fix-centralmass}
M(t,r)\ge\Phi(r)\ge\tfrac12k(0)r^2,
\qquad \rho(t,0)\ge k(0).
\end{equation}
Combining this with~\eqref{eq:fix-densitymodulus}, choose fixed
$r_c>0$ and $c_c>0$ such that
\begin{equation}\label{eq:fix-centralpositive}
\rho(t,r)\ge c_c\qquad(0\le r\le r_c,\ 0\le t\le\tau).
\end{equation}

It remains to estimate the density outside the central disk.
Choose $Q>(3\alpha-1)/(1-\alpha)$ sufficiently close to this threshold
that $2Q-1-\beta^2Q^2>0$. Such a choice is possible because
$|\beta|<(1-\alpha)/\alpha$ and
\[
2\frac{3\alpha-1}{1-\alpha}-1
-\beta^2\left(\frac{3\alpha-1}{1-\alpha}\right)^2
>
\frac{(2\alpha-1)(8\alpha^2-5\alpha+1)}
{\alpha^2(1-\alpha)}>0.
\]
In the momentum test with $|u|^{2Q-2}ur$, the normalized viscous
form has first diagonal entry $2Q-1$ and determinant
$2Q-1-\beta^2Q^2>0$. For the pressure remainder we use
$\rho^{2\gamma-\alpha}\le C\rho$, since $2\gamma-\alpha\ge1$ and
$\rho$ has a uniform upper bound. Hence
\begin{equation}\label{eq:fix-outermoment}
\sup_{t\le\tau}\int_0^R\rho|u|^{2Q}r\dd r
+\int_0^\tau\mathcal D_Q(t)\dd t\le C_T,
\qquad
\mathcal D_Q(t)=\int_0^R\rho^\alpha|u|^{2Q-2}
(u_r^2+u^2/r^2)r\dd r.
\end{equation}
Indeed, the moment satisfies
$\frac\dd{\dd t}\int\rho|u|^{2Q}r\dd r
+c_Q\mathcal D_Q\le C_T(1+\int\rho|u|^{2Q}r\dd r)$.

For a pointwise bound on $u$, set $y=\log r$. By
\eqref{eq:fix-lowerinput}, $\partial_y\rho^{\alpha-1/2}$ is bounded
in $L^2(\dd y)$. Starting at $y=\log r$, take an interval to the right
of length $c_T\rho(t,r)^{2\alpha-1}/(1+\rho(t,r)^{2\alpha-1})$.
For $c_T$ sufficiently small, Cauchy--Schwarz gives
$\rho(t,e^{y'})\ge c_\alpha\rho(t,r)$ throughout this interval.
The one-dimensional trace inequality applied to $|u|^Q$ on this interval
therefore yields
\[
\rho^\alpha|u|^{2Q}
\le C_T\mathcal D_Q(t)(1+\rho^{1-2\alpha}).
\]
At the wall, use the zero extension of $u$ and the constant extension of
$\rho$ in the logarithmic coordinate. Since $(1-\alpha)Q>3\alpha-1$
and $\rho$ is bounded above, the last estimate implies
\begin{equation}\label{eq:fix-velocitytrace}
|u(t,r)|^2\le C_T\mathcal D_Q(t)^{1/Q}\rho(t,r)^{\alpha-1}.
\end{equation}
Write $N(t)=1+\max_{0\le r\le R}\rho(t,r)^{\alpha-1}$.
By \eqref{eq:fix-outermoment}, $1+\mathcal D_Q^{1/Q}$ is integrable on
$[0,\tau]$ with integral bounded by $C_T$.

Fix $0<m_*<c_cr_c^2/2$ and let $r(t,m)$ be the trajectory with
$M(t,r(t,m))=m$. The inner trajectory $r(t,m_*)$ lies in the disk
of~\eqref{eq:fix-centralpositive}. Moreover, the density upper bound and
conservation of $M$ give
\begin{equation}\label{eq:fix-outertrajectories}
r(t,m)\ge\left(\frac{2m_*}{\sup\rho}\right)^{1/2}>0
\qquad(m\ge m_*).
\end{equation}
Integrate \eqref{eq:wz} along these trajectories. The radius stays
away from zero, $K$ is bounded, and \eqref{eq:fix-velocitytrace}
controls $u^2$. Hence, uniformly for $m\ge m_*$,
\begin{equation}\label{eq:fix-wlowerouter}
w_-(t,r(t,m))
\le C_T+C_T\int_0^t (1+\mathcal D_Q(s)^{1/Q})N(s)\dd s.
\end{equation}
By \eqref{eq:fix-centralpositive}, $\rho^{\alpha-1}$ is bounded at the
inner endpoint. Integrating $(\rho^{\alpha-1})_r=(1-\alpha)(u-w)/(2a+b)$
outwards from that endpoint
and using~\eqref{eq:fix-wlowerouter} gives
\[
N(t)\le C_T+C_T\int_{r(t,m_*)}^R|u(t,r)|\dd r
+C_T\int_0^t(1+\mathcal D_Q(s)^{1/Q})N(s)\dd s.
\]
The moment~\eqref{eq:fix-outermoment} and the fixed positive lower radius
in~\eqref{eq:fix-outertrajectories} imply
\[
\int_{r(t,m_*)}^R|u|\dd r
\le C_TN(t)^{1/[2Q(1-\alpha)]}.
\]
The exponent is strictly smaller than one. Absorbing this power and
applying Gronwall therefore bounds $N$ on $[0,\tau]$. Together with the
central bound, this establishes
\begin{equation}\label{eq:fix-densitylower}
\inf_{0\le t\le\tau,\,0\le r\le R}\rho(t,r)\ge c_T>0.
\end{equation}
The two density bounds and the spatial modulus in
\eqref{eq:fix-densitymodulus} give the uniform cutoff radius needed
near the origin in the proof of \eqref{eq:qbound}.
The continuation argument in Section~\ref{sec:higher-continuation} then
gives global existence. The uniform spatial modulus
\eqref{eq:fix-densitymodulus} and the exponential energy bound
\eqref{ext:energy-decay} imply $\|\rho-\bar\rho\|_\infty\le Ce^{-ct}$
by interpolation on the disk. In particular $\rho\ge\bar\rho/2$ for
all sufficiently large times. Combining this with
\eqref{eq:fix-densitylower} gives a positive lower bound for all time.

\subsection{The Bresch--Desjardins relation}\label{bdlt:section}
Suppose $1/2\le\alpha\le1$ and $b=2a(\alpha-1)$. Then $\theta=0$,
$2a+b=2a\alpha$ and $\beta=(\alpha-1)/\alpha$. On each positive
strong-solution interval, the effective equation reduces to
\begin{equation}\label{bdlt:relaxation}
 \Dt w+Kw=Ku.
\end{equation}
Testing in the transported mass measure and using $\rho K(w-u)=P_r$ gives
\begin{equation}\label{bdlt:relative-entropy}
 \frac{\dd}{\dd t}\int
       \left(\frac12\rho w^2+h_\gamma(\rho)\right)r\dd r
 +\int\rho K(w-u)^2r\dd r=0 .
\end{equation}
Both integrals are over the full radial domain. The term
$-\int uP_r r\dd r=\int P(u_r+u/r)r\dd r$ cancels the derivative of the
internal energy on the disk, where $u(t,R)=0$. In the plane use
$P-P(\bar\rho)$ and a spatial cutoff; the pressure difference and
velocity are in $L^2$, and the relative energy flux is integrable
on every finite strong interval.

Since $w-u=(2a+b)\rho^{\alpha-1}(\log\rho)_r$, mechanical energy
and \eqref{bdlt:relative-entropy} give
\begin{equation}\label{bdlt:root-energy}
 \sup_{t\le\tau}\int_0^R\rho^{2\alpha-1}|(\log\rho)_r|^2r\dd r\le C.
\end{equation}
By mass conservation there is $r_a\in[R/2,R]$ with $\rho(t,r_a)\le4\bar\rho/3$.
For $\alpha>1/2$, integration from $r_a$ and Cauchy--Schwarz bound the
oscillation of $\rho^{\alpha-1/2}$ by $C|\log(r/r_a)|^{1/2}$.
At $\alpha=1/2$ the same estimate holds for $\log\rho$. In both cases,
for every $\delta>0$,
\begin{equation}\label{bdlt:uniform-logcoarse}
 \rho(t,r)\le C_\delta r^{-\delta},\qquad 0<r\le R .
\end{equation}
Indeed, both a power of $1+|\log r|$ and $\exp(C\sqrt{|\log r|})$
grow more slowly than $r^{-\delta}$ at zero. Since
$\nu^{-1}=\rho^{1-\alpha}/(2a+b)$ and $K=P'(\rho)/\nu$, this also gives
\begin{equation}\label{bdlt:weighted-K}
 K+\nu^{-1}\le C_\delta r^{-\delta}\qquad(\delta>0).
\end{equation}

For $\alpha=1$ we have $\beta=0$ and $\nu=2a$. The logarithmic bound
\eqref{bdlt:uniform-logcoarse} gives $\sup r^{1/2}K\le C$. Let
$g(r)=M(r/R)^{1/2}/4$. Choose $M$ so that $g'\ge K$ and the initial
inequalities below are strict. Compare the physical velocity with
$f_\pm=\pm M+g$ and the effective velocity with $-M$ and $2M$.
At a lower contact with $u=f_-$ while $w\ge-M$, the residual of
the scalar velocity equation applied to $f_-$ is at most
\[
 2a(-M+3g/4)/r^2+(2g-M)g'-Kg<0 .
\]
At an upper contact with $u=f_+$ while $w\le2M$, its residual is at least
\[
 2a(M+3g/4)/r^2+2gg'+K(M-g)>0 .
\]
At $w=-M$ and $w=2M$, \eqref{bdlt:relaxation} gives respectively
$\Dt w>0$ and $\Dt w<0$ as long as $f_-\le u\le f_+$ and $r>0$.
At the wall $u=0$ and $w_t=-Kw$, and at the center $u=w=0$.
The four inequalities thus persist by comparison on annuli and passage
to zero inner radius. Hence $|u|+|w|\le C$ and
$|(\log\rho)_r|=|w-u|/(2a)\le C$. Integrating from a radius where
$\rho=\bar\rho$ gives positive lower and finite upper bounds uniform in time.

For $1/2\le\alpha<1$, first bound $u$ and $w$ from below.
Choose $M$ larger than the initial negative parts of $u,w$ and so that
$(1-\alpha)M\ge\sup Kr$. Compare simultaneously
$u,w\ge-M_\varepsilon(t)$, where $M_\varepsilon=M+\varepsilon t$.
At a velocity contact with $f=-M_\varepsilon$,
\[
 \alpha(f_r+\beta f/r)-K
 =\frac{(1-\alpha)M_\varepsilon}{r}-K\ge0,
 \qquad w-f\ge0.
\]
The physical comparison residual is at most
$-\varepsilon-\nu M_\varepsilon/r^2<0$. At a contact with
$w=-M_\varepsilon$, \eqref{bdlt:relaxation} gives
$\Dt w=K(M_\varepsilon+u)\ge0$. The comparison holds at the center
since $u=w=0$, and at the wall since $u=0$ and $w_t=-Kw$. Letting
$\varepsilon\downarrow0$ on each finite interval gives
\begin{equation}\label{bdlt:joint-lower}
 u(t,r),w(t,r)\ge-M.
\end{equation}

To bound $u$ from above take $k=1-\alpha\in(0,1/2]$ and $f=Br^k$.
At an upper contact the velocity residual in \eqref{eq:velocity-residual} is
\begin{equation}\label{bdlt:power-barrier}
 \mathscr P[f]=\nu(1-k^2)\frac{f}{r^2}+\alpha k\frac{f^2}{r}-Kf
       +\left(K+k^2\frac{f}{r}\right)w .
\end{equation}
Indeed, $\alpha(k+\beta)=-k^2$. By \eqref{bdlt:joint-lower}, the
negative terms are bounded in magnitude by $K(f+M)+Mk^2 f/r$.
For $B\ge1$, \eqref{bdlt:weighted-K} with $\delta=1/4$ gives
\[
 \frac{K(f+M)+Mk^2 f/r}{\nu f/r^2}
 \le C(r^{3/2}+r^{3/2-k}+r^{3/4})\qquad(0<r\le\min\{1,R\}).
\]
Thus diffusion dominates on a fixed inner disk. On the remaining
annulus the positive quadratic term dominates for large $B$.
Increase $B$ to enclose the initial velocity. The wall value is zero,
and $u=O_\tau(r)$ justifies comparison down to the center because $k<1$.
Hence $-M\le u\le Br^k$, and \eqref{bdlt:relaxation} bounds $w$
along every trajectory.

It follows that $|(\log\rho)_r|=|w-u|/\nu\le Cr^{-1/4}$.
Integration from a radius where $\rho=\bar\rho$ bounds $\log\rho$
uniformly. Now $\nu^{-1}$ is bounded, so the same identity bounds
$(\log\rho)_r$ uniformly as well. Together with the case $\alpha=1$,
this proves throughout the BD range
\begin{equation}\label{bdlt:uniform-density}
 0<m\le\rho(t,x)\le M<\infty,
 \qquad
 \norm{\nabla\log\rho(t)}_{L^\infty}
      +\norm{\rho(t)}_{W^{1,\infty}}+\norm{U(t)}_{L^\infty}
      +\norm{w(t)}_{L^\infty}\le C.
\end{equation}
This holds for every $\gamma\ge1$. The exponential energy bound
\eqref{ext:energy-decay} and the uniform Lipschitz bound for $\rho$
then give $\|\rho(t)-\bar\rho\|_\infty\le Ce^{-ct}$ by interpolation.

\subsection{Negative geometric coefficient}\label{ngall:section}
Let $1/2<\alpha\le1$ and $-a\le b<2a(\alpha-1)$. Write
\[
 z=2\alpha-1,\qquad h=-\theta.
\]
Then
\begin{equation}\label{ngall:parameters}
 0<h\le z\le1,\qquad -1\le\beta<(\alpha-1)/\alpha\le0,
 \qquad 1-\beta=(1+h)/\alpha.
\end{equation}
The equality $h=z$ is equivalent to $b=-a$ and $\beta=-1$.
Fix a positive strong-solution interval $[0,\tau]$ with $\tau\le T$. Constants
denoted by $C$ are independent of $T$ and $\tau$.
The mechanical energy identity and Lemma~\ref{ext:uniform-weighted} give
\begin{equation}\label{ngall:inputs}
 \sup_{t\le\tau}\left\{
  \int_0^R\rho u^2r\dd r+\int_0^R\rho^\gamma r\dd r
  +\rho(t,R)\right\}+\int_0^\tau \mathscr D\dd t\le C.
\end{equation}
For $\gamma=1$ the second integral is the conserved mass and the relative
potential is $h_1(s)=\kappa[s\log(s/\bar\rho)-s+\bar\rho]$.
When $h<z$, we have $1-\beta<2$,
and Lemma~\ref{ext:early-decay} gives
\begin{equation}\label{ngall:early-decay}
 \mathcal E(t)\le
 \begin{cases}
 Ce^{-ct},&\gamma>1,\\
 Ce^{-c\sqrt{1+t}},&\gamma=1 .
 \end{cases}
\end{equation}
When $h=z$, we first prove the density bounds using the mechanical
energy identity and then apply Lemma~\ref{ext:early-decay}.

\paragraph{Pressure trace.}
Let $s=\gamma+\alpha-1$ and set
\[
 H_\tau=\sup_{t\le\tau,\ 0<r\le R}\rho(t,r)^s r^{1+h}.
\]
On $[0,\tau]$ this supremum is finite by local regularity. Since $s\ge\alpha>1/2$,
\begin{equation}\label{ngall:K-trace}
 Kr^{1+h}\le C H_\tau^{(\gamma-\alpha)/s},\qquad
 \frac{\gamma-\alpha}{s}=1-\frac z s<1.
\end{equation}
Indeed, $r^{1+h}\rho^{\gamma-\alpha}
=(r^{1+h}\rho^s)^{1-z/s}r^{(1+h)z/s}$.

Choose
\[
 M=C_0\left(1+H_\tau^{(\gamma-\alpha)/s}\right),
 \qquad f(r)=-Mr^{-h},
\]
where $C_0$ depends only on the fixed parameters and the initial norms.
We compare $f$ simultaneously with $w$ and $u$ to show that
\begin{equation}\label{ngall:joint-lower}
 w(t,r)\ge-Mr^{-h},\qquad u(t,r)\ge-Mr^{-h}.
\end{equation}
At a first contact $u=f$ with $w\ge f$, the factor
$\alpha(f_r+\beta f/r)-K=-\alpha(\theta+\beta)Mr^{\theta-1}-K$ is nonnegative
by \eqref{ngall:K-trace} and $\theta+\beta<0$. Hence \eqref{eq:velocity-residual} gives
\[
 \mathscr P[f]\le
 \nu M(\theta^2-1)r^{\theta-2}
 +\theta M^2r^{2\theta-1}<0.
\]
Conversely, at a contact $r^hw=-M$ while $u\ge-Mr^{-h}$, the effective
equation gives
\[
 \Dt(r^hw)=K(M+r^hu)+h r^{h-1}u^2>0.
\]
Increasing $C_0$ makes both comparisons strict at $t=0$. At the wall
$u=0$ and $w_t=-Kw$. At the center $u,w\to0$ while $f\to-\infty$.
We may thus apply the comparisons on $[\epsilon,R]$ and let
$\epsilon\downarrow0$, which proves \eqref{ngall:joint-lower}.

To bound $H_\tau$, differentiate the weighted density power. The choice
$s=\gamma+\alpha-1$ gives
\[
 \partial_r(\rho^s r^{1+h})
 =\frac{s}{2a+b}\rho^\gamma(w-u)r^{1+h}
 +(1+h)\rho^s r^h.
\]
Integrating from $r$ to $R$, using the wall density bound and
\eqref{ngall:joint-lower}, we obtain
\[
 H_\tau\le C+CM\sup_{t\le\tau}\int_0^R\rho^\gamma r\dd r
 +C\sup_{t\le\tau}\int_0^R\rho^\gamma|u|r^{1+h}\dd r.
\]
Cauchy--Schwarz and the kinetic energy bound give
\begin{align*}
 \int_0^R\rho^\gamma|u|r^{1+h}\dd r
 &\le C\left(\int_0^R\rho^{2\gamma-1}r^{1+2h}\dd r\right)^{1/2}\\
 &\le C H_\tau^{(2\gamma-1)/(2s)}
 \left(\int_0^Rr^{-1+(1+h)z/s}\dd r\right)^{1/2}
 \le C H_\tau^{1-z/(2s)}.
\end{align*}
Both exponents in
\[
 H_\tau\le C\left(1+H_\tau^{1-z/s}
                         +H_\tau^{1-z/(2s)}\right)
\]
are strictly below one.  Thus
\begin{equation}\label{ngall:uniform-pressure-lower}
 \sup_{t\le\tau,r}\rho^s r^{1+h}\le C,\qquad
 w,u\ge-Cr^{-h},\qquad Kr^{1+h}\le C.
\end{equation}
The bound is uniform in $\tau$ for every $\gamma\ge1$.

For $\alpha=1$ and $b=-a$ we have $\nu=a$, $\theta=\beta=-1$, and
\begin{equation}\label{enp:linear-equations}
\begin{gathered}
 w=u+a(\log\rho)_r,\qquad Z=rw,\qquad
 K=\frac{\kappa\gamma}{a}\rho^{\gamma-1},\\
 \Dt u-a(u_{rr}+u_r/r-u/r^2)
       =(u_r-u/r-K)(w-u),\\
 \Dt Z+KZ=Kru+u^2,\qquad
 \mathscr D=a\int\rho(u_r-u/r)^2r\dd r.
\end{gathered}
\end{equation}
The mechanical energy identity remains valid, since $\mathscr D\ge0$.
In the plane we use the exterior material circle in \eqref{enp:exterior-trace}.

On the disk, \eqref{eq:wall} gives $W=\int_0^R\rho u r^2\dd r$
whenever $b=-a$. The continuity equation gives $(\rho^\alpha)_t=-\alpha\rho^\alpha u_r$
at the wall. Substitution in the moment identity yields, for all
$1/2\le\alpha\le1$ at $b=-a$,
\begin{equation}\label{enp:wall-density}
 \frac{\dd}{\dd t}\left(W+\frac{aR^2}{\alpha}\rho(t,R)^\alpha\right)
 =\int_0^R(\rho u^2+2P)r\dd r-R^2P(t,R),
 \qquad \sup_t\rho(t,R)\le C.
\end{equation}
The wall density bound follows from \eqref{ext:weighted-bound} at $r=R$.
Also $w_t=-Kw$ at the wall, so $|w(t,R)|\le |w(0,R)|$.

Let $R(t)$ be the fixed wall or the material boundary from
\eqref{enp:exterior-trace}. In both cases the radius stays between
two positive constants, the enclosed mass is fixed, and mechanical
energy gives, for every $b\ge-a$,
\begin{equation}\label{enp:core-inputs}
 \sup_t\int_0^{R(t)}(\rho+\rho^\gamma+\rho u^2)r\dd r\le C,
 \qquad
 \lim_{\varepsilon\downarrow0}\sup_t
      \int_0^\varepsilon\rho r\dd r=0.
\end{equation}
For $\gamma=1$, the second statement follows from the bound for
$\int\rho\log(1+\rho)r\dd r$; for $\gamma>1$ it follows from
H\"older's inequality. These statements also hold on every truncated
positive strong-solution interval, with the same constants.

With $h=z=1$ and $s=\gamma$, \eqref{ngall:uniform-pressure-lower}
gives the disk estimate below. On the material core the same proof uses
\eqref{enp:core-inputs} and \eqref{enp:exterior-trace} for its integral
and boundary terms. At $h=1$ the velocity residual at $-M/r$ is still
strictly negative because it contains $-M^2/r^3$. Thus
\begin{equation}\label{enp:pressure-trace}
 \rho^\gamma r^2\le C,\qquad Z\ge-C,\qquad ru\ge-C .
\end{equation}

We record a velocity comparison for the linear case $\alpha=1$, $b=-a$.
Suppose on a time interval that $K\sqrt r$ is bounded, $Z\ge-C$ on
the core, and $Z\ge-a/4$ on a fixed neighborhood of the center.
Then, provided the initial $C^1$ velocity and the outer velocity trace
are bounded,
\begin{equation}\label{enp:sqrt-comparison}
 |u(t,r)|\le C\sqrt r .
\end{equation}
Indeed, choose $c>0$ so that $c\sqrt{R(t)}\le1/4$ and set
$f_\pm=A(\pm\sqrt r+cr)$. The stationary residual in
\eqref{enp:linear-equations} is
\[
 \mathscr P[f_\pm]
 =\pm A(3a/4+Z/2)r^{-3/2}+(cA-K)f_\pm+KZ/r .
\]
For $A>2\sup K\sqrt r$, the coefficient of $Z$ in each signed residual
$\pm\mathscr P[f_\pm]$ is nonnegative. Near zero, insert $Z=-a/4$
as a lower bound. The leading positive term is $5aAr^{-3/2}/8$.
It dominates $K|f_\pm|+aK/(4r)$ by first choosing the radius small
($Kr^2\to0$ uniformly) and then increasing $A$.
On the remaining annulus $Z\ge-C$ and $K$ is bounded, so
$cA|f_\pm|$ dominates the terms of order $A+1$. Increasing $A$
also encloses the initial and outer velocities. Since $u=O_\tau(r)$
on each truncated strong interval, comparison on annuli gives
\eqref{enp:sqrt-comparison}.
\paragraph{Uniform upper bound.}
Suppose $0<h<z$, and write $N=1+\max_{[0,R]}\rho(t,\cdot)$.
Since $s=\gamma+\alpha-1$,
\[
 ((\rho^s)_r)_-\le C\rho^\gamma(u_++Cr^{-h})
\]
by \eqref{eq:w} and \eqref{ngall:uniform-pressure-lower}.
Splitting at $\ell=RN^{-1/2}$ and using mass conservation gives
\[
 \int_0^R\rho r^{-h}\dd r
 \le C N\ell^{1-h}+\ell^{-1-h}\int_0^R\rho r\dd r
 \le C N^{(1+h)/2}.
\]
The angular dissipation and \eqref{ngall:inputs} give
\[
 \int_0^R\rho^\gamma u_+\dd r
 \le C\mathscr D^{1/2}
       \left(\int_0^R\rho^{2\gamma-\alpha}r\dd r\right)^{1/2}
 \le C\mathscr D^{1/2}N^{(\gamma-\alpha)/2}.
\]
Integrate the one-sided derivative bound towards the wall.
The term $CN^{\gamma-1+(1+h)/2}$ can be absorbed into $N^s$
since the exponents differ by $(z-h)/2>0$. Hence
\begin{equation}\label{ngall:N-D}
 N^{R_0}\le C(1+\mathscr D),\qquad R_0=2s-(\gamma-\alpha)=\gamma+3\alpha-2 .
\end{equation}

Choose $2<p<4$ sufficiently close to two that
\begin{equation}\label{ngall:p-choice}
 \eta=\alpha-1+1/p>0,\qquad
 a_p=(p-1)(\gamma-\alpha)<R_0,\qquad
 p-1-p^2\beta^2/4>0 .
\end{equation}
At $p=2$ these inequalities follow from $\alpha>1/2$,
$R_0-(\gamma-\alpha)=2z>0$, and $|\beta|<1$.
Set
\[
 M_p=\int_0^R\rho u_+^p r\dd r,\qquad
 D_p=\int_0^R\rho^\alpha u_+^{p-2}(u_r^2+u^2/r^2)r\dd r .
\]
Testing the momentum equation with $u_+^{p-1}r$ gives the viscous
form $(p-1)X^2+p\beta XY+Y^2$ in $(X,Y)=(u_r,u/r)$.
It is positive definite by \eqref{ngall:p-choice}. The boundary
terms vanish at the center and the wall. Young's inequality gives
\begin{equation}\label{ngall:moment-first}
 M_p'+c_pD_p\le C_p\int_0^R\rho^{2\gamma-\alpha}u_+^{p-2}r\dd r .
\end{equation}
On $\{\rho\le L\}$, H\"older's inequality bounds the last integral
by $C_L(\int_0^R\rho u^2r\dd r)^{(p-2)/2}$ because $p<4$ and
$2\gamma-\alpha\ge(p-2)/2$. On $\{\rho>L\}$ Young's inequality
absorbs a small multiple of $\int\rho^\alpha u_+^p\dd r/r$ into $D_p$.
The remaining integral satisfies
\[
 \int_{\rho>L}\rho^{\gamma+a_p}r^{p-1}\dd r
 \le C N^{a_p}\int_{\rho>L}\rho^\gamma r\dd r
 \le C N^{a_p}\mathcal E(t),
\]
where $L$ is large enough that $h_\gamma(s)\ge c_Ls^\gamma$ on $s>L$,
also for $\gamma=1$. We obtain
\begin{equation}\label{ngall:moment-decaying-forcing}
 M_p'+c_pD_p\le C\mathcal E(t)^{(p-2)/2}+C \mathcal E(t)N^{a_p}.
\end{equation}
By \eqref{ngall:N-D}, $N^{a_p}\le C(1+\mathscr D^{a_p/R_0})$.
If $a_p>0$, apply H\"older's inequality with exponents $R_0/a_p$ and
$R_0/(R_0-a_p)$ to $\mathcal E\,\mathscr D^{a_p/R_0}$.
The integrals are bounded by \eqref{ngall:inputs} and
\eqref{ngall:early-decay}, independently of $\tau$. If $a_p=0$, only the
energy decay is needed. Thus \eqref{ngall:moment-decaying-forcing} gives
\begin{equation}\label{ngall:uniform-moment}
 \sup_{t\le\tau}M_p(t)+\int_0^\tau D_p\dd t\le C.
\end{equation}

Return to the derivative of $\rho^s$. Since $p>2$, H\"older's inequality gives
\[
 \int_0^R\rho^\gamma u_+\dd r
 \le M_p^{1/p}N^{\gamma-1/p}
       \left(\int_0^R r^{-1/(p-1)}\dd r\right)^{1-1/p}
 \le C M_p^{1/p}N^{s-\eta}.
\]
Absorbing the $r^{-h}$ contribution as before and dividing by $N^{s-\eta}$ gives
\begin{equation}\label{ngall:uniform-upper}
 N^{p\eta}\le C(1+M_p)\le C.
\end{equation}
In particular $K\le K_*$ and $\nu\ge\nu_*>0$ uniformly in time.

\paragraph{Positivity and density convergence.}
We first compare $u$ with powers of $r$ to prove
\begin{equation}\label{ngall:uniform-velocity}
 |u(t,r)|\le V r^{-\beta}
\end{equation}
with $V$ independent of time. For the upper barrier $Br^{-\beta}$,
\eqref{ngall:uniform-pressure-lower} and \eqref{ngall:uniform-upper}
give the following lower bound for the residual,
\[
 \nu_*B(1-\beta^2)r^{-\beta-2}
 -\beta B^2r^{-2\beta-1}-K_*Br^{-\beta}-CK_*r^{-h}.
\]
Since $2+\beta-h>0$, the positive diffusion term dominates the two
negative terms on a sufficiently small fixed disk for $B\ge1$ sufficiently
large. On the remaining annulus the positive quadratic term dominates
when $B$ is large.
For the lower barrier take $f=-Ar^{-\beta}+B_0r$ with
$B_0=K_*/[\alpha(1+\beta)]$. The preceding upper barrier implies
$w-u\ge-C_0r^{-h}$, while
$\alpha(f_r+\beta f/r)-K=K_*-K\ge0$.  At a contact its residual is
at most
\[
 -\nu_*A(1-\beta^2)r^{-\beta-2}
 -\beta AMr^{-\beta-h-1}
 +B_0Br^{-\beta}+C_0K_*r^{-h}.
\]
Here one uses the coarse lower bound $u\ge-Mr^{-h}$ in the term
$A\beta r^{-\beta-1}u$, and the upper bound in $B_0u$.
First choose a fixed $r_0>0$ so small that
$(-\beta)Mr_0^{1-h}\le\nu_*(1-\beta^2)/4$.
Then choose $A$ large enough to absorb the last two displayed terms
on $(0,r_0]$, and to have $f<-Mr^{-h}$ on $[r_0,R]$.
This excludes a first contact on the annulus and makes the residual
strictly negative on the inner disk. Both barriers enclose the initial
velocity after enlarging $A,B$ (choosing $B$ first). Since $-\beta<1$,
$u=O_\tau(r)$ justifies passage to zero inner radius.

The following argument applies on the fixed disk and on the material
core. Suppose $\rho\le C$, $w\ge-Cr^{-h}$,
$|u|\le Cr^k$ for some $h\le k<1$, and
$\int_0^{R(t)}\rho^\alpha(u_r^2+u^2/r^2)r\dd r\le C\mathscr D$.
The outer value of $w$ is bounded, and the outer velocity vanishes
on the disk or satisfies $|u(t,R(t))|^2\le C\mathscr D$ in the plane.
On the disk, $(\log\rho)_r\ge-Cr^{-h}$ gives
$\rho(t,R)\ge c_0>0$ by integration from a point where $\rho=\bar\rho$.
In the plane the lower bound is in \eqref{enp:exterior-trace}.

Test the equation for $r^h w$ in \eqref{eq:wz} by $r^h w$ in the
transported mass measure. By Young's inequality, the term containing
$Kr^hu$ is bounded by half the damping term plus $C\mathscr D$, since
$\rho K\le C\rho^\alpha$. In the source
$h\int\rho w u^2r^{2h}\dd r$, substitute
$\rho w=\rho u+(2a+b)(\rho^\alpha)_r/\alpha$.
Integration by parts gives a nonpositive term proportional to
$-\int\rho^\alpha u^2r^{2h-1}\dd r$, a mixed term bounded by
$C\mathscr D$, and an outer term bounded by $C|u(t,R(t))|^2$.
The center term vanishes. The cubic integral is also bounded by
$C\mathscr D$ since $\rho$ and $u$ are bounded. With
$\rho K(w-u)^2=\kappa\gamma(2a+b)\rho^{\gamma+\alpha-1}|(\log\rho)_r|^2$,
the damping and the angular dissipation give
\begin{equation}\label{ws:negative-pressure-dissipation}
 \int_0^\tau D_\rho(t)\dd t\le C,\qquad
 D_\rho(t)=\int_0^{R(t)}
        |\partial_r\rho^{(\gamma+\alpha-1)/2}|^2r^{1+2h}\dd r .
\end{equation}

Let $N(t)$ be the running maximum of $r^h w$ up to time $t$, with a
fixed lower cutoff $N_0\ge1$ above the initial and outer values and
$2\sup r^h|u|$. Then $|r^h w|\le CN$ and
$|(\log\rho)_r|\le CN r^{-h}$. Set
$\chi=\max\{0,1-2k/(1-h)\}$. If $\chi>0$, compare $u$ with
$\pm A N(t)^{\chi/2}r^{(1-h)/2}$ on $[0,t]$, keeping $N(t)$ fixed.
At any contact the bound $|u|\le Cr^k$ gives
$N(t)r^{1-h}\le (C/A)^{2/\chi}$. Writing
$f=A N(t)^{\chi/2}r^{(1-h)/2}$, the signed residual is bounded below
by $f/r^2$ times
\[
 \nu_*\bigl(1-(1-h)^2/4\bigr)
       -C N(t)r^{1-h}-\frac C A(N(t)r^{1-h})^{1-\chi/2}.
\]
Indeed, the convection and the $w-u$ term in
\eqref{eq:velocity-residual} are bounded using
$|w|\le CN(t)r^{-h}$, $K\le C$ and $|u|\le Cr^k$ at contact.
For large $A$ the residual is positive and the barriers enclose the
initial and outer data. Hence $r^{h-1}u^2\le CN^\chi$.
When $\chi=0$, this bound follows directly from $|u|\le Cr^k$.

At an increasing maximum attained at $r$, the effective equation gives
$N'\le-K(N-r^hu)+hr^{h-1}u^2$. The center limit and outer trace
confine these maxima to the interior, so the inequality follows along
material trajectories. If $\rho(t,r)\ge c_0/4$, then $K\ge c>0$
and $\chi<1$ exclude an increase for large $N_0$. Otherwise let
$s>r$ be the first radius with $\rho(t,s)=c_0/2$.
Integrating the logarithmic gradient from $r$ to $s$ gives
$Ns^{1-h}\ge c$. Weighted Cauchy--Schwarz from $s$ to $R(t)$ gives
$D_\rho(t)\ge cs^{2h}$. If $\chi>0$, then $k\ge h$ implies
$N^\chi\le C N s^{2h}$. If $\chi=0$, then $h-1+2k\ge0$ and
\[
 r^{h-1}u^2\le Cr^{h-1+2k}
       \le C s^{3h-1}\le C N s^{2h}.
\]
Thus $N'\le C D_\rho(t)N$, and \eqref{ws:negative-pressure-dissipation}
bounds $N$ uniformly. Integration from the outer circle yields
\begin{equation}\label{ngall:uniform-log-gradient}
 0<m\le\rho\le M,\qquad |(\log\rho)_r|\le Cr^{-h}.
\end{equation}

For $0<h<z$ the hypotheses hold with $k=-\beta$ by
\eqref{ngall:uniform-pressure-lower}, \eqref{ngall:uniform-upper},
\eqref{ngall:uniform-velocity} and positive bulk viscosity.
Since $r^{-h}\in L^q(r\dd r)$ for $2<q<2/h$, these bounds give
continuation by Section~\ref{sec:higher-continuation}. On the disk,
\eqref{ngall:uniform-log-gradient} gives a time-independent H\"older
modulus. Interpolation with the exponential $L^2$ decay from
Lemma~\ref{ext:early-decay} gives $\|\rho-\bar\rho\|_\infty\le Ce^{-ct}$.

For $b=-a$ and $1/2<\alpha<1$ the limiting case $h=z$ changes the velocity comparison. On either domain put
\[
 h=2\alpha-1\in(0,1),\qquad p=1-\alpha>0.
\]
Then $\theta=-h$ and $\beta=-1$. We use the transported equation
\eqref{eq:wz} for $r^hw$.
On the disk, \eqref{ext:weighted-bound} bounds the wall density
and $w_t=-Kw$ bounds $w$ at the wall. In the plane use
\eqref{enp:exterior-trace}. The core energy and uniform integrability
of the mass at the center follow from \eqref{enp:core-inputs}.

Take $\ell=\gamma+\alpha-1$. Applying
\eqref{ngall:K-trace}--\eqref{ngall:uniform-pressure-lower} with $s=\ell$
and $h=z$, using \eqref{enp:core-inputs} and the outer values above,
gives
\begin{equation}\label{enp:pressure-lower}
 \rho^\ell r^{1+h}\le C,\qquad
 Kr^{1+h}\le C,\qquad u,w\ge-Mr^{-h},
\end{equation}
with $M$ now independent of $\tau$.

Fix $\varepsilon$ small and a time $t$. The mass $m(t,2\varepsilon)$
tends to zero uniformly in $t$ by \eqref{enp:core-inputs}.
Averaging with $\dd r/r$ gives $r_a\in[\varepsilon,2\varepsilon]$ with
$\rho(t,r_a)r_a^2\le m(t,2\varepsilon)/\log2$. Since $1+h=2\alpha$,
\[
 \partial_r(\rho r^2)^\alpha
 =\frac{\alpha}{a}\rho(w-u)r^{1+h}+(1+h)\rho^\alpha r^h .
\]
Integrate from $r\le\varepsilon$ to $r_a$ and use
$w\ge-Mr^{-h}$. Cauchy--Schwarz and the kinetic energy give
$\int_0^{r_a}\rho|u|r^{1+h}\dd r\le C\varepsilon^h m(t,2\varepsilon)^{1/2}$.
Thus
\[
 (\rho(t,r)r^2)^\alpha\le C\bigl[
 m(t,2\varepsilon)^\alpha+m(t,2\varepsilon)
 +\varepsilon^h m(t,2\varepsilon)^{1/2}\bigr].
\]
Consequently
\begin{equation}\label{enp:center-small}
 \sup_{t\le\tau,\ r\le\varepsilon}\rho(t,r)r^2
 \longrightarrow0\qquad(\varepsilon\downarrow0),
\end{equation}
uniformly in $\tau$.

Take the upper velocity barrier $f=r^{-h}+Br$. Its residual at
$u=f$, using $w\ge-Mr^{-h}$, is bounded below by
\[
 \begin{aligned}
 &\nu(1-h^2)r^{-h-2}
 +[-h-\alpha(1+h)]r^{-2h-1}\\
 &\quad+2(1-\alpha-\alpha^2)Br^{-h}+B^2r
 -\alpha(1+h)Mr^{-2h-1}
 -K(1+M)r^{-h}-KBr .
 \end{aligned}
\]
Young's inequality for the $Br^{-h}$ and $KBr$ terms, followed by
$Kr^{1+h}\le C$, gives
\begin{equation}\label{enp:composite-residual}
 \mathscr P[f]\ge \frac12B^2r+
 \left[a(1-h^2)(\rho r^2)^{-p}-C\right]r^{-2h-1}.
\end{equation}
Since $p>0$ and $h<1$, \eqref{enp:center-small} makes the bracket
in \eqref{enp:composite-residual} positive on a fixed small disk. The small radius is chosen before
$B$. On the complementary annulus the negative part of the
bracket is bounded, and $B$ large makes the residual positive.
The initial and outer values are enclosed by increasing $B$.
At the center $f\to+\infty$, so comparison gives
\begin{equation}\label{enp:rough-velocity}
 -Mr^{-h}\le u(t,r)\le r^{-h}+Br\le Cr^{-h}.
\end{equation}

Set $v=\rho^{-p}$. From \eqref{enp:pressure-lower} and
\eqref{enp:rough-velocity},
\begin{equation}\label{enp:one-sided-power}
 v_r=\frac pa(u-w)\le Cr^{-h}.
\end{equation}
Let $C$ be the constant in \eqref{enp:one-sided-power}. By
\eqref{enp:center-small}, we can fix $0<r_1<\inf_t R(t)$ such that
$(\rho(t,r_1)r_1^2)^{-p}\ge1+C/(1-h)$ for every $t$.
Since $1-h=2p$, integration gives, for $r\le r_1$,
\[
 v(t,r)\ge v(t,r_1)-C r_1^{1-h}/(1-h)\ge r_1^{2p}.
\]
For $r\ge r_1$ the pressure bound \eqref{enp:pressure-lower} controls
$\rho$ directly. Thus $\rho$ is uniformly bounded above.

We can now recover the missing mechanical coercivity before
proving a density lower bound. Since $v\ge c>0$,
\eqref{enp:one-sided-power} implies, for $r\le s$,
\[
 v(s)\le v(r)+C\le Cv(r),\qquad
 \rho(r)^\alpha\le C\rho(s)^\alpha .
\]
On a time slice write $R=R(t)$ and $\widetilde u=u-ru(R)/R$.
More generally, if a positive weight satisfies
$\omega(r)/\omega(s)\le C(s/r)^d$ for $r\le s$ and some $0\le d<2$, then
\begin{equation}\label{enp:weighted-radial-korn}
 \int_0^R\omega(\widetilde u_r^2+\widetilde u^2/r^2)r\dd r
 \le C_d\int_0^R\omega(u_r-u/r)^2r\dd r .
\end{equation}
Indeed, $\widetilde u(r)=-r\int_r^R(u_s-u/s)\dd s/s$. In the variable
$x=\log(R/r)$ the kernel from $\sqrt\omega\,r(u_r-u/r)$ to
$\sqrt\omega\,\widetilde u$ is bounded in absolute value by
$Ce^{-(1-d/2)x}{\bf1}_{x\ge0}$. Young's inequality gives the angular
bound, and $\widetilde u_r=u_r-u/r+\widetilde u/r$ gives the derivative bound.
Apply this with $\omega=\rho^\alpha$ and $d=0$ by the preceding ratio.
The linear term vanishes on the disk. In the plane its contribution is bounded by
$C|u(R)|^2\int_0^R\rho^\alpha r\dd r\le C\mathscr D(t)$ by
\eqref{enp:exterior-trace} and the core mass bound. Adding the exterior estimate gives
\begin{equation}\label{enp:full-coercivity}
 \int \rho^\alpha(u_r^2+u^2/r^2)r\dd r
 \le C\mathscr D(t).
\end{equation}

The upper density bound gives $K\le K_*$ and $\nu\ge\nu_*>0$.
We next improve \eqref{enp:rough-velocity} to
\begin{equation}\label{enp:positive-power}
 |u(t,r)|\le Cr^\alpha\qquad(0<r\le R(t)).
\end{equation}
For the upper barrier $f=Br^\alpha$,
\[
 \begin{aligned}
 \mathscr P[f]\ge{}&
 \alpha^2B^2r^{2\alpha-1}
 +\nu_*(1-\alpha^2)Br^{\alpha-2}\\
 &-\alpha pMB r^{\alpha-h-1}
 -K_*Br^\alpha-K_*Mr^{-h}.
 \end{aligned}
\]
Since $h<1$ and $\alpha<1$, diffusion absorbs every negative
term near the center for $B\ge1$, and the positive quadratic
term handles the remaining annulus for large $B$.

For the lower barrier take $f=-Ar^\alpha$. On a fixed small disk,
for all $A\ge1$,
\[
 c_f=\alpha pAr^{\alpha-1}-K\ge0.
\]
At a lower contact the rough lower comparison gives
$Ar^{\alpha+h}\le M$, so
\[
 uf_r\le\alpha AMr^{\alpha-h-1},\qquad
 -c_f(w-u)\le\alpha pAMr^{\alpha-h-1}.
\]
Consequently
\[
 \mathscr P[f]\le
 -\nu_*A(1-\alpha^2)r^{\alpha-2}
 +\alpha(1+p)AMr^{\alpha-h-1}<0
\]
on a fixed sufficiently small disk, independently of $A$.
On its complementary annulus, choosing $A$ large places $f$
strictly below $-Mr^{-h}$, so a first contact is impossible there.
Both powers enclose the initial data since $u_0=O(r)$, and enclose
the outer traces after increasing their constants. This proves
\eqref{enp:positive-power}.

The bounds \eqref{enp:pressure-lower}, \eqref{enp:full-coercivity}
and \eqref{enp:positive-power} give the hypotheses of the argument
leading to \eqref{ngall:uniform-log-gradient}, now with $k=\alpha\ge h$.
It follows on either domain that
\begin{equation}\label{enp:subcritical-density}
 0<m\le\rho\le M,\qquad
 |\rho_r|+|(\log\rho)_r|\le Cr^{-h}.
\end{equation}

\section{The half-power endpoint}\label{pr:half-all-pressure}
Let $\alpha=1/2$, $\gamma\ge1$ and $b\ge-a$. Then
$\theta=(a+b)/(2a+b)\in[0,1)$ and $\beta=2\theta-1$.
At $\theta=0$, \eqref{bdlt:uniform-density} already gives both density
bounds. Assume $0<\theta<1$. Then $\mathscr D$ controls
$\int_0^R\sqrt\rho\,(u_r^2+u^2/r^2)r\dd r$. The weighted estimate
below bounds the oscillation of $\log\rho$, and mass conservation fixes
its value at one radius.
\subsection*{A positive weighted quadratic energy}
Take $\omega=r^{-2\theta}$. Giving the physical kinetic energy twice the
weight of the effective kinetic energy yields a sum of squares for
every $0<\theta<1$. Set
\begin{equation}\label{gh-energy}
\begin{aligned}
 E&=\int_0^R\omega
       \left[\frac{\rho}{2}(2u^2+w^2)+3\Pi\right]r\dd r,\\
 F&=\int_0^R\omega|(\log\rho)_r|^2r\dd r,\qquad
 M=\int_0^R\omega\rho u^2r\dd r,\\
 D_\omega&=\int_0^R\omega\sqrt\rho\,(u_r^2+u^2/r^2)r\dd r,\\
 D_P&=\int_0^R\omega\rho K(w-u)^2r\dd r
 =\kappa\gamma (2a+b)\int_0^R
       \omega\rho^{\gamma-1/2}|(\log\rho)_r|^2r\dd r .
\end{aligned}
\end{equation}
At $\gamma=1$ add $1+3\kappa\bar\rho e^{-1}\int_0^R\omega r\dd r$
to $E$, so that its shifted potential
part is positive. Since $\rho(w-u)^2=(2a+b)^2|(\log\rho)_r|^2$,
$F+M\le CE$.
All initial integrals are finite by \eqref{eq:centerregularity}, since
$\theta<1$ and we may choose $\chi>\theta$.

Test the physical and effective equations by $2\omega u$ and $\omega w$
in the mass measure. Substituting $\alpha=1/2$, $\beta=2\theta-1$ and
$(c_u,c_w)=(2,1)$ into \eqref{eq:weighted-viscous-form} gives
\begin{equation}\label{gh-identity}
\begin{aligned}
 E'+D_P&+a\int_0^R\omega\sqrt\rho\left[
 \frac{2}{1-\theta}\left(u_r-\frac ur\right)^2
 +4\theta\frac{u^2}{r^2}\right]r\dd r\\
 &=-3\theta\int_0^R\omega\rho u^3\dd r
   -6\theta\int_0^R\omega(P+\Pi)u\dd r .
\end{aligned}
\end{equation}
For each $\theta>0$ the second square controls $u/r$, and the first
then controls $u_r$. Thus the viscous integral bounds $c_\theta D_\omega$.
The pressure terms use $\rho\Pi'-\Pi=P$. No slip removes the wall fluxes,
and \eqref{eq:centerregularity} removes the central fluxes.

\subsection*{Velocity traces and the cubic term}
In the coordinate $x=\log(R/r)$, $F=\int_0^\infty\omega
|\partial_x\log\rho|^2\dd x$ and $\omega$ increases towards the center.
Fix $x$ and take $I=[x,x+\ell]$ with $\ell=(1+F/\omega(x))^{-1}$.
For $y\in I$, Cauchy--Schwarz gives
\[
 |\log\rho(y)-\log\rho(x)|^2\le \ell F/\omega(x)\le1 .
\]
Thus $\sqrt{\rho(y)}\ge e^{-1/2}\sqrt{\rho(x)}$ on $I$.
In this coordinate $\mathscr D$ controls
$\int_I\sqrt\rho\,(u^2+|\partial_yu|^2)\dd y$,
and $D_\omega$ controls the same integral with weight $\omega$.
Applying $|u(x)|^2\le 2\ell^{-1}\int_I u^2\dd y+
2\ell\int_I|\partial_y u|^2\dd y$ gives
\begin{equation}\label{gh-traces}
\begin{aligned}
 \sqrt\rho\,u^2&\le C \mathscr D(1+Fr^{2\theta}),\\
 \sqrt\rho\,u^2&\le C r^{2\theta}D_\omega(1+Fr^{2\theta}).
\end{aligned}
\end{equation}
Both constants are independent of the minimum density.

The first trace in \eqref{gh-traces} and Cauchy--Schwarz imply
\[
\begin{aligned}
 \int_0^R\omega\rho|u|^3\dd r
 &\le C\sqrt{\mathscr D}\left[
   \int_0^R\omega\rho^{3/4}u^2\dd r
   +\sqrt F\int_0^R\sqrt\omega\,\rho^{3/4}u^2\dd r\right]\\
 &\le C\sqrt{\mathscr D\,D_\omega M}+C\mathscr D\sqrt{FM}\\
 &\le\varepsilon D_\omega+C_\varepsilon \mathscr D(F+M).
\end{aligned}
\]
On $\{Fr^{2\theta}\le1\}$, interpolate the two estimates in
\eqref{gh-traces} with weights $1-\theta$ and $\theta$ to obtain
\[
 |u|\le C\rho^{-1/4}\mathscr D^{(1-\theta)/2}
                       D_\omega^{\theta/2}r^{\theta^2}.
\]

Write $\ell(s)=1$ for $\gamma>1$ and
$\ell(s)=1+|\log(s/\bar\rho)|$ for $\gamma=1$, so
$|P+\Pi|\le C\rho^\gamma\ell(\rho)$. For
$b_*=\theta$ or $b_*=(2-\theta)\theta=1-(1-\theta)^2<1$, integration by parts gives
\begin{equation}\label{gh-pressure-ibp}
 \int_0^R r^{-b_*}\rho^{\gamma-1/4}\ell(\rho)\dd r
 \le C+C\sqrt{D_P}
 \left(\int_0^R \rho^\gamma\ell(\rho)^2 r^{1+2\theta-2b_*}\dd r\right)^{1/2}.
\end{equation}
Indeed, $s^{\gamma-1/4}\ell(s)$ is locally Lipschitz on $(0,\infty)$
and $|s\,\partial_s(s^{\gamma-1/4}\ell(s))|\le C s^{\gamma-1/4}\ell(s)$
almost everywhere. The center term vanishes since $b_*<1$, and the wall density
is bounded. Cauchy--Schwarz with $D_P$ gives the displayed exponent
$2(\gamma-1/4)-(\gamma-1/2)=\gamma$.

For $b_*=\theta$ the last integral is bounded by $C$, and for
$b_*=(2-\theta)\theta$ by $C(1+E)^{1-\theta}$. When $\gamma>1$ these
bounds follow by interpolation between the mechanical bound for
$\int\rho^\gamma r\dd r$ and the weighted bound
$\int\omega\rho^\gamma r\dd r\le CE$. When $\gamma=1$ both integrals
are bounded by $C$. On $\rho\le1$ the density factor is bounded,
and on $\rho>1$ the coarse bound $\rho\le Cr^{-2(1-\theta)}$ leaves the
integrable weight $Cr^{-1+4\theta-2b_*}(1+|\log(R/r)|^2)$.

On $Fr^{2\theta}\le1$, the preceding velocity interpolation and
\eqref{gh-pressure-ibp} with $b_*=(2-\theta)\theta$ bound the pressure
integral by
\[
 C\mathscr D^{(1-\theta)/2}D_\omega^{\theta/2}
 \bigl[1+\sqrt{D_P}(1+E)^{(1-\theta)/2}\bigr].
\]
On $Fr^{2\theta}>1$, the first trace in \eqref{gh-traces} and
\eqref{gh-pressure-ibp} with $b_*=\theta$ give
$C\sqrt{\mathscr D F}(1+\sqrt{D_P})$. Young's inequality and $F\le CE$
therefore give, for every $\gamma\ge1$,
\begin{equation}\label{pr:half-pressure}
 \int_0^R\omega |P+\Pi|\,|u|\dd r
 \le\varepsilon(D_\omega+D_P)+C_\varepsilon\bigl[
       \mathscr D(1+E)+\mathscr D^{(1-\theta)/(2-\theta)}\bigr]+C\sqrt{\mathscr D F}.
\end{equation}
Substituting this and the cubic bound into \eqref{gh-identity} gives
\begin{equation}\label{pr:half-integrable-entropy-ode}
 E'+c_*(D_\omega+D_P)
 \le C(\mathscr D+\mathscr D^{1/2})(1+E)+C\mathscr D^{(1-\theta)/(2-\theta)}.
\end{equation}
By \eqref{pr:iso-fractional-coefficients}, each power of $\mathscr D$ on the
right is integrable in time with a bound independent of $\tau$. Applying
Gronwall's inequality first and then integrating \eqref{pr:half-integrable-entropy-ode} gives
\[
 \sup_{0\le t\le\tau}(E+F+M)(t)+\int_0^\tau(D_\omega+D_P)\dd t\le C.
\]
For any radii $r,s$,
\[
 |\log\rho(t,r)-\log\rho(t,s)|\le
 F^{1/2}\left|\int_s^r y^{2\theta-1}\dd y\right|^{1/2}
 \le C|r-s|^{\min(\theta,1/2)} .
\]
At each time conservation of mass gives a radius $s$ with
$\rho(t,s)=\bar\rho$. Thus $\bar\rho e^{-C}\le\rho(t,r)\le\bar\rho e^C$,
and the same spatial H\"older estimate holds for $\rho$.
Interpolating this estimate with \eqref{ext:energy-decay} gives
$\|\rho(t)-\bar\rho\|_\infty\le Ce^{-ct}$. The density bounds and
the spatial modulus also give the hypotheses of the gradient estimate
in Section~\ref{sec:higher-continuation}.

\section{Higher-order estimates and continuation on the disk}
\label{sec:higher-continuation}
Fix $T<\infty$ and let $[0,\tau]$, $\tau\le T$, be an interval of
existence. The density estimates, including \eqref{enp:finite-time-density}
for the linear endpoint, give
\begin{equation}\label{eq:densitybounds}
 0<m_T\le\rho(t,r)\le M_T,
\end{equation}
with $m_T,M_T$ independent of $\tau$. We first bound
$\nabla\log\rho$ in some $L^q$, $q>2$, and then estimate the
derivatives required to continue the solution.

\subsection{The density gradient}
For $\theta=0$, \eqref{bdlt:uniform-density} bounds $u$ and
$(\log\rho)_r$ throughout $1/2\le\alpha\le1$. If $-1<\theta<0$,
\eqref{ngall:uniform-log-gradient} and \eqref{enp:subcritical-density}
give $(\log\rho)_r\in L^\infty(0,\tau;L^q)$ for $2<q<2/|\theta|$,
with $u$ bounded in both cases. At $\theta=-1$, use
\eqref{enp:finite-time-density} on each short interval before continuation.
It remains to treat $\theta>0$.
By \eqref{eq:densitybounds}, \eqref{eq:radial-sup}, and the mechanical
energy identity, $\int_0^\tau\|U\|_\infty^2\dd t\le C_T$.
Since $0<\theta<1$, choose $q>2$ sufficiently close to $2$ that
\begin{equation}\label{eq:highq}
 2<q<3,\qquad \sigma=(1+\theta)q-2\in(0,2),\qquad \sigma^2<2q.
\end{equation}
We choose this relation between $q$ and $\sigma$ so that the term
in $u^2/r^{1+\theta}$ in \eqref{eq:wz} can be estimated by the
angular part of the viscous dissipation.
Testing momentum by $|u|^{2q-2}u$ with radial measure $r^{1-\sigma}\dd r$
gives the viscous form
\begin{equation}\label{eq:common-moment-form}
 \mathcal Q_{2q,\sigma}(X,Y)
 =(2a+b)[(2q-1)X^2+(2\beta q-\sigma)XY+(1-\beta\sigma)Y^2],
 \qquad X=u_r,\quad Y=u/r .
\end{equation}
Write $M_{2q,\sigma}=\int_0^R\rho|u|^{2q}r^{1-\sigma}\dd r$ and
$D^\alpha_{2q,\sigma}=\int_0^R\rho^\alpha|u|^{2q-2}
(X^2+Y^2)r^{1-\sigma}\dd r$.

When the density is constant, integration by parts eliminates the mixed term.
Indeed, $2q|u|^{2q-2}uu_r=(|u|^{2q})_r$, so
\begin{equation}\label{ext:constant-density-weighted-coercivity}
\begin{aligned}
 &\int_0^R |u|^{2q-2}\mathcal Q_{2q,\sigma}(X,Y)r^{1-\sigma}\dd r\\
 &\quad=(2a+b)\int_0^R |u|^{2q-2}
 \left[(2q-1)X^2+\left(1-\frac{\sigma^2}{2q}\right)Y^2\right]
 r^{1-\sigma}\dd r .
\end{aligned}
\end{equation}
The endpoint term $|u|^{2q}r^{-\sigma}$ vanishes by radial regularity
and the wall condition. Both coefficients on the right are positive.

For variable density, take a smooth cutoff supported on $[0,2r_c]$,
equal to one on $[0,r_c]$. On its support replace $\rho^\alpha$ by
$\rho(t,0)^\alpha$. Since $|\mathcal Q_{2q,\sigma}(X,Y)|\le C(X^2+Y^2)$,
the error is bounded by the oscillation of $\rho^\alpha$ times a constant
multiple of the local dissipation. By \eqref{eq:densitybounds} and the
spatial modulus of $\rho$, we may choose $r_c$ uniformly in $t$ so that
this error is absorbed in \eqref{ext:constant-density-weighted-coercivity}.
The derivative of the cutoff contributes at most $C_T M_{2q,\sigma}$.
On $r\ge r_c$, Young's inequality bounds the quadratic form below by
$cX^2-CY^2$. Since $Y^2\le r_c^{-2}|u|^2$ there, the full viscous
integral is bounded below by $c_T D^\alpha_{2q,\sigma}-C_T M_{2q,\sigma}$.

Young's inequality bounds the pressure term by
$\varepsilon D^\alpha_{2q,\sigma}+C_{T,\varepsilon}(1+M_{2q,\sigma})$,
since $P$ is bounded and $\int_0^R r^{1-\sigma}\dd r<\infty$.
For the transport term, Cauchy--Schwarz gives
\[
 \|U\|_{L^\infty}\int_0^R\rho|u|^{2q}r^{-\sigma}\dd r
 \le C_T\|U\|_{L^\infty} M_{2q,\sigma}^{1/2}
                    (D^\alpha_{2q,\sigma})^{1/2}.
\]
Absorbing a small multiple of the dissipation gives
\[
 M_{2q,\sigma}'+c_TD^\alpha_{2q,\sigma}
 \le C_T(1+\|U\|_{L^\infty}^2)(1+M_{2q,\sigma}).
\]
For $Z_p=\int_0^R\rho|r^{-\theta}w|^pr\dd r$, testing
\eqref{eq:wz} in the transported mass measure gives
\begin{equation}\label{eq:Zidentity}
\begin{aligned}
Z_p'+p\int_0^R\rho K|r^{-\theta}w|^p r\dd r
&=p\int_0^R\rho K|r^{-\theta}w|^{p-2}
(r^{-\theta}w)(r^{-\theta}u)r\dd r\\
&\quad-p\theta\int_0^R\rho|r^{-\theta}w|^{p-2}
(r^{-\theta}w)u^2r^{-\theta}\dd r.
\end{aligned}
\end{equation}
Take $p=q$. Young's inequality absorbs half the damping in the first
term on the right. Since $\rho$ and $K$ are bounded above and below,
the remainder is at most
\[
 C_T M_{2q,\sigma}^{1/2}
 \left(\int_0^R r^{1+\sigma-2q\theta}\dd r\right)^{1/2}.
\]
The integral is finite because $2+\sigma-2q\theta=q(1-\theta)>0$.
For the second term, H\"older's inequality leaves the velocity integral
$\int_0^R\rho |u|^{2q}r^{1-(1+\theta)q}\dd r$. Its exponent is
$-1-\sigma$, so it is bounded by $C_T D^\alpha_{2q,\sigma}$. Thus
the geometric term is at most $C_T Z_q^{1-1/q}(D^\alpha_{2q,\sigma})^{1/q}$.
Young's and Gronwall's inequalities now give
\[
 \sup_{t\le\tau}(M_{2q,\sigma}+Z_q)
      +\int_0^\tau D^\alpha_{2q,\sigma}\dd t\le C_T.
\]
These moments and $(\log\rho)_r=(2a+b)^{-1}\rho^{1-\alpha}(w-u)$
give \eqref{eq:qbound}.

\subsection{Parabolic estimates for the velocity}
We have now found, for each admissible choice of the coefficients,
an exponent $q>2$ such that
\begin{equation}\label{eq:qbound}
\sup_{t\le\tau}\left(\norm{\nabla\log\rho}_{L^q}
+\norm\rho_{W^{1,q}}+\norm U_{L^{2q}}\right)\le C_T.
\end{equation}
By \eqref{eq:densitybounds}, \eqref{eq:qbound}, and
$W^{1,q}\hookrightarrow C^{0,1-2/q}$, the diffusion coefficient $\nu$
is bounded above and away from zero and is H\"older continuous in
space, uniformly for $\tau\le T$.
Take a smooth partition of unity on $\overline\Omega$. For a cutoff
$\chi_j$ supported near $x_j$, write the equation for $\chi_j U$
with diffusion $\nu(t,x_j)$. The change of time
$s=\int_0^t\nu(z,x_j)\dd z$ reduces it to the Dirichlet heat equation.
The bounds on $\nu$ make the $L^m$ norms in $s$ and $t$ equivalent.
The localized initial data have zero trace and belong to $W^{2,q}$
since $U_0\in H^3$. We may therefore apply the constant-coefficient
estimate of \cite[Theorem~2.3]{DHP2007} on the disk. Extending the
transformed forcing by zero to a fixed interval of length
$T\sup\nu$ makes the constant independent of $\tau$. The spatial H\"older bound allows us to
choose the supports small enough to absorb
$(\nu(t,x)-\nu(t,x_j))\Delta(\chi_j U)$. The terms with derivatives
of $\chi_j$ contain at most one derivative of $U$ and are bounded by
interpolation after summation.
Thus, for every finite $m>1$,
\[
\norm{U_t}_{L^mL^q}+\norm U_{L^mW^{2,q}}
\le C_{T,m}\bigl(\norm{U_0}_{W^{2,q}}
+\norm{F_\alpha}_{L^mL^q}+\norm U_{L^mL^q}\bigr).
\]
The forcing satisfies
$\norm{F_\alpha}_{L^q}\le C_T(1+\norm{\nabla U}_{L^\infty})
\le\varepsilon\norm U_{W^{2,q}}+C_{T,\varepsilon}$.
Absorbing the first term and taking $m=2,4$ gives
\begin{equation}\label{eq:L}
L(t)=\norm{U(t)}_{W^{1,\infty}},\qquad
\int_0^\tau (L^2+L^4)\dd t\le C_T .
\end{equation}

\subsection{Estimates through order three}\label{sec:third-order-disk}
Dividing \eqref{eq:heat} by $\nu$ gives the Dirichlet elliptic equation
\begin{equation}\label{eq:ellipticidentity}
\Delta U=\nu^{-1}U_t+\nu^{-1}U\cdot\nabla U
-\frac{\alpha}{2a+b}\mathcal S(U)\nabla\log\rho
+\frac{P'}\nu\nabla\log\rho.
\end{equation}
The coefficient of $\mathcal S(U)\nabla\log\rho$ is constant.
Elliptic regularity and \eqref{eq:qbound} give
\begin{equation}\label{eq:ellipticbounds}
\begin{aligned}
\norm U_{H^2}&\le C_T(1+\norm{U_t}_{L^2}),\\
\norm U_{H^3}&\le C_T\left[\norm{U_t}_{H^1}
+(1+L(t))(1+\norm{\log\rho}_{H^2}+\norm{U_t}_{L^2})\right].
\end{aligned}
\end{equation}
Indeed, set $p=2q/(q-2)$. At order zero, \eqref{eq:qbound} bounds
the stress term by $C_T\|\nabla U\|_{L^p}$ and convection by
$C_T\|\nabla U\|_{L^{2q/(q-1)}}$. Both gradient norms are bounded
by $\varepsilon\|U\|_{H^2}+C_{T,\varepsilon}$ by interpolation.
After one derivative, $\|\nabla(\nu^{-1})U_t\|_{L^2}\le C_T\|U_t\|_{H^1}$
and $\|\nabla^2U\,\nabla\log\rho\|_{L^2}\le C_T\|\nabla^2U\|_{L^p}$.
Interpolating the last norm between $H^2$ and $H^3$ of $U$ allows
us to absorb its $H^3$ contribution. In the pressure derivative we use
$\|\nabla\log\rho\|_{L^4}^2\le C\|\nabla\log\rho\|_{L^2}
\|\nabla\log\rho\|_{H^1}$.

Differentiating \eqref{eq:heat} in time gives
\begin{equation}\label{eq:timeequation}
 U_{tt}-\nu\Delta U_t=(F_\alpha)_t+\nu_t\Delta U .
\end{equation}
Use $\nu\Delta U=U_t-F_\alpha$ and
$\nu_t/\nu=(\alpha-1)(\log\rho)_t$ to cancel the differentiated
coefficient of $\mathcal S(U)$. The source becomes
\begin{equation}\label{ext:forcing-time-derivative}
 \begin{aligned}
 (F_\alpha)_t+\nu_t\Delta U
 &= -\partial_t(U\cdot\nabla U)
   +\alpha\rho^{\alpha-1}\mathcal S(U_t)\nabla\log\rho\\
 &\quad+[\alpha\rho^{\alpha-1}\mathcal S(U)-P'I]\nabla(\log\rho)_t\\
 &\quad+(\log\rho)_t\bigl[(\alpha-1)(U_t+U\cdot\nabla U)
                 -(\gamma-\alpha)P'\nabla\log\rho\bigr].
 \end{aligned}
\end{equation}
On either domain let $L$ bound $\|U\|_{W^{1,\infty}}$ and
$\|(\log\rho)_t\|_{L^q}$ up to fixed constants. On the disk the
latter bound follows from continuity and \eqref{eq:qbound}.
Continuity also gives
\[
 \|\nabla(\log\rho)_t\|_2
 \le C\bigl(\|U\|_{H^2}+L\|\nabla\log\rho\|_{H^1}\bigr).
\]
Test \eqref{eq:timeequation} with $\nu^{-1}U_t$. On the disk $U_t$ has
zero trace. The terms quadratic in $U_t$ that contain
$(\log\rho)_t$, including the derivative of the test weight, are
bounded by $CL\|U_t\|_{L^{2q/(q-1)}}^2$. The term containing
$\mathcal S(U_t)$ is bounded by
$C\|\nabla\log\rho\|_{L^q}\|\nabla U_t\|_2\|U_t\|_{L^{2q/(q-2)}}$.
Interpolation and $q>2$ bound both by a small multiple of
$\|\nabla U_t\|_2^2$ and $C(1+L+L^2)\|U_t\|_2^2$.
For the remaining products in the last line of
\eqref{ext:forcing-time-derivative}, use
$L^q\cdot L^q\subset L^{q/2}$ and $H^1\hookrightarrow L^{q/(q-2)}$.
Estimating convection by $\|U\|_{W^{1,\infty}}$ and using the displayed
bound for $\nabla(\log\rho)_t$ yields
\begin{equation}\label{ext:acceleration-L2-energy}
 \begin{aligned}
 \frac{\dd}{\dd t}\int_\Omega\nu^{-1}|U_t|^2\dd x+c\|\nabla U_t\|_2^2
 &\le C(1+L+L^2)(\|U_t\|_2^2+\|U\|_{H^2}^2)
       +C(L+L^2)\|\nabla\log\rho\|_{H^1}^2\\
 &\quad+C(L^2+L^4)(\|U\|_{L^q}+\|\nabla\log\rho\|_{L^q})^2 .
 \end{aligned}
\end{equation}
Here $C$ depends on the density band and the bound for
$\|\nabla\log\rho\|_{L^q}$. On the disk \eqref{eq:ellipticbounds}
bounds the right side by
$C_T(1+L^4)(1+\|U_t\|_2^2+\|\log\rho\|_{H^2}^2)$.
In the twice differentiated equation for $\log\rho$, the terms with the
most derivatives on $U$ satisfy
\[
 \norm{\nabla^2\diver U}_{L^2}
 +\norm{\nabla^2U}_{L^{2q/(q-2)}}\norm{\nabla\log\rho}_{L^q}
 \le C_T\norm U_{H^3}.
\]
The remaining commutator terms are bounded in $L^2$ by
$C L(t)\norm{\log\rho}_{H^2}$. Testing the transport equation through
order two and substituting \eqref{eq:ellipticbounds} thus gives
\[
\begin{aligned}
 \frac{\mathrm d}{\mathrm dt}\norm{\log\rho}_{H^2}^2
 &\le \varepsilon\norm{\nabla U_t}_{L^2}^2\\
 &\quad+C_{T,\varepsilon}(1+L(t))
       \bigl(1+\norm{\log\rho}_{H^2}^2+\norm{U_t}_{L^2}^2\bigr).
\end{aligned}
\]
Adding the energy inequality for $U_t$ and taking $\varepsilon$ small gives
\[
E'(t)+d\norm{\nabla U_t}_{L^2}^2\le C_T(1+L(t)^4)E(t),
\quad E=1+\norm{\log\rho}_{H^2}^2
+\int_\Om\nu^{-1}|U_t|^2\dd x.
\]
Thus $\log\rho$ is uniformly bounded in $H^2$, $U_t$ in $L^2$,
and $U$ in $H^2$, with $U\in L^2H^3$.

Once $\|U\|_{H^2}+\|\nabla\log\rho\|_{H^1}$ is bounded, the continuity equation
gives $\|(\log\rho)_t\|_{H^1}\le C\|U\|_{H^2}$. In
\eqref{ext:forcing-time-derivative}, the only term requiring a second
derivative of $U_t$ satisfies
\[
 \|\mathcal S(U_t)\nabla\log\rho\|_2
 \le C\|\nabla\log\rho\|_{H^1}\|U_t\|_{H^1}^{1/2}\|U_t\|_{H^2}^{1/2}.
\]
Using $H^2\hookrightarrow L^\infty$ and $H^1\hookrightarrow L^4$
for the other products, then Young's inequality, gives on either domain
\begin{equation}\label{eq:combined-source}
 \|(F_\alpha)_t+\nu_t\Delta U\|_2
 \le\varepsilon\|U_t\|_{H^2}
 +C_\varepsilon(\|U\|_{H^2}+\|\nabla\log\rho\|_{H^1})
 (1+\|\nabla U\|_\infty+\|U_t\|_{H^1}).
\end{equation}
Testing \eqref{eq:timeequation} with $-\Delta U_t$, we absorb the
first term and use the time integrals of $L^2$ and $\|U_t\|_{H^1}^2$ to obtain
\[
\sup_t\norm{U_t}_{H^1}^2
+\int_0^\tau(\norm{U_t}_{H^2}^2+\norm{U_{tt}}_{L^2}^2)\dd t\le C_T.
\]
The second bound in \eqref{eq:ellipticbounds} now gives
$\norm U_{H^3}\le C_T(1+L)$. Interpolation with the bounded $H^2$
norm gives $L\le\varepsilon\norm U_{H^3}+C_{T,\varepsilon}$,
so $\sup_t\norm U_{H^3}\le C_T$.

The $H^4$ velocity estimate now follows directly from
\eqref{eq:ellipticidentity}. Both $\nu^{-1}$ and $P'/\nu$ are bounded
in $H^2$ by the density bounds and the $H^2$ estimate for $\log\rho$.
Since $H^2$ is an algebra and $\|U\|_{H^3}\le C_T$,
\[
 \|\mathcal S(U)\nabla\log\rho\|_{H^2}
 \le C\|U\|_{H^3}\|\nabla\log\rho\|_{H^2}.
\]
The convection term is bounded in $H^2$ by $C_T$. Elliptic regularity gives
\[
\norm U_{H^4}\le C_T(1+\norm{U_t}_{H^2}+\norm{\log\rho}_{H^3}).
\]
Write $g=\log\rho$. For $|\beta|=3$, the logarithmic continuity equation gives
\[
 z_t+U\cdot\nabla z=f_\beta,\qquad z=\partial^\beta g,\qquad
 f_\beta=-\partial^\beta\diver U-[\partial^\beta,U\cdot\nabla]g .
\]
The commutator is bounded in $L^2$ by $C\|U\|_{H^3}\|g\|_{H^3}$, using $L^\infty\cdot L^2$ for the terms with one derivative on either
factor and $L^4\cdot L^4$ for the two Hessians. Since $\|U\|_{H^3}\le C_T$,
the preceding elliptic estimate gives
$\|f_\beta\|_2\le C_T(1+\|g\|_{H^3}+\|U_t\|_{H^2})$.

Suppose $U\in L^1W^{1,\infty}$, with zero trace in the disk case, and
$z\in L^\infty L^2$ solves $z_t+U\cdot\nabla z=f\in L^1L^2$
with weak $L^2$ initial trace $z_0$. Then $z\in CL^2$ and, almost everywhere in time,
\begin{equation}\label{eq:transport-energy}
 \tfrac12\frac{\dd}{\dd t}\|z\|_2^2
 =\int_\Omega \bigl(zf+\tfrac12(\diver U)|z|^2\bigr)\dd x .
\end{equation}
On the disk extend $U,z,f$ by zero. A cutoff in the boundary layer
$\Lambda_\delta$ of width $\delta$ gives an error bounded, up to the test
function, by $C|\Lambda_\delta|^{1/2}\|z\|_{L^\infty L^2}
\|\nabla U\|_{L^1L^\infty}$, since $|U|\le\delta\|\nabla U\|_\infty$ there.
Thus the extended equation holds on $\R^2$. Spatial mollification has a
commutator tending to zero in $L^1L^2$. Along the flow $X(t,x)$ of $U$,
starting at $x$, it follows that
\[
 z(t,X(t,x))=z_0(x)+\int_0^t f(s,X(s,x))\dd s .
\]
The weak trace fixes the integration constant. The flow has bounded Jacobian
and inverse Jacobian on finite time intervals. Approximation of $z_0$ by
smooth functions now gives strong $L^2$ continuity, and differentiation
of the Jacobian gives \eqref{eq:transport-energy}.
Applying the identity to the derivatives of $\log\rho$ through order three gives
\[
 \frac{\dd}{\dd t}\|\log\rho\|_{H^3}^2
 \le C_T\bigl(1+\|\log\rho\|_{H^3}^2+\|U_t\|_{H^2}^2\bigr).
\]
Gronwall's inequality bounds $\sup_{t\le\tau}\|\log\rho\|_{H^3}$.
Squaring the elliptic bound and integrating gives
$\int_0^\tau\|U\|_{H^4}^2\dd t\le C_T$.
Since $\rho=\exp(\log\rho)$ and \eqref{eq:densitybounds} holds, we obtain
\begin{equation}\label{eq:globalbounds}
 \sup_{t\le\tau}
 \left(\norm{\rho(t)}_{H^3}^2+\norm{U(t)}_{H^3}^2+\norm{U_t(t)}_{H^1}^2\right)
 +\int_0^\tau
 \left(\norm U_{H^4}^2+\norm{U_t}_{H^2}^2+\norm{U_{tt}}_{L^2}^2\right)\dd t
 \le C_T .
\end{equation}

Suppose the maximal existence time $T_*<\infty$ and fix $T>T_*$.
By \eqref{eq:densitybounds} and \eqref{eq:globalbounds}, the $H^3$
norms and the reciprocal density are bounded uniformly for $t<T_*$.
Both sides of \eqref{eq:heat} are continuous in $H^1$ on every compact
subinterval of $[0,T_*)$. Their distributional equality therefore gives
$\nu(\rho)\Delta U+F_\alpha(\rho,U)=U_t\in H^1_0$ at every $t<T_*$.
The data at each such time satisfy \eqref{eq:compat}. By
Proposition~\ref{prop:local}, the uniform density lower bound and the
uniform $H^3$ bounds give a common lifespan $\delta>0$.
Restarting at a time $t<T_*$ with $T_*-t<\delta$ and using uniqueness
on the overlap extends the solution beyond $T_*$, a contradiction.

\section{Time-uniform estimates and exponential decay}
\label{ext:higher-section}
The finite-time estimates and Section~\ref{sec:higher-continuation}
give global existence. For all $b\ge-a$, the density bounds follow
from Sections~\ref{sec:positive-density}--\ref{pr:half-all-pressure}
and, at $\theta=-1$, from \eqref{enp:center-power} and the outer density bound.
On this density range the relative internal energy controls
$\|\rho-\bar\rho\|_2^2$. In each coefficient range the preceding
density estimates give a spatial H\"older seminorm bounded by
$C(1+t)^A$ for some fixed $A\ge0$. For $-1<\theta<0$ use
\eqref{ngall:uniform-log-gradient} or \eqref{enp:subcritical-density},
and at $\theta=-1$ use \eqref{enp:center-power}.
Interpolation with Lemma~\ref{ext:early-decay} gives
\begin{equation}\label{ext:uniform-density}
 0<m\le\rho(t,x)\le M,\qquad
 \norm{\rho(t)-\bar\rho}_{L^\infty}\le Ce^{-ct}
 \quad (t\ge0,\ x\in\overline\Om).
\end{equation}
We now derive the decay of the derivatives.
\subsection{Decay of the velocity and the density gradient}

By \eqref{ext:uniform-density}, $\nu$ and $K$ are uniformly
bounded above and below by positive constants.
In this section $C$ and $c>0$ may change from line to line and are
independent of $t$. By \eqref{eq:globalbounds} it suffices to prove the
estimates for sufficiently large $t$.

For $-1<\theta<1$, we use weighted moments of $u$ and $w$ to estimate
$\nabla\log\rho$ in $L^q$ for some $q>2$.
Choose $s>\max\{2,2/(1+\theta)\}$ sufficiently close to this maximum that
\begin{equation}\label{ext:late-moment-parameters}
\sigma=(1+\theta)s-2\in(0,2),\qquad \sigma^2<2s.
\end{equation}
At the lower endpoint $s=2/(1+\theta)$ for $\theta\le0$, we have
$\sigma=0$. At $s=2$ for $\theta>0$, we have $\sigma=2\theta<2$.
The desired inequalities therefore hold just above the respective endpoints.
Set
\[
M(t)=\int_0^R\rho|u|^{2s}r^{1-\sigma}\dd r,\qquad
D(t)=\int_0^R\rho^\alpha|u|^{2s-2}
       (u_r^2+u^2/r^2)r^{1-\sigma}\dd r.
\]
In \eqref{ext:constant-density-weighted-coercivity} with $q=s$, the
mixed term integrates out and $2a+b\ge a>0$. Thus the frozen viscous
form controls both terms in $D$, also when $\beta=-1$.
Replacing $\bar\rho^\alpha$ by $\rho^\alpha$ introduces an error
bounded by $C\|\rho^\alpha-\bar\rho^\alpha\|_\infty D$.
Thus the viscous integral is at least $d_0D$ for all large $t$.

Since $|u|^{2s-2}ur^{1-\sigma}$ vanishes at both endpoints, we may
subtract $P(\bar\rho)$ in the pressure integral. The density bounds,
$\sigma<2$, and Young's inequality give
\[
\begin{aligned}
&\left|\int_0^R [P(\rho)-P(\bar\rho)]
       \partial_r\bigl(|u|^{2s-2}u r^{1-\sigma}\bigr)\dd r\right|\\
&\qquad
\le\varepsilon D+
 C_\varepsilon\norm{P(\rho)-P(\bar\rho)}_{L^\infty}^2M^{1-1/s}
\le 2\varepsilon D+C_\varepsilon e^{-ct}.
\end{aligned}
\]
Here the angular part of $D$ and the density bounds give
$M\le C R^2D$. Apply Young's inequality to $M^{1-1/s}$ and use
the exponential decay of $\|P(\rho)-P(\bar\rho)\|_\infty$.
The weighted transport term is at most
\[
C\norm U_{L^\infty}M^{1/2}D^{1/2}
\le\varepsilon D+C_\varepsilon\norm U_{L^\infty}^2M.
\]
By \eqref{eq:radial-sup} on the positive density band and
\eqref{ext:energy-decay},
\begin{equation}\label{ext:integrable-velocity-amplitude}
 \int_t^\infty\|U(s)\|_\infty^2\dd s
 \le C\int_t^\infty\mathscr D(s)\dd s\le Ce^{-ct}.
\end{equation}

Consequently, after absorption in the weighted momentum identity,
\[
M'+d_1D\le C\norm U_{L^\infty}^2M+Ce^{-ct}.
\]
Since $D\ge cM$ and $\|U\|_{L^\infty}^2$ is integrable, an integrating
factor gives $M(t)\le Ce^{-ct}$. Integrating the inequality from $t$ onward,
using \eqref{ext:integrable-velocity-amplitude}, then gives
\begin{equation}\label{ext:late-weighted-velocity-decay}
M(t)+\int_t^\infty D(\tau)\dd\tau\le Ce^{-ct}.
\end{equation}

Recall $w=u+(2a+b)\rho^{\alpha-1}(\log\rho)_r$ and put
$Z(t)=\int_0^R\rho|r^{-\theta}w|^s r\dd r$.
Since $K\ge k_0>0$, \eqref{eq:Zidentity} with $p=s$ gives
\[
Z'+k_1Z\le C Z^{1-1/s}D^{1/s}+CD^{1/2}.
\]
For the geometric term, the identity $s(1+\theta)=2+\sigma$ gives the
same power of $r$ as in the angular part of $D$. For the term containing
$Kr^{-\theta}u$, Young's inequality absorbs part of the damping and leaves
\[
C\int_0^R\rho |u|^s r^{1-s\theta}\dd r
\le CD^{1/2}
 \left(\int_0^R\rho r^{1+s(1-\theta)}\dd r\right)^{1/2}
\le CD^{1/2}.
\]
Thus $Z'+k_2Z\le C(D+D^{1/2})$. By
\eqref{ext:late-weighted-velocity-decay} and Cauchy--Schwarz, the
integral of the right-hand side on $[t,t+1]$ is bounded by $Ce^{-ct}$.
Integrating against $e^{-k_2(t-\tau)}$ over successive unit intervals gives
\begin{equation}\label{ext:late-effective-moment-decay}
Z(t)\le Ce^{-ct}.
\end{equation}

By the choice of $s$ in \eqref{ext:late-moment-parameters}, we can
choose $q$ close enough to two that
\begin{equation}\label{ext:late-gradient-exponent}
 2<q<3,\qquad \frac1q>\frac1s+\tfrac12\max\{-\theta,0\}.
\end{equation}
Then $q<s$ and $r^\theta\in L^{sq/(s-q)}(r\dd r)$ on every bounded
core. H\"older's inequality and the lower density bound give
$\|w\|_{L^q}\le C Z^{1/s}$ for either sign of $\theta$. Since $M$
controls $\|U\|_{L^{2s}}^{2s}$ on the disk and
$(\log\rho)_r=(w-u)/\nu$, we obtain
\begin{equation}\label{ext:lower-order-decay}
\norm U_{L^{2q}}+\norm{\nabla\log\rho}_{L^q}
+\norm{\rho-\bar\rho}_{W^{1,q}}\le Ce^{-ct}.
\end{equation}
The zero-order density estimate follows from mass conservation and
Poincar\'e's inequality.

It remains to consider $\theta=-1$, that is, $\alpha=1$ and $b=-a$.
On the disk or the bounded material core in the plane, let $Z=rw$.
By \eqref{enp:center-power}, $|w|\le Cr^{-1+\eta_*}$ for some
$\eta_*>0$. Both $u$ and $K$ are bounded, and $K\ge k>0$.
In the equation $\Dt Z+KZ=Kru+u^2$, the source has $L^2(\dd m)$ norm
at most $C\|U(t)\|_2$, since the radius is bounded and $\dd m=\rho r\dd r$.
Integration along mass labels therefore gives
\[
 \|Z(t)\|_{L^2(\dd m)}
 \le e^{-kt}\|Z(0)\|_{L^2(\dd m)}
    +C\int_0^t e^{-k(t-s)}\|U(s)\|_2\dd s .
\]
This decays exponentially on the disk by \eqref{ext:energy-decay}
and tends to zero in the plane once \eqref{ws:zero-decay} holds.
Choose $2<q<\min\{3,2+\eta_*\}$. Since $Z=rw$,
\[
 |w|^q=|Z|^{q-2}|w|^2r^{2-q}\le C|Z|^{q-2}r^{2\eta_*-q}.
\]
The last weight belongs to $L^{2/(4-q)}(r\dd r)$ on the bounded core
because $q<2+\eta_*$. H\"older's inequality thus gives
$\|w\|_{L^q}^q\le C\|Z\|_{L^2(\dd m)}^{q-2}$. Both $\|w\|_{L^q}$
and $\|w\|_{L^2}$ decay exponentially on the disk and tend to zero on
the material core. Interpolation of $\|U\|_{L^2}$ with the pointwise
bound gives the corresponding decay in $L^{2q}$.
Together with \eqref{eq:w} and mass conservation, this proves
\eqref{ext:lower-order-decay} on the disk.

\subsection{Parabolic estimates on fixed time intervals}

On $[t,t+2]$, choose a smooth cutoff $\chi$ with $\chi(t)=0$,
$\chi=1$ on $[t+1,t+2]$, and $|\chi'|\le C$. By \eqref{eq:heat}, $\chi U$ solves the Dirichlet
heat equation with diffusion $\nu(\bar\rho)$ and forcing
\[
(\nu-\nu(\bar\rho))\Delta(\chi U)-U\cdot\nabla(\chi U)
+\alpha\rho^{\alpha-1}\mathcal S(\chi U)\nabla\log\rho
-\chi P'(\rho)\nabla\log\rho+\chi' U.
\]
For $m=2,4$ its $L^mL^q$ norm is at most
$Ce^{-ct}\norm{\chi U}_{L^mW^{2,q}}+Ce^{-ct}$, by
\eqref{ext:uniform-density}, \eqref{ext:lower-order-decay}, and
$W^{2,q}\hookrightarrow W^{1,\infty}$.
The constant in the Dirichlet heat estimate is independent of $t$.
For large $t$ the first term can be absorbed. On shifting the resulting
estimate on $[t+1,t+2]$ by one unit and using the finite-time
maximal-regularity estimate for bounded $t$, we obtain
\begin{equation}\label{ext:decaying-maximal-regularity}
\int_t^{t+1}
 (\norm{U_t}_{L^q}^m+\norm U_{W^{2,q}}^m)\dd\tau\le Ce^{-ct},
\qquad m=2,4.
\end{equation}
With $L$ as in \eqref{eq:L}, summation over unit intervals and
H\"older's inequality give
\begin{equation}\label{ext:integrable-parabolic-coefficients}
\int_t^\infty
 (L+L^2+L^4+\norm{U_t}_{L^2}^2)\dd\tau\le Ce^{-ct}.
\end{equation}

We record a damping estimate used also in the plane. Suppose
$y'+dy+G\le Ay+f$ with $d>0$ and $y,G,A,f\ge0$.
If $\int_T^\infty A<\infty$ and $\int_t^{t+1}f\le Ce^{-c_0t}$, then
\begin{equation}\label{ext:damped-gronwall}
 y(t)+\int_t^\infty G(s)\dd s\le Ce^{-c_1t},
 \qquad 0<c_1<\min\{d,c_0\}.
\end{equation}
If instead $\int_t^{t+1}(A+f)\dd s\to0$, then
$y(t)+\int_t^{t+1}G(s)\dd s\to0$.
Indeed, the integrating factor from $s$ to $t$ is at most $Ce^{-d(t-s)}$
in the first case and $Ce^{-d(t-s)/2}$ in the second, since eventually
$\int_s^{s+1}A\le d/2$. Splitting the convolution with $f$ into unit
intervals gives the bound for $y$. Integration of the original
inequality then gives the asserted bound for $G$.

\subsection{Second-order density decay}

Differentiate the logarithmic continuity equation. Since $U$ is radial
without swirl, $\nabla\diver U=\Delta U$ and $\nabla U$ is symmetric.
Substitution of \eqref{eq:ellipticidentity} gives
\begin{equation}\label{ext:damped-density-gradient}
\begin{aligned}
(\partial_t+U\cdot\nabla)\nabla\log\rho+K\nabla\log\rho
={}&-\nu^{-1}(U_t+U\cdot\nabla U)\\
&+[\alpha\beta(\diver U)I-\theta\nabla U]\nabla\log\rho.
\end{aligned}
\end{equation}
The no-slip condition makes the boundary integral from transport vanish.

By \eqref{ext:lower-order-decay}, the coefficients
$\|\nabla\log\rho\|_{L^q}$ and $\|U\|_{L^{2q}}$ in the elliptic
estimates tend to zero. With $K=P'/\nu$, \eqref{eq:ellipticidentity}
therefore gives, for all sufficiently large $t$,
\begin{equation}\label{ext:late-elliptic}
\begin{aligned}
\norm U_{H^2}&\le C(\norm{U_t}_{L^2}+e^{-ct}),\\
\norm U_{H^3}&\le C\bigl[\norm{U_t}_{H^1}
 +(1+L)\norm{\nabla\log\rho}_{H^1}+L(\norm{U_t}_{L^2}+e^{-ct})\bigr].
\end{aligned}
\end{equation}
For the second line, use
$\norm{\nabla^2U\,\nabla\log\rho}_{L^2}\le C\norm{\nabla\log\rho}_{L^q}\norm U_{H^3}$
and absorb it. The variable-coefficient term obeys
$\norm{\nabla(\nu^{-1})U_t}_{L^2}
\le C\norm{\nabla\log\rho}_{L^q}\norm{U_t}_{H^1}$.
The remaining pressure derivative is controlled by
$\norm{\nabla\log\rho}_{L^4}^2\le C\norm{\nabla\log\rho}_{L^2}\norm{\nabla\log\rho}_{H^1}$.

Testing \eqref{ext:damped-density-gradient} at orders zero and one,
and using \eqref{ext:late-elliptic}, yields
\begin{equation}\label{ext:density-H1-energy}
\begin{aligned}
&\frac{\dd}{\dd t}\norm{\nabla\log\rho}_{H^1}^2
 +d\norm{\nabla\log\rho}_{H^1}^2\\
&\quad\le C\norm{U_t}_{H^1}^2+Ce^{-ct}(1+L^2)
  +C(L+L^2+e^{-ct})(\norm{\nabla\log\rho}_{H^1}^2+\norm{U_t}_{L^2}^2).
\end{aligned}
\end{equation}
In the differentiated stress, $\nabla^2 U\,\nabla\log\rho$ is
bounded in $L^2$ by $Ce^{-ct}\|U\|_{H^3}$, for which we use the
second line of \eqref{ext:late-elliptic}. The transport commutator
is bounded by $CL\|\nabla\log\rho\|_{H^1}$. For the pressure commutator, use
$\nabla K=(\gamma-\alpha)K\nabla\log\rho$ and the $L^4$ bound above.

In \eqref{ext:acceleration-L2-energy}, both $\|U\|_{L^q}$ and
$\|\nabla\log\rho\|_{L^q}$ decay exponentially by
\eqref{ext:lower-order-decay} and the mechanical energy estimate.
The first line of \eqref{ext:late-elliptic} then bounds the right side by
\[
 C\|U_t\|_2^2
 +C(L+L^2)(\|\nabla\log\rho\|_{H^1}^2+\|U_t\|_2^2)
 +Ce^{-ct}(1+L^2+L^4).
\]
Multiply \eqref{ext:density-H1-energy} by a sufficiently small fixed
positive number and add the acceleration estimate.
Poincar\'e's inequality gives damping for $U_t$, and
\eqref{ext:integrable-parabolic-coefficients} bounds the time integral
of $L+L^2$. The term $\norm{U_t}_{L^2}^2$ on the right has
exponentially decaying integrals over $[t,t+1]$ by
\eqref{ext:decaying-maximal-regularity}.
Applying \eqref{ext:damped-gronwall} and \eqref{ext:late-elliptic} gives
\begin{equation}\label{ext:second-order-density-decay}
\norm{\nabla\log\rho}_{H^1}^2+\norm{U_t}_{L^2}^2+\norm U_{H^2}^2
+\int_t^\infty\norm{U_t}_{H^1}^2\dd\tau\le Ce^{-ct}.
\end{equation}

\subsection{Higher velocity and density derivatives}

By \eqref{ext:second-order-density-decay}, the factor
$\|U\|_{H^2}+\|\nabla\log\rho\|_{H^1}$ in \eqref{eq:combined-source}
is at most $Ce^{-ct}$. Testing \eqref{eq:timeequation} with $-\Delta U_t$
and absorbing the $\varepsilon\|U_t\|_{H^2}$ term into the Dirichlet
dissipation gives
\[
\frac{\dd}{\dd t}\norm{\nabla U_t}_{L^2}^2+d\norm{U_t}_{H^2}^2
\le Ce^{-ct}\norm{U_t}_{H^1}^2+Ce^{-ct}(1+L^2).
\]
Absorb the first term for large $t$, and apply
\eqref{ext:damped-gronwall}. The time-differentiated equation then
also controls $U_{tt}$. Finally \eqref{ext:late-elliptic} yields
$\norm U_{H^3}\le Ce^{-ct}(1+L)$, and $H^3\hookrightarrow W^{1,\infty}$
allows us to absorb $Ce^{-ct}L$ for large $t$.
We have proved
\begin{equation}\label{ext:higher-velocity-decay}
\norm{U_t}_{H^1}^2+\norm U_{H^3}^2
+\int_t^\infty
 (\norm{U_t}_{H^2}^2+\norm{U_{tt}}_{L^2}^2)\dd\tau\le Ce^{-ct}.
\end{equation}

It remains to control $\nabla\log\rho$ in $H^2$. The transport commutator and
geometric source in \eqref{ext:damped-density-gradient} are bounded
in $L^2$, after at most two derivatives, by
$C\norm U_{H^3}\norm{\nabla\log\rho}_{H^2}$.
The pressure commutator is bounded by
\[
\begin{aligned}
&\sum_{1\le|\zeta|\le2}
 \norm{\partial^\zeta(K\nabla\log\rho)-K\partial^\zeta\nabla\log\rho}_{L^2}\\
&\qquad\le C\bigl(\norm{\nabla\log\rho}_{H^1}+\norm{\nabla\log\rho}_{H^1}^2\bigr)
              \norm{\nabla\log\rho}_{H^2}.
\end{aligned}
\]
By the chain rule, the commutator is a sum of products of two gradients, a gradient
and a Hessian, or three gradients of $\log\rho$. For the mixed product
use $L^4\cdot L^4\subset L^2$ and for the cubic use $H^1\hookrightarrow L^6$.
The density bounds and \eqref{ext:second-order-density-decay} give
\[
\norm{\nu^{-1}}_{H^2}\le C,\qquad
\norm{\nu^{-1}U_t}_{H^2}\le C\norm{U_t}_{H^2}.
\]
The algebra property of $H^2$ also gives
$\|\nu^{-1}U\cdot\nabla U\|_{H^2}\le C\|U\|_{H^3}^2$.
By \eqref{ext:second-order-density-decay} and
\eqref{ext:higher-velocity-decay}, the coefficients of
$\|\nabla\log\rho\|_{H^2}$ in the commutator bounds tend to zero.
They can therefore be absorbed by the damping $K\ge k_0$.
Apply \eqref{eq:transport-energy} to each component of
$\partial^\zeta\nabla\log\rho$, $|\zeta|\le2$, in \eqref{ext:damped-density-gradient}.
Keeping the damping term on the left and summing gives, for all sufficiently large times,
\begin{equation}\label{ext:highest-density-damping}
\frac{\dd}{\dd t}\norm{\nabla\log\rho}_{H^2}^2+d\norm{\nabla\log\rho}_{H^2}^2
\le C\norm{U_t}_{H^2}^2+Ce^{-ct}.
\end{equation}
The integral of the right-hand side over $[t,t+1]$ is bounded by
$Ce^{-ct}$. Applying \eqref{ext:damped-gronwall} gives exponential
decay of $\nabla\log\rho$ in $H^2$. Mass conservation,
\eqref{ext:lower-order-decay} and composition with the exponential
then give $\|\rho-\bar\rho\|_{H^3}^2\le Ce^{-ct}$.

Applying the $H^2$ algebra estimate to \eqref{eq:ellipticidentity}
now gives
\[
\norm U_{H^4}
\le C(\norm{U_t}_{H^2}+\norm{\nabla\log\rho}_{H^2}+\norm U_{H^3}^2).
\]
Squaring and integrating gives $\int_t^\infty\|U\|_{H^4}^2\dd s\le Ce^{-ct}$.
The continuity equation gives
$\norm{\rho_t}_{H^2}\le C\norm\rho_{H^3}\norm U_{H^3}\le Ce^{-ct}$.
Together with \eqref{ext:higher-velocity-decay}, these estimates give
\eqref{eq:maindecay} for all sufficiently large $t$. Each absorption
time depends only on the data and fixed parameters, so
\eqref{eq:globalbounds} covers the finite initial interval.
Taking $t=0$ in the dissipation tail gives \eqref{eq:mainH3}.
The density bounds follow from \eqref{ext:uniform-density}.

The constants on the disk can be chosen uniformly for data with bounded
$H^3$ norm and a common positive lower density bound. For $\alpha=1$,
$b=-a$, the argument for fixed data used annular compactness.
The finite-time bounds through \eqref{enp:finite-time-density} are uniform on each fixed
interval for such a family: the length $\ell$ is common, the initial
$C^{1,\lambda}$ bounds and lower density give a common $r_1$ on the
first window, and a finite iteration uses the resulting local $H^3$
bounds on each subsequent window.

Given any sequence of these initial data, take a weak $H^3$ and strong
$H^{5/2+\varepsilon}$ subsequence on the disk, with $0<\varepsilon<1/2$.
The strong convergence of the second derivative traces preserves
\eqref{eq:compat}. The limit solution tends to its equilibrium in $H^3$
by the result for fixed data. On each fixed time interval the difference
estimate in Appendix~\ref{app:uniqueness}, with the initial difference
retained, and interpolation with the common $H^3$ bound give $C^1$
convergence in space uniformly in time. Since the limit solution decays
in $H^3$, for any $\varepsilon>0$ we can choose $T$ and then a tail of
the subsequence so that $\|\rho(T)-\bar\rho\|_\infty+\mathcal E(T)<\varepsilon$.
At the same time we may impose the fixed bounds
$\|U(T)\|_{C^1}\le1$, $\|w(T)\|_\infty\le1$ and
$\|Z(T)\|_\infty\le a/8$. In particular $|Z(T,r)|\le r$.

On the provisional band $\bar\rho/2\le\rho\le2\bar\rho$ for $t\ge T$,
the means are bounded above and away from zero uniformly in the family.
By \eqref{enp:deviatoric-korn} and the radial trace estimate,
$\int_T^t\|u\|_\infty^2\dd s\le C\mathcal E(T)$.
Since $K$ has fixed positive upper and lower bounds, the integrating
factor in \eqref{enp:linear-equations} gives
\[
 \|Z(t)\|_\infty\le \|Z(T)\|_\infty+
 C\bigl(\mathcal E(T)^{1/2}+\mathcal E(T)\bigr).
\]
For sufficiently small $\varepsilon$ this gives $|Z|\le a/4$.
Estimate \eqref{enp:sqrt-comparison} now gives $|u|\le C\sqrt r$
with a common constant. On the density band, conservation of enclosed
mass makes all radii of a particle comparable. Thus $|Kru+u^2|\le Cr$
along its path, where $r$ is its final radius. Using $|Z(T,r)|\le r$
and the positive lower bound for $K$ in \eqref{enp:linear-equations}
gives $|Z(t,r)|\le Cr$. Hence $|(\log\rho)_r|=|Z/r-u|/a\le C$,
with a constant depending only on the band and the fixed entry bounds.
We can therefore choose $\varepsilon$ independently of $T$ so that
interpolation with $\|\rho-\bar\rho\|_2^2\le C\mathcal E(T)$ strictly
improves the band. Continuity preserves it for all later times.
The pressure test and subsequent
high-order estimates now have common constants and an exponential
rate on this subsequence. The finite interval up to $T$ has common
bounds as well. Since this extraction applies to every sequence in
the bounded data class, it gives the asserted uniform constants.

This completes the proof of Theorem~\ref{thm:global}.

\section{The Cauchy problem in the plane}\label{ws:section}

In this section $\Omega=\R^2$ and $\bar\rho$ is the far-field density
in \eqref{ws:farfield}. Let $[0,\tau]$ be an interval on which the
strong solution has positive density. The Cauchy problem is locally
solvable by the construction in Appendix~\ref{app:localconstruction}.

\subsection{The exterior density and a material circle}

The mechanical balance with the relative potential of
Section~\ref{sec:preliminary} gives, for every $b\ge-a$,
\begin{equation}\label{ws:energy}
\begin{gathered}
 \mathcal E(t)+\int_0^t\mathscr D(s)\dd s=\mathcal E(0),\\
 \mathscr D(t)=\int_0^\infty\rho^\alpha
 [a(u_r-u/r)^2+(a+b)(u_r+u/r)^2]r\dd r
 \ge a\int_0^\infty\rho^\alpha(u_r-u/r)^2r\dd r .
\end{gathered}
\end{equation}
For positive bulk viscosity the decomposition also gives
$\int_0^\infty\rho^\alpha(u_r^2+u^2/r^2)r\dd r\le C\mathscr D(t)$.
Here the radial integrals defining $\mathcal E$ and $\mathscr D$ run
from zero to infinity. In the momentum equation we subtract the constant
$P(\bar\rho)$ from the pressure before removing the spatial cutoff.
The pressure difference and the viscous stress belong to $L^2$ on each
finite strong interval. Their products with $U$ are therefore
integrable, also when $\gamma=1$.

Let $m(t,r)=\int_0^r\rho(t,y)y\dd y$. Since
$|s-\bar\rho|\le C(\sqrt{h_\gamma(s)}+h_\gamma(s))$, we have
\begin{equation}\label{ws:mass}
 \Dt m=0,\qquad
 \left|m(t,r)-\tfrac12\bar\rho r^2\right|
 \le C\bigl(r\sqrt{\mathcal E(0)}+\mathcal E(0)\bigr).
\end{equation}
Let $S'=u(t,S)$ and $S(0)=R_0$. Apply \eqref{ws:mass} at $(t,S(t))$
and $(0,R_0)$. Since $m(t,S(t))=m(0,R_0)$, division by $S(t)+R_0$ gives
\[
 |S(t)-R_0|\le C\bigl(\sqrt{\mathcal E(0)}+\mathcal E(0)/R_0\bigr).
\]
In particular $R_0/2<S(t)<2R_0$ for all times on the solution interval
once $R_0$ is sufficiently large.

At $b=-a$ only the trace-free part of the strain contributes to
$\mathscr D$. For every $b\ge-a$, on the provisional exterior density band,
\begin{equation}\label{ws:exterior-korn}
 \int_r^\infty(u_y-u/y)^2y\dd y
 =\int_r^\infty(u_y^2+u^2/y^2)y\dd y+u(t,r)^2,
 \qquad r\ge S(t).
\end{equation}
Thus the exterior radial and angular derivatives, and the square of the
pointwise velocity, are controlled by $C\mathscr D(t)$. This uses only
the lower bound for $\rho^\alpha$ on the exterior band and the last
inequality in \eqref{ws:energy}. The following exterior argument
therefore applies for every $b\ge-a$.

We claim that for $r\ge S(t)$,
\begin{equation}\label{ws:exterior}
 |\rho(t,r)-\bar\rho|\le C r^{-1/2},\qquad
 \int_{S(t)}^\infty (w^2+u^2+|\rho_r|^2)r\dd r\le C,\qquad
 |w(t,r)|\le C .
\end{equation}
To prove this, first assume $\bar\rho/2\le\rho\le2\bar\rho$ in the
material exterior. It holds initially for large $R_0$ and gives
$0<k\le K\le C$ there. If $B=\frac12\int_{S(t)}^\infty\rho w^2r\dd r$,
\eqref{eq:weq} gives
\[
 B'+\int_{S(t)}^\infty\rho Kw^2r\dd r
 =\int_{S(t)}^\infty\rho
 \left[Kuw+\theta (u/r)w(w-u)\right]r\dd r .
\]
There is no transport flux at $S(t)$. Integration of $(u^2)_r$ from
infinity gives $\sup_{r\ge S(t)}|u|^2\le C \mathscr D(t)$ by Cauchy--Schwarz
with the weights $r$ and $1/r$. Young's inequality yields
\[
 B'+kB/2\le C\mathcal E(0)
 +C R_0^{-2}\mathscr D(t)(B+\mathcal E(0)).
\]
Since $\int_0^\tau\mathscr D\dd t\le\mathcal E(0)$, the integrating
factor gives $\sup_{0\le t\le\tau}B(t)\le C$, uniformly in the choice
of sufficiently large $R_0$. On this density range, $h_\gamma$ controls
$|\rho-\bar\rho|^2$ and $|\rho_r|\le C|w-u|$. Thus
\[
 r|\rho(t,r)-\bar\rho|^2\le
 2\left(\int_r^\infty|\rho-\bar\rho|^2 y\dd y\right)^{1/2}
  \left(\int_r^\infty|\rho_y|^2 y\dd y\right)^{1/2}\le C .
\]
Since $r\ge S(t)>R_0/2$, the last bound gives
$|\rho-\bar\rho|\le(2C/R_0)^{1/2}$. Choose $R_0$ so large that this
is less than $\bar\rho/4$ and that the initial exterior density lies
strictly inside the assumed band. Continuity then preserves the band.
Along exterior trajectories the coefficient
of $w$ on the left of \eqref{eq:weq} is $K-\theta u/r\ge k/2-C\mathscr D$.
The source is bounded by $C(\sqrt{\mathscr D}+\mathscr D)$.
Equation \eqref{ws:energy} bounds the convolution of $\mathscr D$
with $e^{-kt/2}$, and Cauchy--Schwarz gives the same bound for
$\sqrt{\mathscr D}$. This proves the estimate for $w$ in
\eqref{ws:exterior}.

To bound the velocity trace, fix $L$ with $L-4>2R_0$. On unit radial
intervals in $r\ge L-4$ and time intervals of length at most four,
\eqref{ws:energy} and \eqref{ws:exterior-korn} give a uniform
$L^\infty_tL^2_r\cap L^2_tH^1_r$ bound for $u$.
Interpolation gives $u\in L^6$ on these cylinders.
The radial velocity satisfies
\[
 u_t-\nu u_{rr}=\nu(u_r/r-u/r^2)-uu_r
            +(w-u)\{\alpha u_r+\alpha\beta u/r-K\}.
\]
Its right-hand side belongs to $L^{3/2}$ since $w$ and $K$ are
bounded there, $u\in L^6$, and $u_r\in L^2$.
Choose $L$ so large that \eqref{ws:exterior} allows us to absorb
$(\nu-\nu(\bar\rho))u_{rr}$ in the localized heat estimates for
both exponents $3/2$ and $2$~\cite{DHP2007}. The first estimate gives
$u\in W^{2,1}_{3/2}$ on a smaller cylinder. By the one-dimensional
parabolic Sobolev embedding, $u_r\in L^3$ there. Thus $u u_r\in L^2$,
and a second localization gives $u\in W^{2,1}_2\hookrightarrow L^\infty$.
For cylinders meeting $t=0$ the initial trace is controlled by
$\|U_0\|_{H^3}$. These bounds are uniform in the location of the cylinder. By \eqref{eq:w}, both $u$ and $\rho_r$ are
therefore uniformly bounded for $r\ge L$.

Increase $L$ so that particles starting in $[2L,2L+1]$ stay in
$[3L/2,3L]$ by \eqref{ws:mass}, and let $I$ be their interval of mass
labels. On the maximal interval $[0,T_*)$, the exterior density band and
\eqref{ws:exterior-korn} give, by monotone convergence,
\[
 \int_I\int_0^{T_*}|u_r(t,r(t,m))|^2\dd t\dd m
 \le C\int_0^{T_*}\mathscr D(t)\dd t\le C\mathcal E(0).
\]
Choose a label whose inner integral is at most the average over $I$.
Its radius $R(t)$ satisfies, for every $\tau<T_*$,
\begin{equation}\label{ws:interface}
 3L/2\le R(t)\le3L,\qquad R'=u(t,R),\qquad
 \int_0^\tau|u_r(t,R(t))|^2\dd t\le C .
\end{equation}
Throughout the proof we use this same circle for all estimates on the
material core. The exterior estimates also give
\begin{equation}\label{enp:exterior-trace}
 \begin{gathered}
  0<c\le\rho(t,R(t))\le C,
  |u(t,R(t))|+|w(t,R(t))|\le C,\qquad
  |u(t,R(t))|^2\le C\mathscr D(t).
 \end{gathered}
\end{equation}
Here the pointwise velocity estimate follows from
\eqref{ws:exterior-korn}. On the material core, relative energy gives
$\int_0^{R(t)}(\rho+\rho^\gamma+\rho u^2)r\dd r\le C_L$,
including $\gamma=1$.

\subsection{The first weighted density bound}
Let $H$ be the primitive used in \eqref{eq:H}, with upper limit
$R(t)$. Its equation on the moving disk is
\begin{equation}\label{ws:primitive}
\begin{aligned}
 \Dt H+r^{1-\theta}P+(1-\theta)\int_r^{R(t)}P y^{-\theta}\dd y
 ={}&\theta\int_r^{R(t)}\rho u^2 y^{-\theta}\dd y\\
 &-R(t)^{1-\theta}\mathscr F(t,R(t))
 +2aR(t)^{-\theta}(\rho^\alpha u)(t,R(t)).
\end{aligned}
\end{equation}
Indeed, $\partial_r(\Dt H)=-r^{1-\theta}P_r-\theta\rho u^2r^{-\theta}$
as in \eqref{eq:H}. At $r=R(t)$,
$H=(2a+b)R^{1-\theta}\rho^\alpha/\alpha$. By continuity and $R'=u(t,R)$,
\[
 \frac{\dd}{\dd t}H(t,R(t))
 =2aR^{-\theta}\rho^\alpha u
 -(2a+b)R^{1-\theta}\rho^\alpha(u_r+u/R),
\]
where all quantities on the right are evaluated at $R(t)$.
The pressure integration contributes $R^{1-\theta}P(t,R)$. These terms in
\eqref{ws:primitive} have absolute value at most $C(1+|u_r(t,R)|)$.

The spatial estimates in Lemma~\ref{ext:uniform-weighted} apply with
upper limit $R(t)$, since $3L/2\le R(t)\le3L$. Set
$S_H=1+\max(0,\max_{0\le r\le R(t)}H(t,r))$. At a maximum of $H$
the pressure terms in \eqref{ws:primitive} control $c_1S_H-C$.
Equation~\eqref{eq:absorption} bounds the positive part of the kinetic
tail by $C\mathscr D S_H$. Since both endpoints are material, we obtain
\[
 D^+S_H+c_1S_H\le C(1+|u_r(t,R)|)+C\mathscr D S_H .
\]
By \eqref{ws:energy} and \eqref{ws:interface} the exponentially weighted
integrals of $\mathscr D$ and $|u_r(t,R(t))|$ are bounded, so
$S_H\le C$ independently of $\tau$. Thus
\begin{equation}\label{ws:coarse}
 \rho(t,r)\le C r^{-(1-\beta)},\qquad 0<r\le R(t).
\end{equation}
For $b>0$ the stronger estimate \eqref{eq:radial-pressure-bound}
also holds on this core. In Lemma~\ref{pos:mass-pressure}, both the
mass weight $m^s$ and the endpoint $R(t)$ are transported. The test
by $m^s(r^{-\theta}w)_+$ therefore has zero boundary flux. The
half-height argument uses only the bound for
$\int_0^{R(t)}\rho^\gamma r\dd r$ and a fixed upper bound for $R(t)$.
Together with \eqref{ws:coarse} in the lower pressure range, this
gives $\sup_{t,r\le R(t)}r^2\rho^\gamma\le C$.

\subsection{A uniform modulus of continuity for \texorpdfstring{$\theta>0$}{positive theta}}
Suppose $\theta>0$. Since $R'=u(t,R)$, the transport flux cancels.
If the kinetic energy density is $\rho(c_u u^2+c_w w^2)/2$, the term
added to the right-hand side of the weighted energy identity is
\[
 \left[
 r^{1-2\theta}u
 \left\{c_u\rho^\alpha\bigl[(2a+b)u_r+b u/r\bigr]-(c_u+c_w)P\right\}
 -\frac{c_w\theta(2a+b)}{\alpha}\,r^{-2\theta}\rho^\alpha u^2
 \right]_{r=R(t)} .
\]
Here $c_w=1$ in \eqref{pht:quadratic-identity},
$(c_u,c_w)=(1+\theta,1-\theta)$ in \eqref{eq:fix-energyidentity},
and $(c_u,c_w)=(2,1)$ in \eqref{gh-identity}.
The first term contains the radial viscous stress and the pressure work. The last term comes from integrating
$-c_w\theta\int_0^{R(t)}r^{-2\theta}\rho(w-u)u^2\dd r$ by parts,
using $\rho(w-u)=(2a+b)(\rho^\alpha)_r/\alpha$.
In the inward moment the stress and pressure are paired with $u_-^{j-1}$.
Since $j>2$ and $\rho,u$ are bounded at $R(t)$, the boundary terms in
both estimates are bounded in absolute value by
\[
 C\bigl(|u(R)|+|u(R)|^2+|u(R)u_r(R)|\bigr).
\]
By \eqref{ws:interface} and Cauchy--Schwarz, the last two terms have
bounded time integrals independently of $\tau$. For the first use
$C|u(R)|\le\varepsilon+C_\varepsilon|u(R)|^2$.

For $\alpha>1/2$ use $Y=1+E+\Lambda M_J$ from \eqref{pht:joint-energy}
when $b>0$ and $Y=1+E+M_j$ from \eqref{dic:budget} when $b\le0$.
For $\alpha=1/2$ take $Y=1+E$ from \eqref{gh-energy}. When $\gamma=1$,
we shift $E$ by a constant large enough to make its potential
part nonnegative for every $R(t)\le3L$. In all three cases write
the pressure dissipation, up to its constant factor, as
\[
 S_P=\int_0^{R(t)}\rho^{\gamma+\alpha-1}
            |(\log\rho)_r|^2r^{1-2\theta}\dd r .
\]
We claim that for every $\varepsilon>0$
\begin{equation}\label{ws:small-remainder}
 Y'+c_*(D_\omega+S_P+D_j+U_j)
 \le\varepsilon+C_\varepsilon \mathscr D\,Y+g_\varepsilon(t),
 \qquad \int_0^\tau g_\varepsilon(t)\dd t\le C_\varepsilon .
\end{equation}
Here $g_\varepsilon\ge0$, $j=J$ in the positive-$b$ case, and the last two dissipations
are omitted at $\alpha=1/2$. The coefficient $c_*>0$ is independent
of $\varepsilon$. All estimates may first be made up to the end of a
positive strong-solution interval, with the same constants.

First let $\alpha>1/2$. For $b>0$ we apply the estimates in
\eqref{pht:joint-energy}, using the bound for $r^2\rho^\gamma$
proved after \eqref{ws:coarse}. For $b\le0$ we apply
\eqref{dic:budget}, using the relative energy bound on the core
and $\mathscr D^{1/2}\le1+\mathscr D$. In both cases the interior
terms on the right are at most $C+C\mathscr D Y$ after retaining a
fixed positive fraction of each dissipation on the left. When $S_P$
exceeds a sufficiently large fixed constant, a fraction of $S_P$
absorbs $C$.

At the remaining times \eqref{pht:pressure-peak} or
\eqref{dic:two-traces}, integrated from the bounded value
$\rho(t,R(t))$, gives an upper bound for the density on the core. In particular $(P+\Pi)^2/\rho^\alpha\le C$, also
at $\gamma=1$. Split the pressure drift at a small fixed $\delta>0$.
Cauchy--Schwarz gives
\[
 \int_0^{R(t)}r^{-2\theta}|P+\Pi|\,|u|\dd r
 \le C\delta^{1-\theta}D_\omega^{1/2}
       +C_\delta \mathscr D^{1/2}.
\]
Young's inequality bounds this by a prescribed small multiple of
$D_\omega$, plus $\varepsilon+C_\varepsilon \mathscr D$, by taking $\delta$
small enough. The high-density cubic term still obeys \eqref{dic:high}.
For $b>0$ the low-density cubic is absorbed by $D_\omega$ and the inward
transport $U_J$ with the same choice of $\Lambda$ as in
\eqref{pht:joint-energy}. For $b\le0$, \eqref{dic:low} applies with
$\int_0^{R(t)}\rho u^2r\dd r\le C\mathscr D$, since both $r$ and
$\rho^{1-\alpha}$ are bounded there.

In the inward moment the pressure leaves, after a fixed small multiple
of $D_j$, the integral
\[
 C\int_0^{R(t)}\rho^{2\gamma-\alpha}u_-^{j-2}r^{1+\delta_j}\dd r ,
\]
where $r^{\delta_j}$ is its original spatial weight. If $2<j<3$,
H\"older bounds this by $C \mathscr D^{(j-2)/2}$: the remaining density
power $2\gamma-\alpha j/2$ is positive, and the remaining radial power
is integrable. Since $(j-2)/2<1$, the resulting bound is at most
$\varepsilon+C_\varepsilon\mathscr D$.
If $j=J>4$, on $u_->V$ its ratio to the angular integrand of $D_J$
is at most $C\rho^{2\gamma-2\alpha}r^2/V^2$. Taking $V$ large makes
this part smaller than the fraction of $D_J$ retained on the left.
On $u_-\le V$, use $u_-^{J-2}\le V^{J-4}u_-^2$ to bound
the integral by $C_V \mathscr D$. These estimates prove
\eqref{ws:small-remainder} for $\alpha>1/2$.

For $\alpha=1/2$, apply the cubic estimate following \eqref{gh-traces}
and the pressure estimate \eqref{pr:half-pressure}. The endpoint term
in \eqref{gh-pressure-ibp} is bounded by the exterior density band.
The two remaining sublinear terms satisfy
\[
 \mathscr D^{(1-\theta)/(2-\theta)}\le\varepsilon+C_\varepsilon \mathscr D,\qquad
 (\mathscr D\,F)^{1/2}\le\varepsilon+C_\varepsilon \mathscr D\,F .
\]
This proves \eqref{ws:small-remainder} also for $\alpha=1/2$.

When $S_P$ is small, the density on the core is bounded above and
away from zero independently of $Y$.
Indeed, weighted Cauchy--Schwarz and $0<\theta<1$ give
\[
 \sup_{r\le R(t)}
 \left|\rho(t,r)^{(\gamma+\alpha-1)/2}
       -\rho(t,R(t))^{(\gamma+\alpha-1)/2}\right|
 \le C S_P(t)^{1/2}.
\]
By \eqref{ws:exterior}, $\rho(t,R(t))\ge\bar\rho/2$ and is bounded above.
Choose $s_0>0$ so that $C\sqrt{s_0}\le
\tfrac12(\bar\rho/2)^{(\gamma+\alpha-1)/2}$.
At each time with $S_P(t)\le s_0$, the preceding inequality gives
a positive lower bound and an upper bound throughout the core.
At such times $F\le C S_P$, the weighted velocity energy is at most $C D_\omega$
by its angular term, and $M_j\le C D_j$ by the same comparison.
The weighted pressure potential is bounded since $\theta<1$. Consequently
$Y\le C(1+D_\omega+S_P+D_j)$ whenever $S_P\le s_0$.
Choose $Y_*$ so large that $Y>Y_*$ and $S_P\le s_0$ imply
$D_\omega+S_P+D_j\ge1$, omitting $D_j$ when $\alpha=1/2$.
On $\{Y>Y_*\}$ this bound holds whenever $S_P\le s_0$, and $S_P$
alone controls the remainder when $S_P>s_0$. Thus the dissipation
in \eqref{ws:small-remainder} absorbs any fixed
$0<\varepsilon<c_*\min\{1,s_0\}$ on $\{Y>Y_*\}$.
On $\{Y\le Y_*\}$ the derivative of $(Y-Y_*)_+$ vanishes almost
everywhere. It follows that
\[
 \frac{\mathrm d}{\mathrm dt}(Y-Y_*)_+
 \le C_\varepsilon\mathscr D\bigl((Y-Y_*)_++Y_*\bigr)+g_\varepsilon .
\]
The integrals of $\mathscr D$ and $g_\varepsilon$ are bounded
independently of $\tau$. Gronwall's inequality gives $\sup_tY(t)\le C$,
and hence $\sup_t F(t)\le C$. For $0\le r_1<r_2\le R(t)$ apply
weighted Cauchy--Schwarz to $\log\rho$ if $\alpha=1/2$ and to
$\rho^{\alpha-1/2}/(\alpha-1/2)$ if $\alpha>1/2$. In both cases
the difference between the two values is bounded by
\[
 C F(t)^{1/2}\left(\int_{r_1}^{r_2}s^{2\theta-1}\dd s\right)^{1/2}
 \le C|r_2-r_1|^{\min\{\theta,1/2\}} .
\]
Taking $r_2=R(t)$ and using the exterior density band bounds $\rho$
above throughout the core, and away from zero if $\alpha=1/2$.
The inverse functions are Lipschitz on the resulting ranges.
Thus $\rho$ has the same spatial modulus, which extends across
$R(t)$ by the exterior gradient bound.

\subsection{Interior bounds for \texorpdfstring{$\theta\le0$}{nonpositive theta}}\label{ws:nonpositive}
\paragraph{The BD curve.}
For $\theta=0$, use the entropy identity \eqref{bdlt:relative-entropy}
in the whole plane. Restricting its bound and the mechanical energy
bound to the material core gives \eqref{bdlt:root-energy} there.
Integrating from the exterior density trace \eqref{enp:exterior-trace}
gives \eqref{bdlt:uniform-logcoarse} and \eqref{bdlt:weighted-K}
with constants uniform in $t$, since $R(t)$ stays in a fixed annulus.
For $\alpha=1$ use the two signed comparisons following
\eqref{bdlt:weighted-K}, with $3L$ in place of $R$. For
$1/2\le\alpha<1$ use the lower comparison \eqref{bdlt:joint-lower}
and the upper barrier \eqref{bdlt:power-barrier}. In both cases the
constants can be increased to enclose the exterior velocity traces.
Relaxation then bounds $w$, and integration of $(\log\rho)_r=(w-u)/\nu$
from $R(t)$ gives the same density and gradient bounds on the core.

\paragraph{Negative geometric coefficient.}
Here $-a\le b<2a(\alpha-1)$. With $h=-\theta$ and $z=2\alpha-1$,
we have $0<h\le z\le1$ and $-1\le\beta\le-h<0$.
The comparison and integration leading to
\eqref{ngall:uniform-pressure-lower} require only the bulk energy
bound and bounded outer traces, and apply also when $h=z$.
At $h=z=1$ the diffusion contribution vanishes, but the strict
negative quadratic residual is retained, as in
\eqref{enp:pressure-trace}. Thus, on the material core,
\[
 w,u\ge-Cr^{-h},\qquad
 r^{1+h}\rho^{\gamma+\alpha-1}\le C,\qquad Kr^{1+h}\le C .
\]
For the next moment estimate suppose $h<z$. The absorption in \eqref{ngall:N-D} uses
$\alpha-(1+h)/2=(z-h)/2>0$, and gives
\[
 \bigl(1+\max_{r\le R(t)}\rho\bigr)^{\gamma+3\alpha-2}
 \le C(1+\mathscr D(t)).
\]
Take $p,a_p$ from \eqref{ngall:p-choice} and use
$M_p,D_p$ from that proof on the material core.
The pressure estimate in \eqref{ngall:moment-decaying-forcing} and
the same stress trace estimate as above yield
\[
 M_p'+c_p D_p\le
 C\big(1+\mathscr D^{\,a_p/(\gamma+3\alpha-2)}+|u_r(t,R(t))|\big).
\]
Here $a_p/(\gamma+3\alpha-2)<1$. Also $M_p\le CD_p$ by angular
dissipation, since the pressure trace bounds
$r^2\rho^{1-\alpha}$ on the core. The right-hand side has bounded
integral on each unit time interval. Retaining the damping of $M_p$
and summing the exponential kernel over these intervals gives
$\sup_t M_p\le C$. The spatial variation estimate for
\eqref{ngall:uniform-upper} now bounds $\rho$ uniformly from above.
The comparisons in \eqref{ngall:uniform-velocity} give
$|u|\le Cr^{-\beta}$, using the exterior velocity at $R(t)$.

Apply the argument leading to \eqref{ngall:uniform-log-gradient}
on the material core with $k=-\beta$. The exterior band and
\eqref{eq:w} extend the core bounds for $\nabla\log\rho$ in
$L^2\cap L^q$, $2<q<2/h$, to the whole plane.
For $h=z<1$ use \eqref{enp:subcritical-density} on the material core.
The case $h=z=1$ uses
\eqref{enp:finite-time-density} for continuation and
\eqref{enp:center-power} in the large-time argument below.

\subsection{Finite-time positivity and continuation}\label{ws:finite-high}

Consider first $\alpha=1$, $b=-a$ on the disk or in the plane.
We establish finite-time positivity before making any use of large-time
convergence. Set $d=1/(16\gamma)$. There is a time length $\ell>0$,
depending only on the constants in \eqref{enp:core-inputs} and
\eqref{enp:pressure-trace}, such that the following estimates hold on
every truncated strong interval $[t_0,\tau]\subset[t_0,t_0+\ell]$.
The resulting constants may depend on the positive $H^3$ data at
$t_0$, but not on $\tau$.

Complete the square in \eqref{enp:linear-equations}:
\[
 \Dt Z+KZ=(u+Kr/2)^2-K^2r^2/4.
\]
Along any trajectory, the negative Duhamel source is bounded by
\[
 \frac14\int_{t_0}^t
 e^{-\int_s^tK}\,K(s)^2r(s)^2\dd s
 \le C\varepsilon^{2/\gamma}+C_\varepsilon\ell.
\]
Indeed, on $r(s)\le\varepsilon$ use
$Kr^2\le Cr^{2/\gamma}$ and $\int e^{-\int K}K\le1$;
on its complement use the uniform bound for $K^2r^2$ there.
Choose $\varepsilon$ and then $\ell$ so that this is less than $ad/2$.
If the current radius is small, \eqref{enp:core-inputs} makes its
conserved mass label small. Its radius at $t_0$ is then small, and
the initial value of $Z$ has negative part less than $ad/2$.
Thus, for a fixed $r_1>0$ depending on the data at $t_0$,
\begin{equation}\label{enp:short-subcritical}
 Z(t,r)\ge-ad,\qquad t_0\le t\le\tau,\quad 0<r\le r_1.
\end{equation}
The choice of $\ell$ is independent of this radius.

Since $Z=r(u+a(\log\rho)_r)\ge-ad$, we have
\[
 (r^{d/2}\sqrt\rho)_r\ge -(2a)^{-1}r^{d/2}\sqrt\rho\,u_+ .
\]
Integrate from $r$ to $r_1$. The endpoint value is bounded by
\eqref{enp:pressure-trace}, and Cauchy--Schwarz bounds the integral by
\[
 C\left(\int_0^{r_1}\rho u^2 s\dd s\right)^{1/2}
  \left(\int_0^{r_1}s^{d-1}\dd s\right)^{1/2}\le C_{t_0}
\]
by \eqref{enp:core-inputs} and $d>0$. On the remaining annulus use
\eqref{enp:pressure-trace}. Thus
\begin{equation}\label{enp:short-density-power}
 \rho\le Cr^{-d},\qquad K\le Cr^{-d(\gamma-1)},\qquad
 d(\gamma-1)<1/16.
\end{equation}

The bounds \eqref{enp:pressure-trace}, \eqref{enp:short-subcritical},
and \eqref{enp:short-density-power} verify the hypotheses of the comparison
leading to \eqref{enp:sqrt-comparison}. Hence
\begin{equation}\label{enp:short-velocity}
 |u(t,r)|\le C_{t_0}\sqrt r .
\end{equation}

Equations \eqref{enp:short-density-power} and
\eqref{enp:short-velocity} imply $|Kru+u^2|\le C_{t_0}r$.
The effective equation therefore gives $|Z|\le C_{t_0}$ on this
short interval. In particular
\[
 (\log\rho)_r\le B_{t_0}/r.
\]
By \eqref{enp:core-inputs}, choose $r_a>0$ so that $[0,r_a]$ contains
less than half the core mass at every time. Since $R(t)$ is bounded,
some $r_b(t)\in[r_a,R(t)]$ satisfies $\rho(t,r_b(t))\ge c>0$.
Integrating $(\log\rho)_r\le B_{t_0}/r$ from $r$ to $r_b(t)$ gives
$\rho(t,r)\ge c_{t_0}r^{B_{t_0}}$ for $0<r\le r_a$.
On $[r_a,R(t)]$, $|Z|+|u|\le C_{t_0}$ bounds both signs of
$(\log\rho)_r$. Enlarging $B_{t_0}$ if necessary, we obtain
\begin{equation}\label{enp:short-mass-powers}
 c_{t_0}r^{Q}\le m(t,r)\le C_{t_0}r^{2-d},\qquad
 m(t,r)=\int_0^r\rho(t,s)s\dd s,\qquad Q>2.
\end{equation}

For a trajectory ending at $(t,r)$, conservation of its mass label and
\eqref{enp:short-mass-powers} give
\[
 r(s)\le C_{t_0}r^{(2-d)/Q},\qquad t_0\le s\le t.
\]
Since $w(t_0)$ is bounded, $|Z(t_0,r)|\le C_{t_0}r$.
The Duhamel formula and $|Kru+u^2|\le C_{t_0}r$ therefore give
\begin{equation}\label{enp:finite-time-density}
 \begin{gathered}
  |Z(t,r)|\le C_{t_0}r^\eta,\qquad
  |(\log\rho)_r(t,r)|\le C_{t_0}r^{-1+\eta},\\
  0<c_{t_0}\le\rho(t,r)\le C_{t_0},\qquad \eta=(2-d)/Q\in(0,1).
 \end{gathered}
\end{equation}
The bounds extend over the whole core by changing the constants.
Choose $2<q<\min\{3,2/(1-\eta)\}$. Then
$\nabla\log\rho\in L^\infty(t_0,\tau;L^q)$ and $U$ is bounded.
In the plane the exterior estimates give the same bounds outside
the core. The higher-order continuation estimates apply to the
constant-coefficient parabolic form in
\eqref{enp:linear-equations}. The positive lower density bound and
the Korn identity \eqref{enp:deviatoric-korn} give control of the full
velocity gradient by the dissipation.

If the maximal time $T_*$ were finite, we could choose $t_0<T_*$
with $T_*-t_0<\ell$, since $\ell$ is independent of the restart data.
These bounds are then uniform for $t_0\le t\le\tau<T_*$. Together
with regularity on $[0,t_0]$, they supply the estimates needed for
continuation on the disk and in the plane.
For $\theta>0$, the bound for $Y$ already gives a uniform upper
density bound on the core. Fix $T<\infty$ and restrict to $\tau\le T$.
Completing the square in \eqref{eq:wz} on the material core now gives
$w\le Cr^\theta$, since $Kr^{1-\theta}$ is uniformly bounded there.
For $b>0$ the upper density bound and the exterior velocity bound
are uniform in time, so the same comparison as in \eqref{eq:barriers}
gives time-independent constants $A,B$ on the core. Consequently, if
$\phi_\alpha'(\rho)=(2a+b)\rho^{\alpha-2}$,
\[
 \phi_\alpha(\rho(t,r))\ge\phi_\alpha(\rho(t,R(t)))-
 \int_0^{R(t)}(C y^\theta+A y^{-\beta})\dd y .
\]
Since $\rho(t,R(t))$ is bounded away from zero, $\beta<1$, and
$\phi_\alpha(0+)=-\infty$, this gives a lower bound for $\rho$ uniform
in time. At $\alpha=1/2$ the oscillation bound for $\log\rho$ gives
the same conclusion for either sign of $b$.

The remaining case is $b\le0$, $\theta>0$, $1/2<\alpha<1$.
Here the mass comparison \eqref{eq:fix-massequation}--\eqref{eq:fix-centralpositive}
gives a fixed positive lower bound on $[0,r_c]$. In its comparison
function use $3L$ for the outer radius and
$\frac12\min\{m_0,2m(0,2L)/(3L)^2\}$ for the initial level.
Since $R(0)\ge2L$ and $R(t)\le3L$, the comparison mass at $R(t)$
is strictly less than $m(0,2L)\le m(t,R(t))$.
On the remaining annulus take $Q$ as in \eqref{eq:fix-outermoment}. The test
on the material core adds only the outer stress and pressure paired
with $R|u|^{2Q-2}u$, bounded by $C(1+|u_r(t,R)|)$. Hence
\eqref{eq:fix-outermoment} holds up to $T$ by \eqref{ws:interface}.
In deriving \eqref{eq:fix-velocitytrace}, take the short interval
in $\log r$ towards the center. Its length is unchanged and the
density comparison follows from the same weighted derivative bound.
Choose $m_*>0$ below the lower bound for the mass inside $r_c$,
as in \eqref{eq:fix-outertrajectories}. The trajectory with label $m_*$
then stays inside $r_c$, while the density upper bound keeps every
trajectory with label in $[m_*,m(0,R(0))]$ a fixed distance from zero. For
$N(t)=1+\max_{r\le R(t)}\rho^{\alpha-1}$ the integration of
\eqref{eq:wz} on these trajectories and the outward integration of
$(\rho^{\alpha-1})_r=(1-\alpha)(u-w)/(2a+b)$ yield, as in
\eqref{eq:fix-wlowerouter},
\[
 N(t)\le C_T+C_T N(t)^{1/[2Q(1-\alpha)]}
    +C_T\int_0^t(1+\mathcal D_Q(s)^{1/Q})N(s)\dd s.
\]
The displayed power is less than one and $\int_0^\tau\mathcal D_Q\dd t
\le C_T$. Absorption and Gronwall give the finite-time lower bound
independently of $\tau$.

On this band, take the weighted velocity and effective-velocity
moments of Section~\ref{sec:higher-continuation} on the material core.
The physical test adds a boundary term bounded by $C(1+|u_r(t,R)|)$.
The cutoff used to freeze the viscosities vanishes near the material
boundary. The effective-velocity test
uses the transported mass measure. On the finite-time density band,
\eqref{enp:deviatoric-korn} and \eqref{ws:energy} bound
$\int_0^\tau\|\nabla U\|_2^2\dd t$ by $C_T$ for all $b\ge-a$.
Radial regularity gives the same bound for $\int_0^\tau\|U\|_\infty^2\dd t$.
Equation~\eqref{eq:Zidentity} then gives $\|\nabla\log\rho\|_{L^q}
+\|U\|_{L^{2q}}\le C_T$ on the core for some $2<q<3$ when $\theta>0$.
When $\theta\le0$ the interior comparisons give these estimates directly.
For $-1<\theta<0$ take $q<2/|\theta|$. At $\theta=-1$ take
$q<2/(1-\eta)$ as in \eqref{enp:finite-time-density}.
On $r\ge L$ interpolation between \eqref{ws:exterior} and the pointwise
velocity bound supplies the exterior norms. Thus for each $T<\infty$ and $\tau\le T$ we have
\begin{equation}\label{ws:finite-input}
 0<m_T\le\rho\le M_T,\qquad
 \sup_{t\le\tau}\bigl(\norm{\nabla\log\rho}_{L^q}
                       +\norm U_{L^{2q}}\bigr)\le C_T
 \quad\text{for some }2<q<3 .
\end{equation}
The constants are independent of $\tau$. On this density range
$h_\gamma(\rho)$ is comparable to $|\rho-\bar\rho|^2$. Hence relative
energy also bounds $\rho-\bar\rho$ and $U$ in $L^\infty(0,\tau;L^2)$.

Use $\rho-\bar\rho$ and $\log(\rho/\bar\rho)$ in the local theory and
in the transport estimates of Section~\ref{sec:higher-continuation}.
By \eqref{ws:finite-input} and Morrey's inequality, $\nu(\rho)$ has a
uniform spatial H\"older modulus on each finite interval. Freeze the
coefficient on finitely many interior patches. By \eqref{ws:exterior},
$|\nu(\rho)-\nu(\bar\rho)|\le Cr^{-1/2}$ for large $r$, so a single
exterior patch with coefficient $\nu(\bar\rho)$ suffices. The heat
estimate and $\|F_\alpha\|_{L^q}\le\varepsilon\|U\|_{W^{2,q}}+
C_{T,\varepsilon}$ give $U\in L^mW^{2,q}$ for $m=2,4$, uniformly in $\tau$.

In \eqref{eq:ellipticbounds} and \eqref{ext:acceleration-L2-energy} take
$L=1+\|U\|_{W^{2,q}}$, which controls $\|(\log\rho)_t\|_{L^q}$.
Whole-space elliptic regularity gives
$\|V\|_{H^{k+2}}\le C(\|\Delta V\|_{H^k}+\|V\|_2)$.
Relative energy bounds $\|U\|_2$, and the time-differentiated energy
controls $\|U_t\|_2$ through $\int\nu^{-1}|U_t|^2\dd x$. Once $U$
is bounded in $H^2$, interpolation gives
$\|U\|_{W^{2,q}}\le\varepsilon\|U\|_{H^3}+C_{T,\varepsilon}$.
The disk argument now proves \eqref{eq:globalbounds} with $\rho-\bar\rho$
in the density norm. This and the lower bound in \eqref{ws:finite-input}
give a common local lifespan up to any finite maximal time.
Restarting within that lifespan proves global existence.

\paragraph{Uniform positivity and the limit of the density.}

For $\alpha=1$ and $b=-a$, we first establish compactness on each annulus
away from zero. We give the argument also for the disk, since the wall
density in \eqref{enp:wall-density} has so far only an upper bound.
Put $S=u_r-u/r$. For any fixed compact annulus $A$,
\begin{equation}\label{enp:annular-flux}
 \sup_{t\ge0}\int_t^{t+1}\int_A\rho|u|^3 r\dd r\dd s<\infty .
\end{equation}
Indeed, writing $e=\rho u^2/2+h_\gamma(\rho)$, the local energy balance is
$\partial_t(re)+\partial_r[r((e+P-P(\bar\rho))u-a\rho S u)]
=-a\rho S^2r$. Multiply by a nonnegative decreasing cutoff with
$-\chi'\ge1$ on $A$, supported derivative away from zero, and integrate
over a unit time interval. At the wall the flux vanishes.
On the support of $\chi'$, \eqref{enp:pressure-trace} bounds $u_-$ and
$\rho$ in the core. In the exterior use \eqref{ws:exterior} and the uniform
velocity bound proved before \eqref{ws:interface}. On this annulus
$|(h_\gamma+P-P(\bar\rho))u|\le\varepsilon\rho|u|^3+C_{\varepsilon,A}$:
for $\gamma>1$ the combination $h_\gamma+P-P(\bar\rho)$ is $O(\rho)$
on the bounded density range, and for $\gamma=1$ we use boundedness of
$\rho|\log(\rho/\bar\rho)|^{3/2}$. The stress flux is bounded by
Cauchy--Schwarz using kinetic energy and $\int_t^{t+1}\mathscr D$.
Absorption of the cubic term proves \eqref{enp:annular-flux}.

There is also no concentration of the time-integrated kinetic energy
at the origin. Fix $R_1<R_2$ inside the domain. Test radial momentum
with $F_\delta(r)e_r$, where $F_\delta=\min(r,\delta)$ up to $R_1$,
cut off smoothly to zero at $R_2$. The Lipschitz test is admissible
by approximation. Its moment satisfies
\[
 \frac{\dd}{\dd t}\int\rho u F_\delta r\dd r
 =\int\bigl[\rho u^2F_\delta'+P(F_\delta'+F_\delta/r)
          -a\rho S(F_\delta'-F_\delta/r)\bigr]r\dd r .
\]
The moment is $O(\delta)$ and the negative pressure and convection terms
are supported on the fixed outer annulus where $|F_\delta'|\le C\delta$.
For $\delta<y<R_1$, the square of the second factor in the stress
estimate is bounded by
\[
 \int_0^{R_2}\rho|F_\delta'-F_\delta/r|^2r\dd r
 \le C\sup_t m(t,y)+C_{R_2}\delta^2/y^2.
\]
First let $\delta\downarrow0$ and then $y\downarrow0$.
By \eqref{enp:core-inputs} the stress error tends to zero uniformly, so
$\sup_t\int_t^{t+1}\int_0^\delta\rho u^2r\dd r\dd s\to0$.
Let $m_c=m(0,R(0))$, the conserved mass in the core, and let
$I\subset[0,m_c]$ be an interval of mass labels. Split its contribution
at radius $\delta$ at each time. On the outer part, H\"older's inequality
and \eqref{enp:annular-flux} bound the integral by $C_\delta |I|^{1/3}$
over a unit time interval. First let $|I|\to0$, then $\delta\downarrow0$.
The estimate at the origin gives
\begin{equation}\label{enp:mass-kinetic}
 \psi(\varepsilon):=\sup_{\substack{t\ge0,\ I\subset[0,m_c]\\ |I|\le\varepsilon}}
 \int_t^{t+1}\int_I u^2\dd m\dd s\longrightarrow0
 \qquad(\varepsilon\downarrow0).
\end{equation}

To obtain a spatial modulus for $\rho$, write $Z=ru+a r^2\rho_m$.
For $0<\eta<1$, integration of $r^2\partial_m\min(\rho,\eta)$ over $I$
leaves an endpoint term bounded by $C\eta$, since the core radius is
bounded, and the nonpositive integral $-2\int_I\min(\rho,\eta)/\rho\dd m$.
Since $Z_-\le C$, Cauchy--Schwarz and kinetic energy give
\[
 \int_{I\cap\{\rho<\eta\}}Z_+\dd m\le C\bigl(\eta+\sqrt{|I|}\bigr).
\]
Set $k_{\eta}=\kappa\gamma \eta^{\gamma-1}/a$. From \eqref{enp:linear-equations},
$H_I=\int_I Z_+\dd m$ satisfies
\[
 H_I'+k_{\eta}H_I\le\int_I(Kr u_++u^2)\dd m+
          C k_{\eta} (\eta+\sqrt{|I|}).
\]
The pressure trace gives
$\int(Kr)^2\dd m\le C\int_0^{R(t)}\rho^\gamma r^{1+2/\gamma}\dd r\le C$.
By \eqref{enp:mass-kinetic}, the first source has unit-interval integral at most
$C\sqrt{\psi(\varepsilon)}+\psi(\varepsilon)$ for $|I|\le\varepsilon$.
Integrate with the exponential kernel at fixed $\eta$, let $\varepsilon$
tend to zero using the bounded initial $Z$, and then let $\eta$ tend to zero.
Thus $\int_I |Z|\dd m\to0$ uniformly in $t$ as $|I|\to0$.
Within the core, take $r_1\le r<s\le r_2$ on a fixed annulus, $r_1>0$.
For $I=[m(t,r),m(t,s)]$ the density bound on the annulus gives $|I|\le C(s-r)$.
The identity for $\rho_m$ then bounds $|\rho(t,s)-\rho(t,r)|$ by
$C_{r_1}(\int_I |Z|\dd m+\sqrt{|I|})$, even if a limit density vanishes.
In the plane, the uniform bound for $\rho_r$ on $r\ge L$ extends this
spatial modulus across $R(t)\ge3L/2$. Spatial mollification of
$\partial_t(r\rho)=-\partial_r(r\rho u)$ gives a uniform time modulus
on each annulus, since $r\rho u$ is locally bounded in $L^\infty_tL^1_r$.
At the wall extend $r\rho$ evenly and its flux oddly using $u(R)=0$.

For $t_j\to\infty$, pass to a subsequence along which
$\rho_j(s,r)=\rho(t_j+s,r)$ converges locally uniformly to $\rho_*(s,r)$
on annuli and bounded time intervals. Write $u_j,S_j$ for the corresponding
translates. In the plane fix $r_a=3L$. By \eqref{ws:exterior-korn},
$|u(t,r_a)|^2\le C\mathscr D(t)$, so $u_j(\cdot,r_a)\to0$ in $L^2$
on every bounded time interval. At a time $s_0$, take a compact radial
interval in $\{\rho_*(s_0,\cdot)>0\}$ joining $r_a$ to a radius $r$.
On a small time cylinder about this interval $\rho_j\ge c>0$ for large $j$.
The tail dissipation gives $S_j\to0$ in $L^2$ there. Since
\[
 \frac{u_j(s,r)}r-\frac{u_j(s,r_a)}{r_a}
 =-\int_r^{r_a} S_j(s,y)\,\frac{\dd y}{y},
\]
Cauchy--Schwarz gives $u_j\to0$ in $L^2_sL^\infty_r$ on the cylinder.
Thus momentum, convective flux and viscous stress all tend to zero in
distributions. The two conservation laws yield $\partial_rP(\rho_*)=0$
and $\partial_s\rho_*=0$. The exterior bound
$|\rho-\bar\rho|\le Cr^{-1/2}$ identifies the constant as $\bar\rho$.
Continuity extends this positive component to every $r>0$.

On the disk, Jensen's inequality in \eqref{enp:wall-density} gives
$\int_t^{t+T}P(\rho(s,R))\dd s\ge T P(\bar\rho)-C$.
On a sufficiently long fixed time window the limiting wall value is
therefore positive at some interior time. Wherever it is positive, the
same cylinder argument with $u(R)=0$ shows that the wall value is
locally constant in time. Its positive time component extends to both
ends of the window by continuity. At the middle time the positive
radial component then extends from the wall to every positive radius.
The mass bound in \eqref{enp:core-inputs} permits passage to the total
mass and identifies the constant as $\bar\rho$ on the disk as well.
Every sequence has the same annular limit. For fixed $r>0$, split the
enclosed mass at $0<\delta<r$. The integral over $[\delta,r]$ tends to
$\bar\rho(r^2-\delta^2)/2$, while the mass in $[0,\delta]$ tends to zero
uniformly in time by \eqref{enp:core-inputs}. Letting $\delta\downarrow0$
and using finite-time positivity gives, on both domains,
\begin{equation}\label{enp:annular-mass}
 m(t,r)\to\bar\rho r^2/2,\qquad \inf_{t\ge0}m(t,r)>0\qquad(r>0).
\end{equation}
Consequently $\sup_{s\ge0}r(s,m)\to0$ as $m\downarrow0$.

Along a particle, $Kr^2\le Cr^{2/\gamma}$ and
\eqref{enp:linear-equations} give
\[
 Z_-(t,m)\le Z_-(0,m)+\tfrac14\sup_{s\le t}K(s,m)r(s,m)^2 .
\]
By \eqref{enp:annular-mass} the right side tends uniformly to zero
with $m$, and \eqref{enp:core-inputs} gives the same conclusion as
the current radius tends to zero.

Let $R_*=R$ on the disk and $R_*=3L$ in the plane, so $R(t)\le R_*$.
The outer density is bounded below by $c_0>0$ by the exterior estimate
in the plane and by annular convergence and finite-time positivity
on the disk. For each $d>0$, the small-label estimate gives
$Z\ge-ad$ near zero. Integrating $(r^{d/2}\sqrt\rho)_r$ as in
\eqref{enp:short-density-power} gives $\rho\le C_d r^{-d}$.
Fix $0<d<1/4$ with $d(\gamma-1)<1/2$. Then $K\sqrt r$ is bounded
and $Z\ge-a/4$ near zero. By \eqref{enp:pressure-trace} and
\eqref{enp:sqrt-comparison}, $|u|\le C\sqrt r$ throughout the core.
The local bound $Z\ge-ad$ and the global bound $Z\ge-C$ imply
$(\log\rho)_r=(Z/r-u)/a\ge-d/r-C$. Integrating gives
\begin{equation}\label{enp:linear-ratio}
 \rho(t,r)/\rho(t,s)\le C(s/r)^d,\qquad 0<r<s\le R(t).
\end{equation}
The radial identity $(u/r)_r=(u_r-u/r)/r$ now gives
\begin{equation}\label{enp:linear-korn}
 \int_0^{R(t)}\rho(u_r^2+u^2/r^2)r\dd r
 \le C\mathscr D(t).
\end{equation}
Apply \eqref{enp:weighted-radial-korn} with $\omega=\rho$ and the
exponent $d<2$ in \eqref{enp:linear-ratio}. By the core mass bound,
the linear part contributes at most $C|u(t,R(t))|^2$. This is zero
on the disk and bounded by $C\mathscr D(t)$ in the plane by \eqref{enp:exterior-trace}.

Using \eqref{enp:linear-korn}, test \eqref{enp:linear-equations} for $Z$ by $Z$ in $\dd m$.
Young's inequality leaves $C\int K r^2 u^2\dd m\le C\mathscr D(t)$ since
$Kr^4$ is bounded. In $\int Zu^2\dd m$, the derivative term is
$a\int_0^{R(t)}\rho_r u^2r^2\dd r$. Integration by parts gives a
boundary term $a\rho(R(t))u(R(t))^2R(t)^2$ and the two integrals of
$-2a\rho uu_r r^2$ and $-2a\rho u^2r$. All are bounded in absolute
value by $C\mathscr D(t)$, as is $\int ru^3\dd m$.
Consequently
\begin{equation}\label{enp:linear-Zenergy}
 \sup_t\int Z^2\dd m+\int_0^\infty D_\rho(t)\dd t\le C,\qquad
 D_\rho(t)=\int_0^{R(t)}|(\rho^{\gamma/2})_r|^2r^3\dd r .
\end{equation}
Indeed, $K(Z-ru)^2\le2KZ^2+2Kr^2u^2$, and
$K(Z-ru)^2\dd m=(4a\kappa/\gamma)|(\rho^{\gamma/2})_r|^2r^3\dd r$.

First take $\gamma>1$ and set $p=\gamma-1$. Integrating
$(r^2\rho)_m=2+(Z-ru)/a$ from the center gives
$r^2\rho\le C\sqrt m$ by Cauchy--Schwarz, \eqref{enp:linear-Zenergy}
and the kinetic energy bound. Choose $0<\vartheta<\min\{1,p\}$ and
use $\rho\le C_d r^{-d}$ with $d$ small enough to obtain
\[
 Kr^2=C(r^2\rho)^\vartheta r^{2(1-\vartheta)}
          \rho^{p-\vartheta}\le C m^{\vartheta/2}.
\]
Completing the square in \eqref{eq:wz} gives
$\dot Z+KZ\ge-K^2r^2/4\ge-CKm^{\vartheta/2}$. Initially
$|Z|\le C\sqrt m$ by positivity and regularity. Integration at fixed
$m$ and $m\le C_d r^{2-d}$ show that
\begin{equation}\label{enp:negative-slope}
 ((\log\rho)_r)_-\le Cr^{-1+\eta_0}
\end{equation}
for some $\eta_0>0$. Integrating to the outer circle bounds $\rho$
uniformly from above, hence also bounds $P'(\rho)$.

For $v=u/r$ the radial momentum equation reads
\[
 v_t=a(v_{rr}+3v_r/r)+\bigl(a(\log\rho)_r-u\bigr)v_r-v^2-P'(\rho)(\log\rho)_r/r .
\]
Compare first with $A-r^e$, where $0<e<\min(\eta_0,1/2)$.
The diffusion term is $-ae(e+2)r^{e-2}$. The coefficient of $(\log\rho)_r$ is
negative, so \eqref{enp:negative-slope} bounds its positive contribution
by $C(r^{\eta_0-2}+r^{e+\eta_0-2})$. The term involving $u$ is lower
order by $|u|\le C\sqrt r$. A first upper contact at zero would force
$v(t,r)\le v(t,0)-r^e$ nearby. This contradicts
$v(t,r)=\int_0^1 u_r(t,sr)\dd s=v(t,0)+O_\tau(r^\lambda)$, since
$U\in C^{1,\lambda}$ on each strong interval for $e<\lambda<1$.
For the lower comparison use $-A+B\log(r/r_0)$ with $aB>\sup P'$.
Now $2aB/r^2$ dominates the negative contributions from
$(aB-P')(\log\rho)_r/r$, $-Bu/r$ and $-v^2$, using \eqref{enp:negative-slope}
and $v^2\le C/r$ at contact. In each comparison choose the inner radius
$r_0$ first and then $A$ to enclose the initial and annular values.
With $\ell(r)=1+\log(R_*/r)\ge1$ on $0<r\le R(t)$, we obtain
\begin{equation}\label{enp:one-log}
 u_+\le Cr,\qquad u_-\le Cr\ell(r).
\end{equation}

Let $F(0)=0$, $F'(r)=C_1r\ell(r)\ge u_-$ and $Y=Z+F(r)$.
The expression $u^2+F'u$ is nonpositive for $u<0$ and is bounded by
$Cr^2\ell(r)$ for $u\ge0$. Thus
\[
 \dot Y+KY\le K(ru+F)+Cr^2\ell(r).
\]
Let $N(t)$ be the larger of a fixed $N_0$ and the running maximum of $Y$
on the material core. Take $N_0\ge1$ larger than the initial and exterior
values and twice $\sup|ru+F|$. At an increasing maximum of radius $r$,
\[
 0\le N'\le-KN/2+Cr^2\ell(r).
\]
This inequality holds almost everywhere by the material ordinary
differential equation. The center value is zero and hence cannot attain
an increasing maximum. Since $r^2\ell(r)\le C r/R_*$, increasing $N_0$
gives $\rho(t,r)\le(c_0/4)(r/R_*)^{1/p}$. Let $s>r$ first satisfy
$\rho(t,s)=c_0/2$. From $(\log\rho)_r\le CN/r$ we obtain
\[
 CN\log(s/r)\ge \log2+p^{-1}\log(R_*/r),\qquad
 D_\rho(t)\ge
 \frac{[c_0^{\gamma/2}-(c_0/2)^{\gamma/2}]^2}
 {\int_s^{R(t)} y^{-3}\dd y}\ge c s^2 .
\]
Since $\ell(r)\le CN\log(s/r)\le CN(s/r)^2$, it follows that
$r^2\ell(r)\le CN s^2\le CN D_\rho(t)$.
Gronwall's inequality and \eqref{enp:linear-Zenergy} give a fixed
upper bound for $Z$.

Fix $B\ge1$ with $(\log\rho)_r\le B/r$ and take $0<\eta<\min\{1/2,1/(pB)\}$.
Where $K>(r/R_*)^{pB\eta}$, \eqref{enp:one-log} gives
$u^2\le C r^2\ell(r)^2\le Cr^\eta K$.
Where $K\le(r/R_*)^{pB\eta}$ and $r$ is small, we have
$\rho(t,r)\le C(r/R_*)^{B\eta}$. Integrating $(\log\rho)_r\le B/r$ up to the
first level $\rho(t,s)=c_0/2$ gives $s\ge cr^{1-\eta}$, hence
$D_\rho(t)\ge cs^2\ge cr^{2-2\eta}$ and $u^2\le Cr^\eta D_\rho(t)$.
Away from zero the annular density band gives $K\ge c>0$. Since
$|ru|\le Cr^\eta$ throughout the core, this proves
\begin{equation}\label{enp:critical-source}
 |Kru|+u^2\le C r^\eta(K+D_\rho(t)).
\end{equation}
For $\gamma=1$, $K=\kappa/a>0$ is constant and $|u|\le C\sqrt r$.
The damped equation \eqref{eq:wz} then bounds $|Z|$ uniformly, so
$(\log\rho)_r\le B/r$ and \eqref{enp:critical-source} holds with $\eta=1$.

For either pressure law, integrating $(\log\rho)_r\le B/r$ from the outer
circle gives $\rho(t,r)\ge cr^B$. The earlier bound
$\rho(t,r)\le C_d r^{-d}$ with $d=1$ gives
$c r^{B+2}\le m(t,r)\le Cr$. Thus every earlier radius of a particle
is at most $C r(t,m)^{1/(B+2)}$. The integrating factor in
\eqref{eq:wz} has integral at most one against $K\dd s$ and at most
$\int_0^\infty D_\rho\dd s$ against $D_\rho\dd s$.
By \eqref{enp:critical-source} and the initial center bound, for some
$0<\eta_*<1$,
\begin{equation}\label{enp:center-power}
 |Z(t,r)|\le C r^{\eta_*},\qquad |(\log\rho)_r(t,r)|\le C r^{-1+\eta_*}.
\end{equation}
Integrating from the outer circle gives both
uniform density bounds, and the same estimate bounds $\nabla\log\rho$
in $L^2\cap L^q$ on the core for some $q>2$.
For every $b\ge-a$ the density now has a uniform upper bound and a
uniform spatial H\"older modulus. In the core this follows from the
weighted estimate for $\theta>0$, the Lipschitz bound for $\theta=0$,
and integration of $|(\log\rho)_r|\le Cr^{-|\theta|}$ for $-1<\theta<0$.
At $\theta=-1$ use the integrable power in \eqref{enp:center-power}.
The exterior gradient estimate gives a modulus on the overlapping region $r\ge L$.

We claim that
\begin{equation}\label{ws:band}
 0<m\le\rho(t,x)\le M,\qquad
 \lim_{t\to\infty}\|\rho(t)-\bar\rho\|_{L^\infty(\R^2)}=0 .
\end{equation}
Let $t_j\to\infty$. On each fixed ball $B_A$ we have
\[
 \int_{B_A}\bigl(\rho |U|^2+|\rho U|^2+
                  |\rho^\alpha\mathcal S(U)|^2\bigr)\dd x
 \le C_A \mathscr D(t).
\]
Indeed, $r^2\rho^{1-\alpha}\le C_A$ on $B_A$. For positive bulk viscosity,
the dissipation decomposition controls $\int\rho^\alpha(u/r)^2r\dd r$.
At $b=-a$ the core bounds from the BD comparison in Section~\ref{ws:nonpositive},
\eqref{enp:subcritical-density} and \eqref{enp:center-power}, together
with the exterior band, already give a global positive lower density bound.
Hence \eqref{enp:deviatoric-korn} controls the same angular integral.
The uniform density upper bound and the trace decomposition bound the
squared stress by $C\mathscr D(t)$ for every $b\ge-a$.
Since $\int_0^\infty\mathscr D\dd t<\infty$, the integrals of all
three terms over $[t_j,t_j+1]$ tend to zero. Extract a locally uniformly convergent subsequence
of $\rho(t_j,\cdot)$ by the spatial modulus. The continuity equation
implies that $\rho(t_j+s,\cdot)$ has the same limit uniformly for
$0\le s\le1$ on each compact set. Indeed, after spatial mollification
at a fixed scale its variation is bounded by
$C_{A,\varepsilon}(\int_{t_j}^{t_j+1}\mathscr D\dd t)^{1/2}$.
The mollification error tends to zero uniformly by the spatial modulus.
Test conservative momentum with compactly supported smooth functions
of $(s,x)$ in $(0,1)\times\R^2$. The momentum, convective flux and
stress vanish in distributions by the preceding bound. Thus the
limit density satisfies $\nabla P(\rho_*)=0$. Since $P$ is strictly
increasing on $[0,\infty)$, $\rho_*$ is constant. Passing to the
limit at each fixed $r\ge L$ in \eqref{ws:exterior} gives
$|\rho_*-\bar\rho|\le Cr^{-1/2}$. Thus $\rho_*=\bar\rho$. The same
tail bound, together with convergence on each fixed ball, gives uniform
convergence in the plane. Every sequence $t_j\to\infty$ has such a subsequence.
It gives $\rho\ge\bar\rho/2$ for all sufficiently large times, and
finite-time positivity covers the remaining interval.

\subsection{Relaxation of the relative energy}\label{ws:relaxation}

By \eqref{ws:band}, $0<m\le\rho\le M$. The nonnegative trace-free
part in \eqref{ws:energy}, followed by \eqref{enp:deviatoric-korn}, gives
\begin{equation}\label{ws:band-korn}
 \|\nabla U(t)\|_{L^2}^2\le C\mathscr D(t),\qquad
 \int_0^\infty\|U(t)\|_{L^\infty}^2\dd t\le C
 \qquad(b\ge-a).
\end{equation}
The second estimate follows from \eqref{eq:radial-sup} and the mechanical
energy identity. Write $n=\rho-\bar\rho$ and $J=(1-\Delta)^{-1}$.
With the planar measure $\dd x$, the mechanical energy is $2\pi\mathcal E$.
Set $\mathcal E_\varepsilon=2\pi\mathcal E+
\varepsilon\int_{\R^2}\rho U\cdot\nabla Jn\dd x$. For small fixed
$\varepsilon>0$, this energy is comparable to
$\|n\|_{L^2}^2+\|U\|_{L^2}^2$. We claim that
\begin{equation}\label{ws:relative-dissipation}
 \mathcal E_\varepsilon'+d_0\left(\|\nabla U\|_{L^2}^2+
       \|\nabla(1-\Delta)^{-1/2}n\|_{L^2}^2\right)\le0 .
\end{equation}
The kernel of $J=\int_0^\infty e^{-s}e^{s\Delta}\dd s$ is even,
positive and has integral one. Thus
\[
\begin{aligned}
\int_{\R^2} [P(\rho)-P(\bar\rho)](I-J)n\dd x
 &=\tfrac12\iint J(x-y)[P(\rho(x))-P(\rho(y))]
                   [n(x)-n(y)]\dd x\dd y\\
 &\ge \min_{[m,M]} P'\,\|\nabla(1-\Delta)^{-1/2}n\|_{L^2}^2 ,
\end{aligned}
\]
where $J(x-y)$ denotes the convolution kernel.
Indeed, the pressure contribution to the derivative of
$\int\rho U\cdot\nabla Jn\dd x$ is
$-\int[P(\rho)-P(\bar\rho)](I-J)n\dd x$, since $\Delta J=J-I$. The Hessian of $Jn$ is bounded in $L^2$ by
$\|\nabla(1-\Delta)^{-1/2}n\|_{L^2}$.
By \eqref{eq:radial-sup} and the kinetic energy bound, the stress and
the convective flux have $L^2$ norms at most $C\|\nabla U\|_2$.
The time derivative of the test
function contributes $\|(1-\Delta)^{-1/2}\diver(\rho U)\|_{L^2}^2$.
Decomposing $\rho U=\bar\rho U+nU$ bounds its square root by
$\bar\rho\|\nabla U\|_{L^2}+\|n\|_{L^2}\|U\|_{L^\infty}$.
Young's inequality and the mechanical identity give
\eqref{ws:relative-dissipation}.

For the low frequencies, use the acoustic system for $(n,\rho U)$
linearized at $(\bar\rho,0)$. Its Fourier evolution preserves
$P'(\bar\rho)|\widehat n|^2+|\widehat{\rho U}|^2$ at each frequency.
On $|\xi|\le h$, the divergences of the convective flux and viscous
stress have Fourier $L^2$ norm at most $Ch\|\nabla U\|_2$ by the flux
bounds above. The pressure remainder satisfies
$\|P(\rho)-P(\bar\rho)-P'(\bar\rho)n\|_1\le C\|n\|_2^2\le C$.
Its gradient contributes at most $Ch^2$, since
$\||\xi|\|_{L^2(|\xi|\le h)}=Ch^2$ in two dimensions.
Duhamel's formula and $\int_0^\infty\|\nabla U\|_2^2\dd s\le C$ give
\[
 \|(\widehat n,\widehat{\rho U})(t)\|_{L^2(|\xi|\le h)}
 \le C\|(\widehat n_0,\widehat{\rho_0U_0})\|_{L^2(|\xi|\le h)}
       +C(h\sqrt t+h^2t).
\]
Take $h=[(e+t)\log(e+t)]^{-1/2}$. Both error terms vanish, and the
initial Fourier integral over $|\xi|\le h$ tends to zero since
$(n_0,\rho_0U_0)\in L^2$. Also $\|nU\|_{L^1}\le C$, so
$\|\widehat{nU}\|_{L^2(|\xi|\le h)}\le Ch$. The identity
$\bar\rho U=\rho U-nU$ therefore gives the same low-frequency
convergence for $U$. For $|\xi|\ge h$, both dissipation multipliers
in \eqref{ws:relative-dissipation} are at least $h^2/2$. Plancherel's
identity and the energy equivalence therefore give
\[
 \mathcal E_\varepsilon'+c h^2\mathcal E_\varepsilon
 \le C h^2\int_{|\xi|\le h}\bigl(|\widehat n|^2+|\widehat U|^2\bigr)\dd\xi
 =h^2 o(1).
\]
Since $\int_0^\infty h^2\dd t=\infty$, an integrating factor gives
\begin{equation}\label{ws:zero-decay}
 \|\rho-\bar\rho\|_{L^2}+\|U\|_{L^2}\longrightarrow0 .
\end{equation}

\subsection{Higher-order convergence in the plane}\label{ws:late-high}

First suppose $\theta>-1$, including $b=-a$ when $\alpha<1$.
Choose $s,\sigma$ as in \eqref{ext:late-moment-parameters} and set
$\omega=1+r^{-\sigma}$. Testing momentum with $\omega |u|^{2s-2}u$
gives the moment and dissipation
\[
 M=\int_0^\infty \rho |u|^{2s}\omega r\dd r,\qquad
 D=\int_0^\infty \rho^\alpha|u|^{2s-2}
              (u_r^2+u^2/r^2)\omega r\dd r .
\]
Apply \eqref{ext:constant-density-weighted-coercivity} with exponent
$\sigma$ and with exponent zero. In each integration by parts the
boundary term vanishes by radial regularity at zero and by $u(t,r)\to0$
at infinity. The coefficients $2s-1$ and $1-\sigma^2/(2s)$ are
positive even at $\beta=-1$. The coercivity persists
for large $t$ since $\|\rho^\alpha-\bar\rho^\alpha\|_\infty\to0$.
Subtract $P(\bar\rho)$ before testing. Young's inequality bounds the
pressure terms by $\varepsilon D$ plus
\[
 C_\varepsilon\int_0^\infty
    |P(\rho)-P(\bar\rho)|^2(1+|u|^{2s})\omega r\dd r .
\]
The term containing $|u|^{2s}$ is $o(1)M$ by \eqref{ws:band}. For
the other term, $\omega r$ is integrable near zero since $\sigma<2$,
and $\omega$ is bounded for $r\ge1$, where \eqref{ws:zero-decay}
applies. Transport contributes at most
$\varepsilon D+C_\varepsilon\|U\|_{L^\infty}^2M$.
The coefficient $\|U\|_{L^\infty}^2$ is integrable by \eqref{ws:band-korn}.

On $r<L$ the angular part of $D$ controls the contribution to $M$,
since $r^2\le L^2$ and the density is bounded above and away from zero.
On $r\ge L$ the weight and $u$ are bounded by the exterior estimates
preceding \eqref{ws:interface}. Hence
\[
 M\le CD+C\|U\|_{L^2}^2 .
\]
Since $\|U\|_{L^2}\to0$ by \eqref{ws:zero-decay}, the term $o(1)M$
can be absorbed at large times, leaving
$M'+cM+cD\le C\|U\|_{L^\infty}^2M+o(1)$.
The unit-interval estimate following \eqref{ext:damped-gronwall} gives
$M(t)+\int_t^{t+1}D(s)\dd s\to0$.

On the material core $r<R(t)$ use the same moment as on the disk,
$Z=\int_0^{R(t)}\rho |r^{-\theta}w|^s r\dd r$. Its mass labels are
fixed, so \eqref{eq:Zidentity} holds without a boundary term.
Since $R(t)\le3L$ by \eqref{ws:interface}, the estimates leading to
\eqref{ext:late-effective-moment-decay} give
$Z'+cZ\le C(D+\sqrt D)$ with the full-plane dissipation $D$ above.
The right side has vanishing integrals on unit intervals, hence $Z\to0$.
H\"older's inequality as in \eqref{ext:late-gradient-exponent} yields
$w\to0$ in $L^2\cap L^q$ on the core for some $2<q<\min(3,s)$.
The velocity moment $M$ and \eqref{ws:zero-decay} give
$\|U\|_{L^{2q}}\to0$ by interpolation.

At $\theta=-1$, the bounded-core estimate for $rw$ following
\eqref{ext:lower-order-decay} gives the same convergence of $w$ by
\eqref{enp:center-power} and \eqref{ws:zero-decay}. The uniform velocity
bound and $\|U\|_2\to0$ give $\|U\|_{L^{2q}}\to0$ in this case as well.

For every $\theta$ the exterior is treated by one unweighted test. Set
$H=\int_{r\ge R(t)}w^2\dd m$. Equation \eqref{eq:weq}, tested by $w$
on these fixed mass labels, gives
\[
 H'+cH\le C\|U\|_\infty^2H+C\|U\|_2^2 .
\]
Here $K$ is bounded above and away from zero, $r$ is bounded below,
and the exterior bound for $u$ gives
$\int_{r\ge R(t)}u^4\dd m\le C\|U\|_2^2$.
By \eqref{ws:band-korn} and \eqref{ws:zero-decay}, $H\to0$.
Interpolation with the exterior bound for $w$ in \eqref{ws:exterior}
gives $L^q$ convergence there. Combining the two regions in \eqref{eq:w} gives
\begin{equation}\label{ws:low-grad-decay}
 \|U\|_{L^2\cap L^{2q}}
 +\|\nabla\log\rho\|_{L^2\cap L^q}\longrightarrow0 .
\end{equation}

On $[t-1,t+1]$ choose a cutoff $\eta$ that vanishes at $t-1$, equals
one on $[t,t+1]$, and satisfies $|\eta'|\le C$. Freeze the diffusion of
$\eta U$ at $\nu(\bar\rho)$ as in \eqref{ext:decaying-maximal-regularity}.
The diffusion error is absorbed by \eqref{ws:band}. In $L^q$ the
convection and stress are bounded by
\[
 C(\|U\|_{L^q}+\|\nabla\log\rho\|_{L^q})\|\nabla(\eta U)\|_\infty
 \le o(1)\|\eta U\|_{W^{2,q}} .
\]
The pressure and cutoff terms tend to zero by \eqref{ws:low-grad-decay}.
Maximal regularity with time exponents $2$ and $4$ gives
$\int_t^{t+1}(L^2+L^4)\dd s\to0$, where $L=\|U\|_{W^{2,q}}$.
We also need a separate $L^2_x$ estimate in the plane. After absorbing
the same diffusion error, the forcing is bounded in $L^2_x$ by
$C(\|U\|_2+\|\nabla\log\rho\|_2)(1+\|\nabla U\|_\infty)$.
Its $L^2_t$ norm tends to zero on $[t-1,t+1]$, so
$\int_t^{t+1}\|U_t\|_2^2\dd s\to0$.

Using \eqref{ws:zero-decay} and \eqref{ws:low-grad-decay} for the lower
norms, the elliptic identity \eqref{eq:ellipticidentity} gives
\[
\begin{aligned}
 \|U\|_{H^2}&\le C\|U_t\|_{L^2}+o(1),\\
 \|U\|_{H^3}&\le C\bigl[\|U_t\|_{H^1}+L\|U_t\|_{L^2}
 +(1+L)\|\nabla\log\rho\|_{H^1}\bigr]+o(1)(1+L).
\end{aligned}
\]
In the second estimate absorb $\|\nabla^2U\,\nabla\log\rho\|_{L^2}
\le C\|\nabla\log\rho\|_{L^q}\|U\|_{H^3}$. Test the damped equation
\eqref{ext:damped-density-gradient} at orders zero and one, and test
\eqref{eq:timeequation} by $\nu^{-1}U_t$. With a sufficiently small fixed $\delta>0$, set
\[
 Y=\delta\|\nabla\log\rho\|_{H^1}^2
       +\int_{\R^2}\nu^{-1}|U_t|^2\dd x .
\]
The second test controls $\|\nabla U_t\|_{L^2}^2$ and absorbs
the term $\delta C\|\nabla U_t\|_{L^2}^2$ from the first.
Adding a multiple of $\|U_t\|_{L^2}^2$ to both sides then gives
\[
 Y'+d_4Y+d_4\|\nabla U_t\|_{L^2}^2
 \le [C(L+L^2)+o(1)]Y+
 C\|U_t\|_{L^2}^2+o(1)(1+L^4).
\]
Here we used \eqref{ext:acceleration-L2-energy}, \eqref{ws:low-grad-decay}
and the elliptic estimate for $\|U\|_{H^2}$. By the whole-space heat estimate,
$\int_t^{t+1}\|U_t\|_2^2\dd s\to0$. The integrals of the remaining
source and of the coefficient of $Y$ over $[t,t+1]$ also tend to zero.
The unit-interval estimate following \eqref{ext:damped-gronwall}, applied to $Y$,
and the $H^2$ elliptic bound therefore yield
\[
 \|\nabla\log\rho(t)\|_{H^1}+\|U_t(t)\|_{L^2}+\|U(t)\|_{H^2}
 +\int_t^{t+1}\|\nabla U_t\|_{L^2}^2\dd s\longrightarrow0 .
\]

The same factor $\|U\|_{H^2}+\|\nabla\log\rho\|_{H^1}$ in
\eqref{eq:combined-source} now tends to zero. Test \eqref{eq:timeequation}
with $-\Delta U_t$ and use
$\|U_t\|_{H^2}\le C(\|\Delta U_t\|_2+\|U_t\|_2)$.
After absorption, $\|U_t\|_2=o(1)$ gives
\[
 \frac{\dd}{\dd t}\|\nabla U_t\|_2^2+c\|\Delta U_t\|_2^2
 \le o(1)(1+L^2+\|U_t\|_{H^1}^2).
\]
Choose $s\in[t-1,t]$ with $\|\nabla U_t(s)\|_2^2$ at most its mean on
that interval. This mean tends to zero, as does the integral of the right
side over $[t-1,t+1]$. Integrating from $s$ to $t$ and to $t+1$, and
using \eqref{eq:timeequation} for $U_{tt}$, gives
\[
 \|U_t(t)\|_{H^1}+\int_t^{t+1}
       (\|U_t\|_{H^2}^2+\|U_{tt}\|_2^2)\dd s\longrightarrow0 .
\]
The elliptic estimate now gives $\|U\|_{H^3}\le o(1)(1+L)$.
Since $H^3(\R^2)\hookrightarrow W^{2,q}(\R^2)$ for finite $q$,
we have $L\le C\|U\|_{H^3}$. For large $t$ we can absorb this term
and obtain $\|U\|_{H^3}\to0$.
Finally, the commutator estimates leading to \eqref{ext:highest-density-damping}
apply on $\R^2$. Indeed, $\nu^{-1}-\nu(\bar\rho)^{-1}$ is bounded in $H^2$,
and $\|U\|_{H^3}+\|\nabla\log\rho\|_{H^1}\to0$ allows us to absorb the
commutators into the pressure damping. Thus
\[
 \frac{\dd}{\dd t}\|\nabla\log\rho\|_{H^2}^2+c\|\nabla\log\rho\|_{H^2}^2
 \le C\|U_t\|_{H^2}^2+o(1).
\]
The right side has vanishing integrals on unit intervals. The damping
estimate following \eqref{ext:damped-gronwall} and \eqref{ws:zero-decay}
therefore give $\|\rho-\bar\rho\|_{H^3}\to0$. The continuity equation gives
$\|\rho_t\|_{H^2}\to0$. From \eqref{eq:ellipticidentity},
\[
 \|U\|_{H^4}\le C\bigl(\|U\|_2+\|U_t\|_{H^2}
       +\|\nabla\log\rho\|_{H^2}+\|U\|_{H^3}^2\bigr).
\]
Squaring and integrating gives $\int_t^{t+1}\|U\|_{H^4}^2\dd s\to0$.
Together with the finite-time bounds of Subsection~\ref{ws:finite-high},
these estimates give \eqref{ws:main-estimate} and \eqref{ws:main-limit}.

\section{Non-power viscosity laws}\label{sec:general-laws}
Under \eqref{eq:viscosity-structure}, set $\nu=A(\rho)/\rho$,
$K=P'(\rho)/\nu$ and $w=u+\nu(\log\rho)_r$. Differentiating continuity
and adding the radial momentum equation gives
\[
 \Dt w=\left[\left(1-\frac{2\rho\mu'(\rho)}{A(\rho)}\right)\frac ur-K\right](w-u)
       =(\theta u/r-K)(w-u).
\]
Thus the weighted equation \eqref{eq:wz} holds for both viscosity classes.

For $\theta>0$, write $X=u_r$, $Y=u/r$ and $\ell=\lambda/A$.
With coefficients $c_u,c_w$ for the physical and effective kinetic
energies, the weighted viscous form divided by $A$ is
\begin{equation}\label{eq:nonpower-kinetic-form}
 c_u\{X^2+2(\ell-\theta)XY+(1-2\theta\ell)Y^2\}
 -\frac{4\theta c_w}{1-\theta}\frac{\mu}{A}(XY-\theta Y^2).
\end{equation}
The last term follows by integrating
$-c_w\theta\int r^{-2\theta}\rho(w-u)u^2\dd r$ by parts, using
$\rho(w-u)=2\partial_r\mu(\rho)/(1-\theta)$.

\begin{proof}[Proof of Proposition~\ref{prop:general-laws}]
Here $A(s)$ is comparable to $s^\alpha$, and $\beta$ is given by \eqref{eq:theta}.

At $b=-a$ the monotonicity in \eqref{eq:general-law-condition}
forces $L$ to be constant by the identity following the proposition.
The power-law theorems then apply with coefficients $aL,bL$.

For $\theta=0$ and $1/2\le\alpha<1$, the same BD identity holds
because $\rho K(w-u)=P_r$. Since $A(s)\asymp s^\alpha$,
$\rho(w-u)^2\asymp\rho^{2\alpha-1}|(\log\rho)_r|^2$, so
\eqref{bdlt:uniform-logcoarse} and \eqref{bdlt:weighted-K} follow
as before. In the plane the exterior comparison values are bounded by
\eqref{ws:exterior-korn}--\eqref{enp:exterior-trace}, since
$\mu(s)\asymp s^\alpha$ and $\mu+\lambda\ge0$.

To check the velocity comparisons, write $d=sA'(s)/A(s)$ and
$k=1-\alpha$. The coefficient of $w-u$ is $d u_r+(d-1)u/r-K$.
Since $|d-\alpha|\le\varepsilon_0$, we can take $\varepsilon_0$ so
small that $1-d\ge k/2$ and
\[
 dk+d-1\le-k^2/2,\qquad (1+d)k+d-1\ge\alpha k/2.
\]
The first comparison with $-M$ follows from $(1-d)M/r\ge K$.
For $f=Br^k$ the coefficient of $w$ in the upper residual is
$K-(dk+d-1)f/r\ge0$, and the coefficient of $f^2/r$ is at least
$\alpha k/2$. Thus the negative terms are bounded by
$K(f+M)+CMf/r$, exactly the orders absorbed in the proof of
\eqref{bdlt:power-barrier}. Relaxation and integration of
$(\log\rho)_r=(w-u)/\nu$ give both density bounds and the gradient
bound, also for the weaker slope assumption at $\alpha=1/2$.
On this density interval \eqref{enp:deviatoric-korn} controls the full
velocity gradient, including when the bulk viscosity vanishes. The disk pressure test and
the plane energy argument in Section~\ref{ws:relaxation} give decay of
$U$ in $L^2$. Since $K$ is bounded above and away from zero, integration
of \eqref{bdlt:relaxation} in the transported $L^2(\dd m)$ norm gives
the corresponding decay for $w$. Interpolation with its uniform bound
gives decay in $L^q$.

The remaining density estimates concern positive reference bulk
viscosity. Integrating $|sL'/L|\le\varepsilon_0$
gives $(y/s)^{\varepsilon_0}\le L(y)/L(s)\le(y/s)^{-\varepsilon_0}$ for $0<y<s$.
Substitution in $\mu(s)=\tfrac{1-\theta}{2}\int_0^s A(y)\dd y/y$ yields
\[
 \frac{1-\theta}{2(\alpha+\varepsilon_0)}
 \le\frac{\mu(s)}{A(s)}
 \le\frac{1-\theta}{2(\alpha-\varepsilon_0)} .
\]
For $\varepsilon_0<\alpha(a+b)/(2a+b)$ the upper bound is strictly less
than one. Thus $\mu$ and $\mu+\lambda=A-\mu$ are comparable to $s^\alpha$,
and $\lambda/A=\beta+O(\varepsilon_0)$ uniformly in $s$.

\emph{Effective flux.}
The primitive in \eqref{eq:H} is now
\[
 H=\frac{2}{1-\theta}r^{1-\theta}\mu(\rho)
 -\int_r^R \rho u\,\xi^{1-\theta}\dd\xi
 +2\int_r^R\mu(\rho)\xi^{-\theta}\dd\xi .
\]
The two terms containing $r^{-\theta}\mu(\rho)$ cancel in $H_r$,
leaving $H_r=\rho w r^{1-\theta}$. The wall flux is
$A(\rho)(u_r+u/r)-P$. To remove it as in \eqref{eq:wall}, set
$G(s)=s\int_1^s2(\mu+\lambda)(y)y^{-2}\dd y$.
Then $sG'-G=2(\mu+\lambda)$. For $\alpha<1$ we have
$|G(s)|\le C(s+s^\alpha)$, and for $\alpha=1$ we have
$|G(s)|\le Cs|\log s|$. Thus mass and relative energy control
the wall correction. Since $\mu(s)$ is comparable to $s^\alpha$,
the proof of \eqref{ext:weighted-bound} applies for $\alpha>1/2$.
At $\alpha=1/2$, \eqref{eq:endpointvolterra} uses
$V=r^{1-\theta}\mu(\rho)$ and $f=r\rho u/\mu(\rho)$:
$H=2V/(1-\theta)+\int_0^x V(2-f)\dd y$ and
$\int f^2\dd x\le C\int\rho u^2r\dd r$.

On the disk, the coarse bound gives $r^2\rho^{1-\alpha}\le C$.
The angular dissipation therefore controls the kinetic energy by
\eqref{pr:iso-kinetic-dissipation}. Since $\mu$ and $\mu+\lambda$
are comparable to $\rho^\alpha$, Cauchy--Schwarz and
\eqref{ext:bounded-test-stress} bound the viscous term in the pressure
test by $C(\mathscr D\mathcal C)^{1/2}$. The continuity equation and
the pressure law give the same remaining terms, so
\eqref{pr:iso-modified-energy} holds. For $\gamma>1$ this gives the
exponential bound \eqref{ext:energy-decay}. For $\gamma=1$ the coarse
moment \eqref{pr:iso-coarse} gives the stretched exponential bound
\eqref{pr:iso-energy}. Both estimates precede the uniform density bounds.

In the plane the material endpoint in \eqref{ws:primitive} contributes
$-R^{1-\theta}\mathscr F+2R^{-\theta}\mu u$.
In the exterior region, $K$ and $\nu$ have positive lower and finite
upper bounds as long as $\bar\rho/2\le\rho\le2\bar\rho$.
The estimates \eqref{ws:mass}--\eqref{ws:interface} therefore give
the same exterior bounds and control of the stress on the material circle.

\emph{Weighted density estimates for $\theta>0$.}
The exact identities
\[
 \rho(w-u)=\frac{2}{1-\theta}\partial_r\mu(\rho),\qquad
 \rho(w-u)^2=\frac{A(\rho)^2}{\rho}|(\log\rho)_r|^2,\qquad
 \rho K(w-u)^2=P'(\rho)A(\rho)|(\log\rho)_r|^2
\]
give the density weights in the kinetic energy and pressure dissipation:
$A^2/\rho\asymp\rho^{2\alpha-1}$ and $P'A\asymp\rho^{\gamma+\alpha-1}$.
For $\alpha>1/2$ and $b>0$, the mass test of
Lemma~\ref{pos:mass-pressure} still uses $\Dt m=0$. Its velocity weight
$\rho K$ is comparable to $\rho^{\gamma+1-\alpha}$, so increasing the
threshold in that proof absorbs the comparison constants.

In \eqref{eq:nonpower-kinetic-form} the ratios $\mu/A$ and $\ell$
differ from their power-law values by
$O(\varepsilon_0)$. With the coefficients of \eqref{pht:quadratic-identity},
\eqref{eq:fix-energyidentity}, or \eqref{gh-identity}, the error is
at most $C\varepsilon_0(X^2+Y^2)$ and is absorbed by the strict
lower bound of the corresponding form. The velocity moment with
weight $r^{1+\delta}$ and exponent $j$ has normalized determinant
\[
 (j-1)(1-\ell^2)-\tfrac14[\delta-\ell(j-2)]^2 .
\]
For $b>0$ the power $J=5+b/a$ is fixed. For $b\le0$ the determinant
at $\ell=\beta$ is at least $7(1+\beta)/16$ for every $j,\delta$ in
\eqref{dic:moment-choice}, with $2<j<3$. Thus $\varepsilon_0$ may be
chosen independently of $\gamma$ in both moment estimates.
The pressure and cubic terms have the
same density powers as above. These estimates yield \eqref{pht:joint-energy}
or \eqref{dic:budget} when $\alpha>1/2$, and
\eqref{pr:half-integrable-entropy-ode} when $\alpha=1/2$. On the
material disk the boundary terms are controlled as in \eqref{ws:small-remainder}.
At $\alpha=1/2$ we have $A^2/\rho\asymp1$. Integrating the weighted
bound for $(\log\rho)_r$ from a radius where $\rho=\bar\rho$ on the disk,
or from the material circle in the plane, gives both density bounds.

For $\alpha>1/2$, the scalar velocity equation now contains
$(w-u)(d u_r+e u/r-K)$, where $d=sA'/A$ and
$e=s\lambda'/A=d-1+\theta$ are evaluated at $\rho$.
Here $d=\alpha+O(\varepsilon_0)$ and $e=\alpha\beta+O(\varepsilon_0)$
on $(0,\infty)$. If $b>0$, $Br^\theta$ is an upper barrier since
$d\theta+e$ is bounded below by a positive constant. For the lower
barrier $-Br^{-k}$ choose
\[
 \sup_{s>0}\frac{e(s)}{1+d(s)}<k<
 \inf_{s>0}\frac{e(s)}{d(s)},\qquad 0<k<1 .
\]
At the reference law the endpoints are $\alpha\beta/(1+\alpha)$
and $\beta$. The strict gap persists for small $\varepsilon_0$.
In the lower-barrier residual, the coefficient of
$B^2r^{-2k-1}$ is $e-(1+d)k<0$. The term containing $w$ is
$(e-dk)Bwr^{-k-1}$ and is controlled by the upper bound $w\le Cr^\theta$.
The diffusion term is negative. Since $K$ is bounded on each
finite interval, increasing $B$ gives the comparison on the disk
or material core. The positivity arguments following \eqref{eq:barriers}
and in Section~\ref{ws:finite-high} now use $\int_0^R r^{-k}\dd r<\infty$
and $\int_s^1 A(y)y^{-2}\dd y\to\infty$ as $s\downarrow0$.

When $b\le0$, $\theta>0$ and $\alpha>1/2$, one has $\alpha<1$.
In the enclosed mass equation \eqref{eq:fix-massequation},
$\nu=A/\rho$. Define the stationary comparison function $\Phi$ as in
\eqref{eq:fix-centralmass}, now taking $\phi_\alpha'(s)=A(s)/s^2$.
Since $\alpha<1$, $\phi_\alpha(s)\to-\infty$ as $s\downarrow0$.
The first-contact argument for \eqref{eq:fix-centralmass}
applies because the common radial derivative at contact gives
$\rho=\Phi'/r$. In the plane the endpoint comparison is
given in Section~\ref{ws:finite-high}. For the annular moment keep $Q$ from
\eqref{eq:fix-outermoment}. Its determinant $2Q-1-\beta^2Q^2$ is
strictly positive and remains so with $\lambda/A$ in place of $\beta$.
Since
\[
 (\rho^{\alpha-1})_r=(1-\alpha)\frac{\rho^\alpha}{A(\rho)}(u-w),
\]
the one-sided variation estimate \eqref{eq:fix-wlowerouter} has the
same powers of the inverse density. Its Gronwall estimate gives
positivity on each finite interval, as needed for continuation.
After continuation, the time-independent modulus \eqref{eq:fix-densitymodulus}
and mechanical energy decay give uniform positivity on the disk.
In the plane the compactness argument in \eqref{ws:band} identifies
every limiting density with the positive far-field density.

\emph{The cases $\theta\le0$.}
Here $\alpha>1/2$. For every $j>0$ one has
\[
 (\rho^j)_r=j\,\frac{\rho^\alpha}{A(\rho)}
                    \rho^{j+1-\alpha}(w-u).
\]
The factor $\rho^\alpha/A(\rho)$ lies between fixed positive constants.
This preserves both the signs and the exponents in the one-sided
estimates of the density.

First suppose $\theta=-h<0$. For $f=-Mr^{-h}$ the coefficient
$df_r+ef/r-K$ equals $(dh-e)Mr^{-h-1}-K$. Here
$dh-e=1+h-d(1-h)\ge2h$ since $d\le\alpha\le1$.
Thus the joint lower comparison and the pressure trace
\eqref{ngall:uniform-pressure-lower} follow with $j=\gamma+\alpha-1$.
In the velocity moment \eqref{ngall:p-choice}, the normalized viscous
form has determinant $1-\beta^2>0$ at $p=2$ for the power law.
It retains a fixed lower bound for $p>2$ sufficiently close to $2$
and for $\lambda/A$ sufficiently close to $\beta$. The pressure
remainder uses $P^2/A\le C\rho^{2\gamma-\alpha}$. The preceding
derivative identity in \eqref{ngall:N-D}--\eqref{ngall:uniform-upper}
therefore gives a uniform density upper bound, hence $K\le K_*$
and $\nu\ge\nu_*>0$. The exterior stress term has the same bound
on the material circle as in Section~\ref{ws:nonpositive}.

To recover the velocity bound, put $k=-\beta$. Then $h\le k<1$ and
\[
 dk+e=(1+k)(d-\alpha)\le0,\qquad
 (1+d)k+e\ge k/2,\qquad d+e\ge\alpha(1+\beta)/2>0
\]
after reducing $\varepsilon_0$. At a contact with $f=Br^k$, the
residual of the velocity equation is
\[
 \nu(1-k^2)f/r^2+[(1+d)k+e]f^2/r-Kf+
                     [K-(dk+e)f/r]w .
\]
Since $w\ge-Mr^{-h}$, its last term is bounded below by
$-CMr^{-h}(K_*+\varepsilon_0 f/r)$. The diffusion absorbs this
error near zero, since $k+h<2$ and $h<1$. On the remaining annulus
the positive quadratic term dominates for large $B$.
For the lower barrier $f=-Br^k+B_0r$, choose
$B_0\inf(d+e)\ge K_*$. Then $df_r+ef/r-K\ge0$.
The coarse bound $u,w\ge-Mr^{-h}$ and the upper barrier give
$w-u\ge-Cr^{-h}$. The potentially positive terms in the residual
are at most $CBr^{k-1-h}+C(r^{-h}+r^k)$. The negative diffusion
$-\nu_*B(1-k^2)r^{k-2}$ dominates on a fixed small disk.
Increasing $B$ puts the lower barrier below $-Mr^{-h}$ on its
complement. Together the comparisons give $|u|\le Cr^k$ with $k=-\beta$.
The argument leading to \eqref{ngall:uniform-log-gradient} now applies
on both domains with $k=-\beta$. In the pressure dissipation test use
$\rho(w-u)=2\partial_r\mu(\rho)/(1-\theta)$ and
$P'(\rho)A(\rho)\asymp\rho^{\gamma+\alpha-1}$. The comparison residual
has the same bound because $d$ and $e$ are bounded and $\nu\ge\nu_*>0$.

It remains to check the endpoint $\alpha=1$, $b=0$. Since
$\theta=0$ and $A(s)\asymp s$, the mechanical and BD
energies give $\int|(\sqrt\rho)_r|^2r\dd r\le C$ uniformly in time.
In the coordinate $x=\log(R_*/r)$ the same bound reads
$\int|(\sqrt\rho)_x|^2\dd x\le C$. Also
$\int\rho(u_x^2+u^2)\dd x\le C\mathscr D(t)$.
Here $R_*$ is an upper bound for the radius of the disk or material
core. At a point of density $s$, Cauchy--Schwarz gives $\rho\ge s/4$
on the inward interval of length $\ell=c\min(1,s)$, for fixed small
$c>0$. The one-dimensional trace inequality on this interval gives
\[
 |u|^2\le C(\ell^{-1}\!\int u^2\dd x+\ell\!\int u_x^2\dd x)
       \le \frac{C\mathscr D(t)}{s\min(1,s)} .
\]
Condition \eqref{eq:general-law-condition} thus gives
$|d(s)-1|^2|u|^2\le C\mathscr D(t)$ before any lower density bound.
In the joint comparison following \eqref{bdlt:weighted-K} take
$M(t)=M_0\exp(C_0\int_0^t\mathscr D)$ and $g=M(t)\sqrt{r/R_*}/4$.
The BD bound still gives $\sup \sqrt r\,K\le C$. Choose $M_0$ so
that $g_r\ge K$ and the initial and exterior comparisons hold.
At a velocity contact $u=\pm M+g$, one has $|u|\ge3M/4$ and
$-M\le w\le2M$. The new term in the residual has absolute value
\[
 |(d-1)(w-u)(g_r+u/r)|
 \le C|d-1|M^2/r
 \le \frac{\nu_*M}{4r^2}+C\mathscr D(t)M .
\]
The first term is absorbed by the signed angular diffusion, which
is at least $\nu_*(M-3g/4)/r^2$. The time derivative of the signed
barrier is at least $3M'/4$ and absorbs the second term when $C_0$
is large. At an effective-velocity contact the relaxation equation
points into $[-M,2M]$. Since $\int_0^\tau\mathscr D(t)\dd t\le\mathcal E(0)$
on every strong-solution interval, these comparisons bound $u,w$
independently of $\tau$. They give a uniform bound
for $|(\log\rho)_r|$ and hence both density bounds by integration
from a mean-density point on the disk or the exterior circle in the plane.

\emph{Continuation and decay.}
On each finite interval the density stays in a compact subinterval of
$(0,\infty)$ and has a uniform modulus of continuity on the material disk
or the fixed disk. For $\theta>0$, choose $q,\sigma$ by \eqref{eq:highq}.
In the physical velocity moment \eqref{eq:common-moment-form}, freeze
both $\mu$ and $\lambda$ at $\rho(t,0)$. In
\eqref{ext:constant-density-weighted-coercivity} the leading coefficient
is now $A(\rho(t,0))>0$. Smoothness on the density band controls the
oscillation of both coefficients on the cutoff support. The same
estimate for this moment, followed by the effective velocity test
\eqref{eq:Zidentity}, gives \eqref{eq:qbound}. For $\theta=0$ the
pointwise bounds suffice. For $-1<\theta=-h<0$, take
$2<q<\min\{3,2/h\}$. At $\theta=-1$ the law is linear and
\eqref{enp:finite-time-density} instead gives $2<q<\min\{3,2/(1-\eta)\}$.
In the plane, interpolate the exterior $L^2$ and pointwise estimates
as in \eqref{ws:finite-input}.

In \eqref{ext:damped-density-gradient} the coefficient of
$(\diver U)\nabla\log\rho$ is $e(\rho)$. With $\nu=A/\rho$ we have
\[
 \nabla K=(\gamma-d)K\nabla\log\rho,\qquad
 \nu_t=(d-1)\nu(\log\rho)_t .
\]
The functions $K,\nu,e$ and their derivatives are bounded on the
density band. In \eqref{ext:forcing-time-derivative}, the time derivative
of the stress coefficient proportional to $\nu$ still cancels.
Differentiating the remaining coefficient adds
$\nu\rho e'(\rho)(\log\rho)_t(\diver U)\nabla\log\rho$.
It is estimated as the other products of two $L^q$ factors in
\eqref{ext:acceleration-L2-energy}, with an extra factor
$\|\nabla U\|_\infty$. After the $H^2$ bound, use $L^4\cdot L^4$
for the two logarithmic derivatives to recover \eqref{eq:combined-source}.
The pressure derivative uses $\gamma-d$ in place of $\gamma-\alpha$.
At large times freeze the weighted forms at
$\bar\rho$ as in \eqref{ext:constant-density-weighted-coercivity}.
The plane flux estimate in Section~\ref{ws:relaxation} uses only
boundedness of the viscosities on the density band.

On the positive density interval, relative energy controls
$\|\rho-\bar\rho\|_2$. The preceding bound for $\|\nabla g\|_{H^1}$
controls the first two derivatives of $\rho$ by the chain rule and
$H^1\hookrightarrow L^4$. Smooth composition therefore bounds
$e(\rho)-e(\bar\rho)$ in $H^2$. The $H^2$ product estimate gives
\[
 \|e(\rho)(\diver U)\nabla g\|_{H^2}
 \le C\|U\|_{H^3}\|\nabla g\|_{H^2}.
\]
The factor $\|U\|_{H^3}$ tends to zero before the final damping estimate.
This proves \eqref{eq:maindecay} on the disk and
\eqref{ws:main-estimate}--\eqref{ws:main-limit} in the plane.
\end{proof}

\begin{proof}[Proof of Corollary~\ref{cor:mixed-laws}]
Write $p=s\mu'(s)/\mu(s)$. By \eqref{eq:mixed-law-condition},
$4/3\le A/\mu=2p/(1-\theta)\le10$ and $\mu(s)\le C(1+s)$.
Thus $\mu+\lambda=A-\mu\ge\mu/3$, so mechanical dissipation controls
both radial velocity gradients with weight $\mu(\rho)$.

In the flux argument take $V=r^{1-\theta}\mu(\rho)$ and
$f=r\rho u/\mu(\rho)$ as in \eqref{eq:endpointvolterra}.
Since $\rho/\mu(\rho)^2\le C$, kinetic energy bounds $\int f^2\dd x$
in the coordinate $x=\log(R/r)$, and the kinetic tail is at most
$C\|V\|_\infty\mathscr D$. The wall correction is bounded by mass and
relative energy, since its integrand grows at most as $C(1+s+s\log_+s)$.
As $P(s)\ge c_1\mu(s)-C_1$, the proof of Lemma~\ref{ext:uniform-weighted}
on the disk and of \eqref{ws:coarse} on the material core gives
\[
 r^{1-\theta}\mu(\rho)\le C,\qquad \rho\le Cr^{-2(1-\theta)} .
\]
On the disk $r^2\rho/\mu(\rho)\le C$ controls kinetic energy by
dissipation. The pressure test \eqref{pr:iso-modified-energy} applies since
$G(s)^4\le C\mu(s)$ for $G(s)=s/(s+\bar\rho)$ and $\mu(s)\le C(1+s)$.
Also $\int_0^R\rho^{3/2}r\dd r\le C$ because $3(1-\theta)<2$.
We thus have \eqref{ext:energy-decay} for $\gamma>1$ and the isothermal
bounds \eqref{pr:iso-energy}--\eqref{pr:iso-fractional-coefficients}.

Use the weights $(c_u,c_w)=(3,1)$ and write $R_c$ for $R$ on the disk
and $R(t)$ on the material core. The corresponding energy is
\[
 E=\int_0^{R_c} r^{1-2\theta}
       \left[\tfrac{\rho}{2}(3u^2+w^2)+4\Pi(\rho)\right]\dd r+C_0 .
\]
Take $C_0=0$ when $\gamma>1$. When $\gamma=1$, choose a fixed $C_0$
so the potential integral plus $C_0$ is nonnegative for all $R_c$.
In \eqref{eq:nonpower-kinetic-form} let $Q$
be the matrix of the form multiplied by $p=(1-\theta)A/(2\mu)$.
Substituting $\ell=1-2\mu/A$ gives $Q_{11}=3p$ and
\[
 \det Q=-9\theta^2p^2+6(3-2\theta)p-(3-2\theta)^2\ge\frac{131}{625}.
\]
The polynomial is concave in each variable separately, so its minimum
on the rectangle \eqref{eq:mixed-law-condition} occurs at a vertex,
with the smallest value at $(\theta,p)=(2/5,2/5)$.
Thus the viscous form controls $\mu(\rho)(X^2+Y^2)$ uniformly.
Moreover,
\[
 \frac{A(s)^2}{s}\ge c_1,\qquad P'(s)A(s)\ge c_2s^{\gamma-1/2}.
\]
Since $\rho(w-u)^2=(A^2/\rho)|(\log\rho)_r|^2$, this energy controls
$F+M$ from \eqref{gh-energy}. With weights $(3,1)$ the cubic and pressure
terms on the right of \eqref{gh-identity} have coefficients $-4\theta$
and $-8\theta$, respectively. The estimates \eqref{gh-traces}--\eqref{pr:half-pressure}
therefore apply with different constants and with $\mu(\rho)$ in place of
$\sqrt\rho$ in $D_\omega$. On the disk they give the differential inequality
\eqref{pr:half-integrable-entropy-ode} for the energy above and hence a
uniform bound on the spatial oscillation of $\log\rho$.

On the material core put
$S_P=\int_0^{R(t)}r^{1-2\theta}P'(\rho)A(\rho)|(\log\rho)_r|^2\dd r$.
With $j=(2\gamma-1)/4>0$, weighted Cauchy--Schwarz gives
\[
 \sup_{r\le R(t)}|\rho(t,r)^j-\rho(t,R(t))^j|\le C\sqrt{S_P}.
\]
Since the exterior density has positive upper and lower bounds,
small $S_P$ puts the entire core in a fixed positive density interval.
There $A(\rho)^2/\rho\le C P'(\rho)A(\rho)$, and the angular part of
$D_\omega$ controls the weighted kinetic energy of $u$. Hence
$E\le C(1+D_\omega+S_P)$ at these times. The argument following
\eqref{ws:small-remainder}, using the boundary estimates on the circle
\eqref{ws:interface}, gives the same bound on the oscillation of $\log\rho$.
On the disk mass conservation fixes $\rho=\bar\rho$ at some radius,
and in the plane \eqref{ws:exterior} bounds $\rho(t,R(t))$.
The weighted gradient bound in \eqref{gh-energy} gives
$|\log\rho(t,r)-\log\rho(t,s)|^2\le C|r^{2\theta}-s^{2\theta}|$
on either core. On this positive density interval the viscosities are
smooth, so the continuation and decay estimates in the proof of
Proposition~\ref{prop:general-laws} apply.
\end{proof}

\appendix

\section{Local existence and uniqueness}\label{app:localconstruction}
\subsection{Local solvability}
We construct the solution by alternating a transport equation for
$\rho$ with a linear parabolic equation for $U$. We first give the
estimate needed for the parabolic step.

\begin{lemma}[Linear diffusion with a variable coefficient]\label{lem:localdiffusion}
On the disk and on $[0,T_0]$, suppose that $\nu\in CH^3$, $\nu_t\in CH^2$,
$\nu\ge\nu_*>0$, $f\in CH^1\cap L^2H^2$, and $f_t\in L^2L^2$.
If $W_0\in H^3\cap H^1_0$ and
$\nu(0)\Delta W_0+f(0)\in H^1_0$, then the Dirichlet problem
$W_t-\nu\Delta W=f$, $W(0)=W_0$, has a unique solution with
\[
W\in C(H^3\cap H^1_0)\cap L^2H^4,\qquad
W_t\in CH^1_0\cap L^2H^2,\qquad W_{tt}\in L^2L^2.
\]
On $\R^2$ replace the assumption $\nu\in CH^3$ by
$\nu-\nu_\infty\in CH^3$ for a constant $\nu_\infty>0$, and replace
$H^1_0$ by $H^1$. The same conclusion holds for the Cauchy problem.
\end{lemma}

\begin{proof}
Solve for $V=W_t$ with initial value
$V(0)=\nu(0)\Delta W_0+f(0)$ and zero boundary data on the disk.
The differentiated equation is
\[
V_t-\nu\Delta V-(\nu_t/\nu)V
=f_t-(\nu_t/\nu)f .
\]
Divide by $\nu$ before Galerkin approximation and test with $V_t$.
Adding $\tfrac12(\|V\|_{L^2}^2)'=(V,V_t)_{L^2}$ and using Young's
inequality gives on either domain
\[
 \frac{\dd}{\dd t}\|V\|_{H^1}^2+c\|V_t\|_{L^2}^2
 \le C(1+\|\nu_t\|_{L^\infty}^2)\|V\|_{L^2}^2+
 C\bigl(\|f_t\|_{L^2}^2+\|\nu_t\|_{L^\infty}^2\|f\|_{L^2}^2\bigr).
\]
Gronwall bounds $V$ in $L^\infty H^1$ and $V_t$ in $L^2L^2$.
In the weak limit the equation gives $\Delta V\in L^2L^2$. Elliptic
regularity and the time-trace theorem give $V\in CH^1\cap L^2H^2$,
with zero trace on the disk.
Set $W=W_0+\int_0^tV\dd s$.
The equation for $V$ gives $\partial_t[(V-f)/\nu]=\Delta V$.
Integrating from $0$ to $t$ yields $(V-f)/\nu=\Delta W$.
Elliptic regularity now gives $W\in CH^3\cap L^2H^4$. In the plane
the additional $L^2$ norm is controlled by
$\|W_0\|_{L^2}+\int_0^t\|V\|_{L^2}\dd s$.
Uniqueness follows from the energy estimate.
\end{proof}

On the disk we start the iteration with the solution of
\[
U^{(0)}_t-\nu(\rho_0)\Delta U^{(0)}
=F_\alpha(\rho_0,U_0),\qquad U^{(0)}(0)=U_0,\qquad
U^{(0)}|_{\partial\Om}=0 .
\]
In the plane we use the corresponding Cauchy problem and take the
inhomogeneous density norms of $\rho-\bar\rho$ and $\log(\rho/\bar\rho)$.
Lemma~\ref{lem:localdiffusion} applies by \eqref{eq:compat} on the
disk and $V_0\in H^1$ in the plane.
Given $U^{(n)}$, solve
\[
\begin{cases}
\rho^{(n+1)}_t+\diver(\rho^{(n+1)}U^{(n)})=0,\\
U^{(n+1)}_t-\nu(\rho^{(n+1)})\Delta U^{(n+1)}
=F_\alpha(\rho^{(n+1)},U^{(n)}),\\
(\rho^{(n+1)},U^{(n+1)})|_{t=0}=(\rho_0,U_0),\qquad
U^{(n+1)}|_{\partial\Om}=0 .
\end{cases}
\]
At $t=0$ each parabolic equation gives $U^{(n+1)}_t(0)=V_0$.
In particular all iterates satisfy the wall compatibility condition on
the disk. Suppose the sum of the velocity norms
in \eqref{eq:class} for $U^{(n)}$ is at most $M$ on $[0,T]$, $T\le1$.
Time integration gives
\[
 \int_0^T\|U^{(n)}\|_{H^4}\dd t\le\sqrt T M,\qquad
 \sup_{t\le T}\|U^{(n)}(t)-U_0\|_{H^2}\le\sqrt T M.
\]
If $\sqrt T M$ is small, the transport estimate and the density formula
along trajectories give $\rho^{(n+1)}\ge m_0/2$ and a data-dependent
$CH^3$ bound for the density (after subtraction of $\bar\rho$ on the
plane). The same choice bounds $\|U^{(n)}\|_{CH^2}$ by a data-dependent
constant. Put $f=F_\alpha(\rho^{(n+1)},U^{(n)})$ and
$\nu=\nu(\rho^{(n+1)})$. The continuity equation gives
$\|\rho^{(n+1)}_t\|_{CH^2}\le C(1+M)$. In $f_t$ the term with
$\nabla(\log\rho)_t$ is bounded in $L^2$ by
$C(1+\|\nabla U^{(n)}\|_{L^4})\|(\log\rho)_t\|_{H^2}\le C(1+M)$.
The other products satisfy the same bound by
$H^2\hookrightarrow L^\infty$ and $H^1\hookrightarrow L^4$. Hence
\[
 \|f\|_{CH^1}\le C,\qquad
 \|f\|_{L^2H^2}+\|f_t\|_{L^2L^2}+\|\nu_t\|_{L^2L^\infty}
 \le C\sqrt T(1+M).
\]
In the equation for $V=U^{(n+1)}_t$, the zeroth-order coefficient
contributes at most $C(1+\|\nu_t\|_\infty^2)\|V\|_{H^1}^2$ after
Young's inequality. Gronwall's inequality and the elliptic estimate
of Lemma~\ref{lem:localdiffusion} bound the next sum of velocity norms
by $C+C\sqrt T(1+M)$ when $\sqrt T(1+M)$ is sufficiently small.
To choose $M$ uniformly, fix a lower bound $m_0>0$ for the initial
density and an upper bound for the initial $H^3$ norms. Write
$\nu_0=\nu(\rho_0)$ and $f_0=F_\alpha(\rho_0,U_0)$.
Composition and product estimates bound $\nu_0$ and $\nu_0^{-1}$ in
$H^3$ (after subtracting their far-field values in the plane),
$f_0$ in $H^2$, and $V_0=\nu_0\Delta U_0+f_0$ in $H^1$.
For the frozen iterate, $V=U^{(0)}_t$ solves $V_t-\nu_0\Delta V=0$
with initial value $V_0$. The energy estimates of
Lemma~\ref{lem:localdiffusion} and the identity
$\Delta U^{(0)}=(V-f_0)/\nu_0$ therefore bound its velocity norms
on $[0,1]$ uniformly for these data. Choose $M$ larger than this bound
and $2C$, then $T_0>0$ so that the preceding smallness conditions hold
and $C+C\sqrt{T_0}(1+M)\le M$.

For the contraction, take $n\ge2$ and set
$r^{(n)}=\rho^{(n)}-\rho^{(n-1)}$ and
$v^{(n)}=U^{(n)}-U^{(n-1)}$.
The differences have zero initial data. Their transport equation gives
\[
\sup_{t\le T}\|r^{(n+1)}(t)\|_{L^2}^2
\le CT\int_0^T\|v^{(n)}\|_{H^1}^2\dd t .
\]
We estimate the forcing in $H^{-1}=(H^1_0)^*$. In two dimensions
$H^2$ functions act as multipliers on $H^1_0$ and hence on $H^{-1}$,
as do constants in the plane. For either viscosity class, the
density-gradient terms in momentum after division by $\rho$ are
\[
 2DU\,\frac{\mu'(\rho)}{\rho}\nabla\rho
 +(\diver U)\frac{\lambda'(\rho)}{\rho}\nabla\rho-\nabla\Pi'(\rho).
\]
Here $s\Pi''(s)=P'(s)$. The primitives of $\mu'(s)/s$ and
$\lambda'(s)/s$, as well as $\Pi'(s)$, are Lipschitz on the common
density interval. Their differences are bounded by $C\|r^{(n+1)}\|_2$ in
$L^2$. Their gradients satisfy the same bound in $H^{-1}$, which persists
after multiplication by the strain. The same multiplier estimate bounds the
velocity differences in stress and convection by $C\|v^{(n)}\|_2$ in $H^{-1}$.
The remaining diffusion error satisfies
\[
 \|(\nu_1-\nu_2)\Delta U^{(n)}\|_{H^{-1}}
 \le C\|r^{(n+1)}\|_2\|\Delta U^{(n)}\|_{L^4},
\]
by $H^1\hookrightarrow L^4$. In the test by $v^{(n+1)}$, Young's inequality
bounds the term from $\nabla\nu$ by a small multiple of
$\|\nabla v^{(n+1)}\|_2^2$ plus $C\|v^{(n+1)}\|_2^2$.
The positive diffusion absorbs the gradient term.
The uniform iterate bounds therefore imply
\[
\sup_{t\le T}\|v^{(n+1)}(t)\|_{L^2}^2
+\int_0^T\|\nabla v^{(n+1)}\|_{L^2}^2\dd t
\le C\int_0^T
\bigl(\|r^{(n+1)}\|_{L^2}^2+\|v^{(n)}\|_{L^2}^2\bigr)\dd t .
\]
Set $X_n=\|v^{(n)}\|_{CL^2}^2+\|v^{(n)}\|_{L^2H^1}^2$.
The transport estimate bounds $\|r^{(n+1)}\|_{CL^2}^2$ by $CT X_n$.
The velocity estimate then gives $X_{n+1}\le CT X_n$ for $T\le1$,
after increasing $C$. Choosing $CT_0<1$ makes the differences
summable in $CL^2$ for the density and in $CL^2\cap L^2H^1$ for the
velocity.
Interpolation with the uniform $H^3$ bounds gives convergence in
$CH^2$ and hence a distributional solution of \eqref{eq:system}.
The weak limits retain $U\in L^2H^4$, $U_t\in L^2(H^2\cap H^1_0)$
and $U_{tt}\in L^2L^2$, so $U_t\in CH^1_0$ by the time-trace theorem.
Let $g=\log\rho$ on the disk and $g=\log(\rho/\bar\rho)$ in the plane.
For $|\beta|=3$, the $L^\infty H^3$ bound and $g\in CH^2$ imply
$z(t)=\partial^\beta g(t)\rightharpoonup\partial^\beta g_0$ in $L^2$
as $t\downarrow0$: testing against smooth compactly supported functions
identifies every weak limit. The transport source satisfies
$\|f_\beta\|_2\le C\|U\|_{H^4}(1+\|g\|_{H^3})\in L^1(0,T_0)$.
In the flow formula for \eqref{eq:transport-energy} the source integral
is $o(1)$ in $L^2$ as $t\downarrow0$. The transported initial value is
strongly continuous, so $z\in CL^2$ and $g\in CH^3$. Now
$F_\alpha(\rho,U)\in CH^1$ by the product estimates above, and
$\Delta U=\nu^{-1}(U_t-F_\alpha)\in CH^1$. Elliptic regularity
gives $U\in CH^3$, completing \eqref{eq:class}.
The equations at each step are invariant under all orthogonal
transformations. Uniqueness at that step preserves the symmetry of
the data. In particular, reflection across the radial line through
$x\ne0$ fixes $U(t,x)$, so its tangential component vanishes.

For non-power viscosities use $\nu(s)=(2\mu(s)+\lambda(s))/s$ and
the corresponding $U_t(0)$. On the common density band $\nu$ is smooth
and bounded away from zero, so the same parabolic estimates apply.

\subsection{Uniqueness in the positive strong class}\label{app:uniqueness}
Let $(\rho_i,U_i)$, $i=1,2$, be positive strong solutions on the disk
or the plane. Put $W=U_1-U_2$ and $q=\rho_1-\rho_2$.
On each finite interval the densities lie in a compact subset of
$(0,\infty)$. The continuity equations give
\[
 q_t+\diver(qU_1)=-\diver(\rho_2W).
\]
Testing with $q$ and using the $H^3$ bounds gives, for any
$\varepsilon>0$,
\[
 \frac{\dd}{\dd t}\|q\|_{L^2}^2
 \le\varepsilon\|\nabla W\|_{L^2}^2+
 C_{T,\varepsilon}(\|q\|_{L^2}^2+\|W\|_{L^2}^2).
\]

Test the difference of the momentum equations in \eqref{eq:system}
with $W$. The continuity equation for $\rho_1$ gives the kinetic
term $\frac12\frac{\dd}{\dd t}\int_\Omega\rho_1|W|^2\dd x$. The trace decomposition, $\mu+\lambda\ge0$ and the positive lower
bound for $\mu$ on this finite-time density band give, by
\eqref{enp:deviatoric-korn},
\[
 \int_\Omega \bigl[2\mu(\rho_1)|DW|^2+
        \lambda(\rho_1)(\diver W)^2\bigr]\dd x
 \ge c\|\nabla W\|_{L^2}^2 .
\]
The pressure and viscosity differences are bounded by $C_T|q|$.
After integration by parts their contributions are bounded by
$C_T\|q\|_{L^2}\|\nabla W\|_{L^2}$, since $\nabla U_2\in L^\infty$.
The two remaining terms come from inertia. Their sum in absolute value
is at most
\[
 C_T\|W\|_{L^2}^2+
 \|q\|_{L^2}\|(\partial_t+U_2\cdot\nabla)U_2\|_{L^4}\|W\|_{L^4}.
\]
The acceleration is bounded in $H^1$, and $H^1\hookrightarrow L^4$.
Young's inequality and the density estimate therefore give
\[
 \frac{\dd}{\dd t}\int_\Omega(q^2+\rho_1|W|^2)\dd x
 +c'\|\nabla W\|_{L^2}^2
 \le C_T\int_\Omega(q^2+\rho_1|W|^2)\dd x .
\]
The wall trace or a spatial cutoff in the plane removes the boundary
terms. Gronwall's inequality bounds the difference energy by
$e^{C_Tt}$ times its initial value. Equal data give uniqueness, including
among nonradial strong solutions.

\section*{Statements and Declarations}
\noindent\textbf{Funding.} This research received no external funding.

\smallskip\noindent\textbf{Competing interests.}
The author has no relevant financial or non-financial interests to disclose.

\smallskip\noindent\textbf{Data availability.}
No datasets were generated or analyzed in this study.

\smallskip\noindent\textbf{Use of generative AI.}
ChatGPT (OpenAI) assisted with drafting and revising parts of the proofs
and with editing the manuscript. Claude (Anthropic) helped check selected
calculations and revise the wording.

\end{document}